\documentclass{tm-amsart-latex}

\begin{document}

\title{
  Around Segal conjecture in $p$-adic geometry
}
  \tmnote{{\citetexmacs{TeXmacs:website}}}

\author{Zhouhang Mao}

\date{\tmdate}

\maketitle

\begin{abstract}
  This article records multiple results coming from interplay between
  de-completed topological periodic cyclic homology, Segal conjecture, and
  $F$-smoothness.
  
  We establish completeness of motivic filtration on de-completed topological
  periodic cyclic homology of commutative rings with weakly finitely generated
  absolute cotangent complex. When the ring in question is in addition
  $F$-smooth, we show that Segal conjecture holds for its topological
  Hochschild homology. We also identify our de-completed topological periodic
  cyclic homology with Manam's Frobenius untwisted topological periodic cyclic
  homology for quasiregular semiperfectoid rings.
  
  We find a crystalline degeneration of Segal conjecture which corresponds to
  such a statement for $F$-smoothness. On the other hand, inspired by
  constructions for topological Hochschild homology, the theory of cyclotomic
  synthetic spectra allows us to produce a relative conjugate filtration on
  Hodge--Tate cohomology and its variants, and in the same time, a relative
  conjugate filtration on topological Hochschild homology and its variants. As
  a consequence, we deduce transitivity of weak and strong $F$-smoothness.
\end{abstract}

\section{Introduction}

This article stems from studies of de-completed topological periodic cyclic
homology in {\cite{Mao2024a}}, which turns out to be closely related to Segal
conjecture. The interplay between de-completed topological periodic cyclic
homology, Segal conjecture, and $F$-smoothness leads us to multiple results on
both the commutative and noncommutative world, which is the main theme of this
article.

Roughly speaking, in {\cite{Bhatt2016}}, the authors produce a motivic
filtration on topological periodic cyclic homology, whose associated graded
pieces are given by Nygaard-completed absolute prismatic cohomology. In this
sense, we say that topological periodic cyclic homology {\tmdfn{corresponds
to}} Nygaard-completed absolute prismatic cohomology. Our work gave a
construction of a localizing invariant corresponding to non-Hodge-completed de
Rham cohomology. Via the same techniques, we extended
Devalapurkar--Hahn--Raksit--Yuan's de-completed topological Hochschild
homology to a localizing invariant as well. In particular, it gives rise to
de-completed topological periodic cyclic homology
$\tmop{TP}_{/\mathbb{Z}}^{\tmop{poly}}$ (in this article, we take tacit $(p,
v_1)$-completion for this de-completion). There are two natural questions out
of it:
\begin{enumerate}
  \item When is the motivic filtration on de-completed topological periodic
  cyclic homology complete?
  
  \item When is de-completed topological periodic cyclic homology already
  ``Nygaard-complete''?
\end{enumerate}
It turns out that the ``Nygaard-completeness'' of de-completed topological
periodic cyclic homology is closely related to Segal conjecture relative to
$\tmop{THH} (\mathbb{Z})$ (\Cref{rem:TP-poly-cpl-Segal-conj}). Our first
result addresses the first question, which says that the motivic filtration on
de-completed topological periodic cyclic homology is complete if the ring in
question satisfies a very weak finite generation condition:

\begin{theorem}[\Cref{lem:ft-ext-pow-unif-bnd-below,thm:mot-cpl-TP-poly}]
  Let $R$ be an (animated) commutative ring such that the cotangent complex
  $L_{R /\mathbb{Z}}$ is, as an $R$-module spectrum, $p$-completely (weakly)
  finitely generated. Then the motivic filtration on
  $\tmop{TP}_{/\mathbb{Z}}^{\tmop{poly}} (R)$ is complete.
\end{theorem}

This includes all (animated) commutative rings $R$ whose derived modulo $p$ is
(weakly) $F$-finite (\Cref{ex:F-fin-mot-cpl}). The concept of $F$-smoothness
was inspired by Segal conjecture. We prove a sufficient condition for Segal
conjecture for $\tmop{THH}$ of an $F$-smooth ring:

\begin{theorem}[\Cref{lem:ft-ext-pow-unif-bnd-below,thm:Segal-conj-F-sm}]
  Let $R$ be an $F$-smooth $p$-quasisyntomic ring such that the cotangent
  complex $L_{R /\mathbb{Z}}$ is, as an $R$-module spectrum, $p$-completely
  (weakly) finitely generated. Then the cyclotomic spectrum $\tmop{THH} (R)$
  satisfies $p$-complete Segal conjecture.
\end{theorem}

This incorporates qualitative versions of known results such as smooth
algebras over a perfectoid ring, or regular rings whose reduction $\tmop{mod}
p$ is $F$-finite (\Cref{ex:F-fin-F-sm-Segal-conj}). This also includes a DVR
which is not $F$-finite (\Cref{ex:non-F-fin-DVR}).

On the other hand, in {\cite{Bhatt2023a}}, it is established that, a
$p$-torsion-free $p$-quasisyntomic ring $R$ is $F$-smooth if its reduction $R
/ p$ is $F$-smooth. We show that, the argument in {\cite{Hahn2022}} implies an
analogue for Segal conjecture.

\begin{theorem}[\Cref{prop:crys-deg-Segal-conj}]
  Let $M$ be a right $\tmop{THH} (\tmop{BP} \langle n \rangle)$-module in
  cyclotomic spectra such that $M / (p, v_1, \ldots, v_n)$ is bounded below.
  Suppose that the cyclotomic spectrum
  \[ M \otimes_{\tmop{THH} (\tmop{BP} \langle n \rangle)}^{\mathbb{L}}
     \tmop{THH} (\mathbb{F}_p) \]
  satisfies Segal conjecture, then the cyclotomic spectrum $M / (p, v_1,
  \ldots, v_n)$ satisfies Segal conjecture.
\end{theorem}

Actually, in this case, its further reduction modulo $v_{n + 1}$ satisfies
Lichtenbaum--Quillen conjecture as well. For a cyclotomic spectrum,
Lichtenbaum--Quillen conjecture implies Segal conjecture, and the reverse
implication requires (weak) canonical vanishing in general. However, if the
cyclotomic spectrum in question admits a ring structure and a Bökstedt
element, in which case we call it a {\tmdfn{Bökstedt ring}}, it was shown in
{\cite{Burklund2023}} that canonical vanishing is automatic. An observation is
that, their argument also adapts to modules over a Bökstedt ring
(\Cref{lem:SCV-mod-Boekstedt}). In particular, the cyclotomic rings
$\tmop{THH} (\mathbb{F}_p)$ and $\tmop{THH} (\mathbb{Z})$ are Bökstedt rings
after reduction modulo $(p, v_1)$, and thus for a DG-category $\mathcal{C}$,
we may deduce boundedness of $\tmop{TR} (\mathcal{C}) / (p, v_1)$ when
$\tmop{THH} (\mathcal{C})$ satisfies Segal conjecture
(\Cref{cor:Segal-BP-n-bnd-TR}).

In addition to automatic canonical vanishing of (Segal-conjecture-satisfying)
modules over a Bökstedt ring, their de-completed $(-)^{t\mathbb{T}}$ has a
simple formula out of $(-)^{h\mathbb{T}}$
(\Cref{prop:decpl-Tate-as-localization}). In the special case
$\tmop{TP}_{/\mathbb{Z}}^{\tmop{poly}}$, we have

\begin{theorem}[\Cref{cor:TP-poly-base-chg-neg-TC,cor:TP-poly-perfd-localize-neg-TC,cor:TP-poly-untwist-TP}]
  Let $\mathcal{C}$ be a $\mathbb{Z}$-linear DG-category. Then we have
  \[ \tmop{TP}_{/\mathbb{Z}}^{\tmop{poly}} (\mathcal{C}) \simeq (\tmop{TC}^-
     (\mathcal{C}) \otimes_{\tmop{TC}^- (\mathbb{Z}), \varphi}^{\mathbb{L}}
     \tmop{TP} (\mathbb{Z}))_{(p, v_1)}^{\wedge}, \]
  and the same if we replace $\mathbb{Z}$ by a perfectoid ring, and in this
  case, we may also replace $v_1$ by $\tmop{Ker} (\theta)$ of Fontaine's map
  on the perfectoid ring in question. Moreover, when $R$ is a quasiregular
  semiperfectoid ring, we can identify $\tmop{TP}_{/\mathbb{Z}}^{\tmop{poly}}
  (R)$ with $\tmop{TP}^{(- 1)} (R)$ in {\cite{Manam2024}} as expected.
\end{theorem}

Moreover, combining with ideas in {\cite{Antieau2023}} on continuity of
syntomic cohomology, we give a conceptual explanation of continuity of
$\tmop{TC}$ for modules over a Bökstedt ring
(\Cref{prop:TC-Boekstedt-continuity,cor:cts-TC-p-v-cpl}).

Finally, the theory of cyclotomic syntomic spectra in {\cite{Antieau2024a}}
allows us to systematically consider variants of prismatic cohomology and
topological Hochschild homology on the same foot. In particular, we produce
relative conjugate filtration on both Hodge--Tate cohomology and
$\tmop{THH}_{\Prism}^{\Phi C_p}$ (a de-completion of $\tmop{THH}^{t C_p}$):
the former is a mixed characteristic analogue of the one in
{\cite{Bhatt2012a}}, while the later is a mixed characteristic and relative
generalization of Kaledin's conjugate filtration.

\begin{theorem}[\Cref{ex:KO-Hdg-Tate,ex:KO-THH-Phi-Cp,prop:KO-cpl-Frob-twist}]
  Let $k$ be an (animated) commutative ring.
  \begin{enumerate}
    \item Let $R$ be an (animated) commutative $k$-algebra. Then there exists
    a complete exhaustive filtration on $\barPrism_R \{ i \}$ with the $j$-th
    associated graded piece given by
    \[ \left( \bigwedgestar_R^{i - j} L_{R / k} \right) \otimes_k^{\mathbb{L}}
       \barPrism_k \{ j \} [j - i] . \]
    \item Let $\mathcal{C}$ be a dualizable presentable stable $k$-linear
    $\infty$-category. Then there exists an exhaustive filtration on
    $\tmop{THH}_{\Prism} (\mathcal{C})^{\Phi C_p}$ with the $j$-th associated
    graded piece given by
    \[ \tmop{HH} (\mathcal{C}/ k) \otimes_k^{\mathbb{L}} \barPrism_k \{ j \}
       [2 j] . \]
    Moreover, if the $k$-module spectrum $L_{k /\mathbb{Z}}$ is weakly
    finitely generated, then this filtration is complete.
  \end{enumerate}
\end{theorem}

A similar construction leads to a relative filtration on $\tmop{gr}_N^{\ast} 
\Prism$ and $\tmop{THH}$:

\begin{theorem}[\Cref{ex:KO-Nygaard-gr,ex:KO-THH,prop:KO-cpl-THH}]
  Let $k$ be an (animated) commutative ring.
  \begin{enumerate}
    \item Let $R$ be an (animated) commutative $k$-algebra. Then there exists
    a complete exhaustive filtration on $\barPrism_R \{ i \}$ with the $j$-th
    associated graded piece given by
    \[ \left( \bigwedgestar_R^{i - j} L_{R / k} \right) \otimes_k^{\mathbb{L}}
       \tmop{gr}_N^j  \Prism_k \{ j \} [j - i] . \]
    \item Let $\mathcal{C}$ be a dualizable presentable stable $k$-linear
    $\infty$-category. Then there exists an exhaustive filtration on
    $\tmop{THH} (\mathcal{C})$ with the $j$-th associated graded piece given
    by
    \[ \tmop{HH} (\mathcal{C}/ k) \otimes_k^{\mathbb{L}} \tmop{gr}_N^j 
       \Prism_k \{ j \} [2 j] . \]
    Moreover, if $\tmop{THH} (\mathcal{C})$ is $t$-bounded below, then this
    filtration is complete.
  \end{enumerate}
\end{theorem}

As a consequence, we get transitivity of weak $F$-smoothness and
$F$-smoothness, which was our original motivation to produce relative
conjugate filtrations, inspired by Tess {\tmname{Bouis}}'s proof of
transitivity of $p$-Cartier smoothness in {\cite[Lem~2.10]{Bouis2023}}. Before
that, we first establish an animated version of transitivity of $p$-Cartier
smoothness.

\begin{theorem}[\Cref{thm:transitivity-Cart-sm}]
  \label{thm:Tor-ampl-transitivity-Cart-sm}Let $k$ be an (animated) ring, and
  $R \rightarrow A$ a map of (animated) $k$-algebras, such that $R
  \otimes_{\mathbb{Z}}^{\mathbb{L}} \mathbb{F}_p$ is $t$-bounded, and $A
  \otimes_{\mathbb{Z}}^{\mathbb{L}} \mathbb{F}_p$ has bounded
  $\tmop{Tor}$-amplitude over $R \otimes_{\mathbb{Z}}^{\mathbb{L}}
  \mathbb{F}_p$. Suppose that
  \begin{itemize}
    \item the $R$-module spectrum $L_{R / k}$ is $p$-completely flat; and that
    
    \item the $A$-module spectrum $L_{A / R}$ is $p$-completely flat.
  \end{itemize}
  If the $p$-completed derived de Rham cohomology $\tmop{dR}_{A / R}$ is
  $p$-completely Hodge-complete, then the $p$-completed derived de Rham
  cohomology $\tmop{dR}_{A / k}$ is $p$-completely ``Hodge-complete relative
  to'' $\tmop{dR}_{R / k}$. In particular, if the $p$-completed derived de
  Rham cohomology $\tmop{dR}_{R / k}$ is $p$-completely Hodge-complete, then
  so is $\tmop{dR}_{A / k}$.
\end{theorem}

Now we phrase our transitivity of weak $F$-smoothness and $F$-smoothness.

\begin{theorem}[\Cref{thm:weak-F-sm-quasism}]
  Let $k \rightarrow R$ be a map of (animated) rings. Suppose that
  \begin{itemize}
    \item the (animated) ring $k$ is weakly $F$-smooth; and that
    
    \item the $R$-module spectrum $L_{R / k}$ is $p$-completely flat.
  \end{itemize}
  Then the (animated) ring $R$ is weakly $F$-smooth as well.
\end{theorem}

Inspired by relative Segal conjecture, we introduce relative
Nygaard-completeness (\Cref{def:rel-Nyg-cpl}), and prove that it is implied by
$p$-Cartier smoothness with weak $F$-smoothness of the base.

\begin{theorem}[\Cref{thm:dR-Hdg-cpl-rel-Nyg-cpl,cor:transitivity-F-sm}]
  \label{thm:Tor-ampl-transitivity-F-sm}Let $k \rightarrow R$ be a map of
  (animated) rings, such that $k \otimes_{\mathbb{Z}}^{\mathbb{L}}
  \mathbb{F}_p$ is $t$-bounded, and $R \otimes_{\mathbb{Z}}^{\mathbb{L}}
  \mathbb{F}_p$ has bounded $\tmop{Tor}$-amplitude over $k
  \otimes_{\mathbb{Z}}^{\mathbb{L}} \mathbb{F}_p$. Suppose that
  \begin{itemize}
    \item the (animated) ring $k$ is weakly $F$-smooth; and that
    
    \item the $R$-module spectrum $L_{R / k}$ is $p$-completely flat.
  \end{itemize}
  If the $p$-completed derived de Rham cohomology $\tmop{dR}_{R / k}$ is
  $p$-completely Hodge-complete, then the absolute prismatic cohomology
  $\Prism_R \{ \ast \}$ is Nygaard-complete relative to $\Prism_k \{ \ast \}$.
  In particular, if $k$ is in addition $F$-smooth, then so is $R$.
\end{theorem}

We have formulated \Cref{thm:Tor-ampl-transitivity-F-sm} to be clearly an
absolute version of \Cref{thm:Tor-ampl-transitivity-Cart-sm}, which
generalizes {\cite[Prop~4.9]{Bhatt2023a}}. There are variants of
\Cref{thm:Tor-ampl-transitivity-Cart-sm,thm:Tor-ampl-transitivity-F-sm}
without bounds on $\tmop{Tor}$-amplitude of ring maps. For this, we introduce
the concept of universal convergence of Hodge and Nygaard filtrations
(\Cref{def:Hdg-filt-univ-conv,def:Nyg-filt-univ-conv}), which generalizes
``having a universal Cartier isomorphism'' in {\cite[Def~9.5.9]{Bhatt2018a}}.
These notions are to explore, but as an exhibition of its utility, we have

\begin{theorem}[\Cref{cor:F-sm-rel-perf}]
  Let $k \rightarrow R$ be a map of rings such that the induced map $k
  \otimes_{\mathbb{Z}}^{\mathbb{L}} \mathbb{F}_p \rightarrow R
  \otimes_{\mathbb{Z}}^{\mathbb{L}} \mathbb{F}_p$ of animated
  $\mathbb{F}_p$-algebras is relatively perfect. If the ring $k$ is weakly
  $F$-smooth, the absolute prismatic cohomology $\Prism_k$ is
  Nygaard-complete, and the Nygaard filtration on $\hatbarPrism_k \{ \ast \}$
  converges universally, then the same is true for $R$.
\end{theorem}

This is a version of {\cite[Prop~4.7]{Bhatt2023a}} without flatness, giving an
answer to the paragraph right before {\tmem{loc. cit.}}

\begin{acknowledgments*}
  We thank Benjamin {\tmname{Antieau}}, Tess {\tmname{Bouis}}, Daniel
  {\tmname{Fink}}, Rok {\tmname{Gregoric}}, Hyungseop {\tmname{Kim}} and Georg
  {\tmname{Tamme}} for discussions. The author acknowledges support from the
  Deutsche Forschungsgemeinschaft (DFG, German Research Foundation) through
  the Collaborative Research Centre TRR 326 {\tmem{Geometry and Arithmetic of
  Uniformized Structures}}, project number 444845124.
\end{acknowledgments*}

\section{De-completed $\tmop{TP}$ revisited}

Devalapurkar--Hahn--Raksit--Yuan defines de-completed topological periodic
cyclic homology of animated rings in terms of genuine equivariant homotopy
theory, and {\cite[§5]{Mao2024a}} extended it to DG-categories. In this
section, we first reformulate it via an equivalent description
(\Cref{def:TP-poly}) without genuine equivariant homotopy theory, explaining
its relation to Segal conjecture
(\Cref{rem:TP-poly-cpl-Segal-conj,prop:rel-abs-segal-conj}). Then in
\Cref{cor:TP-poly-untwist-TP}, we identify it with $\tmop{TP}^{(- 1)}$
introduced in {\cite{Manam2024}} (of which {\cite[§1]{Gregoric2025}} gives an
excellent exposition) in the general context of modules over a Bökstedt ring,
a concept introduced in {\cite{Burklund2023}}. Finally, we give a conceptual
explanation (\Cref{cor:cts-TC-p-v-cpl}) of continuity of $\tmop{TC}$, which
applies to filtered systems which are not uniformly $t$-bounded below.

\begin{definition}
  \label{def:TP-poly}Let $\mathcal{C}$ be a dualizable presentable stable
  $\mathbb{Z}$-linear $\infty$-category. Then its ($(p, v_1)$-completed)
  {\tmdfn{topological polynomial periodic cyclic homology}}
  $\tmop{TP}^{\tmop{poly}}_{/\mathbb{Z}} (\mathcal{C})$ is defined to be the
  $(p, v_1)$-completion of the spectrum
  \[ \tmop{THH}_{\Prism} (\mathcal{C})^{h\mathbb{T}} \assign (\tmop{THH}
     (\mathcal{C}) \otimes_{\tmop{THH} (\mathbb{Z}), \varphi_p}^{\mathbb{L}}
     \tmop{THH} (\mathbb{Z})^{t C_p})^{h\mathbb{T}}, \]
  where $\tmop{THH}_{\Prism}$ is de-completed $\tmop{THH}$ à la
  Devalapurkar--Hahn--Raksit--Yuan, extended to DG-categories in
  {\cite[§5]{Mao2024a}}. In this paper, we only use
  \[ \tmop{THH}_{\Prism} (\mathcal{C})^{\Phi C_p} = \tmop{THH} (\mathcal{C})
     \otimes_{\tmop{THH} (\mathbb{Z})}^{\mathbb{L}} \tmop{THH} (\mathbb{Z})^{t
     C_p} \]
  to simplify the cumbersome expression on the right hand side.
\end{definition}

\begin{remark}
  \label{rem:asm-TP-poly}Let $\mathcal{C}$ be a dualizable presentable stable
  $\mathbb{Z}$-linear $\infty$-category. Then the {\tmdfn{assembly map}} out
  of $\tmop{TP}_{/\mathbb{Z}}^{\tmop{poly}} (\mathcal{C})$ is given by the
  $(p, v_1)$-completion of the map
  \[ \left( \tmop{THH}_{\Prism} (\mathcal{C})^{\Phi C_p} \right)^{h\mathbb{T}}
     \longrightarrow (\tmop{THH} (\mathcal{C})^{t C_p})^{h\mathbb{T}} \]
  induced by the cyclotomic Frobenius map $\tmop{THH} (\mathcal{C})
  \rightarrow \tmop{THH} (\mathcal{C})^{t C_p}$, where the target can be $(p,
  v_1)$-completely identified with $\tmop{TP} (\mathcal{C})$ by the Tate orbit
  lemma\footnote{In full generality, we need the $t$-boundedness of $\tmop{TR}
  (\mathbb{Z}) / (p, v_1)$, a special case of Lichtenbaum--Quillen for
  $\tmop{TR} (\tmop{BP} \langle n \rangle) / (p, v_1, \ldots, v_{n + 1})$
  established in {\cite{Hahn2022}}.}.
\end{remark}

\begin{remark}
  In \Cref{def:TP-poly,rem:asm-TP-poly}, we may replace the base $\tmop{THH}
  (\mathbb{Z})$ by $\tmop{THH} (\tmop{BP} \langle n \rangle)$, and replace
  $(p, v_1)$-completion by $(p, v_1, \ldots, v_{n + 1})$-completion.
\end{remark}

{\construction{\label{cons:mot-fil}Let $R$ be an animated ring. Then by
{\cite{Antieau2024a}} (on polynomial rings, from which one can left Kan extend
to general animated rings), the filtered cyclotomic spectrum $\tmop{Fil}_M
\tmop{THH} (R)$ acquires an $\mathbb{E}_{\infty}$-cyclotomic synthetic
spectrum structure. This gives rise to a synthetic spectrum structure on
\[ \tmop{Fil}_M \tmop{THH}_{\Prism} (R)^{\Phi C_p} \assign \tmop{Fil}_M
   \tmop{THH} (R) \otimes_{\tmop{Fil}_{\tmop{ev}} \tmop{THH} (\mathbb{Z}),
   \varphi_p}^{\mathbb{L}} \tmop{Fil}_{\tmop{ev}} \tmop{THH} (\mathbb{Z})^{t
   C_p} \]
with a $\mathbb{T}_{\tmop{ev}}$-action, where $\tmop{Fil}_M \tmop{THH}
(\mathbb{Z})^{t C_p} \simeq (\tmop{Fil}_M \tmop{THH} (\mathbb{Z}))^{t C_{p,
\tmop{ev}}}$. In particular, it gives rise to filtrations on
$\tmop{THH}_{\Prism} (R)^{\Phi C_p}$ and on
$\tmop{TP}_{/\mathbb{Z}}^{\tmop{poly}} (R)$ as well, called {\tmdfn{motivic
filtrations}}. By construction, the first filtration is exhaustive. To see
that the second filtration is also exhaustive, we need some preparations for
$\tmop{Fil}_M \tmop{THH} (\mathbb{Z})$.}}

\begin{remark}
  \label{rem:Frob-THH-Z-inv-mu1}Let $R$ be an animated ring. Recall that
  \[ \pi_{\ast} (\tmop{THH} (\mathbb{Z}) \otimes_{\mathbb{Z}}^{\mathbb{L}}
     \mathbb{F}_p) =\mathbb{F}_p [\mu_1] \otimes \Lambda (\lambda_1) \]
  with $\deg \mu_1 = 2 p$ and $\deg \lambda_1 = 2 p - 1$, and the Frobenius
  map
  \[ \tmop{THH} (\mathbb{Z}) \otimes_{\mathbb{Z}}^{\mathbb{L}} \mathbb{F}_p
     \rightarrow \tmop{THH} (\mathbb{Z})^{t C_p}
     \otimes_{\mathbb{Z}}^{\mathbb{L}} \mathbb{F}_p \]
  is inverting the class $\mu_1$. As in {\cite[§6]{Hahn2022a}}, these classes
  lifts to synthetic classes (non-equivariantly), still denoted by $\lambda_1$
  and $\mu_1$ respectively, and the Frobenius map is inverting $\mu_1$ as
  well. This realizes $\tmop{Fil}_M \left( \tmop{THH}_{\Prism} (R)^{\Phi C_p}
  \otimes_{\mathbb{Z}}^{\mathbb{L}} \mathbb{F}_p \right)$ as inverting the
  class $\mu_1$ in $\tmop{Fil}_M (\tmop{THH} (R)
  \otimes_{\mathbb{Z}}^{\mathbb{L}} \mathbb{F}_p)$, and in particular,
  realizes $\tmop{THH}_{\Prism} (R)^{\Phi C_p}
  \otimes_{\mathbb{Z}}^{\mathbb{L}} \mathbb{F}_p$ as inverting the class
  $\mu_1$ in $\tmop{THH} (R) \otimes_{\mathbb{Z}}^{\mathbb{L}} \mathbb{F}_p$.
\end{remark}

\begin{lemma}
  \label{lem:frob-twist-htpy-fixed-pts-pres-colim}Let $M$ be a $\tmop{Fil}_M
  \tmop{THH} (\mathbb{Z})$-module in $\tmop{SynSp}_{\mathbb{T}_{\tmop{ev}}}$,
  and
  \[ N \assign M \otimes_{\tmop{Fil}_M \tmop{THH} (\mathbb{Z})}^{\mathbb{L}}
     \tmop{Fil}_M \tmop{THH} (\mathbb{Z})^{t C_p} . \]
  Then the norm map
  \[ \Sigma \mathbb{S}_{\tmop{ev}} (1) \otimes_{\mathbb{S}_{\tmop{ev}}}
     N_{h\mathbb{T}_{\tmop{ev}}} \longrightarrow N^{h\mathbb{T}_{\tmop{ev}}}
  \]
  is an equivalence of synthetic spectra after $(v_0,
  \tilde{v}_1)$-completion. Consequently, the functor
  \begin{eqnarray*}
    \tmop{Mod}_{\tmop{Fil}_M \tmop{THH} (\mathbb{Z})}
    (\tmop{SynSp}_{\mathbb{T}_{\tmop{ev}}}) & \longrightarrow &
    \tmop{SynSp}_{(v_0, \tilde{v}_1)}^{\wedge}\\
    M & \longmapsto & (N^{h\mathbb{T}_{\tmop{ev}}})_{(v_0,
    \tilde{v}_1)}^{\wedge}
  \end{eqnarray*}
  preserves small colimits, where $v_0$ and $\tilde{v}_1$ are the synthetic
  analogues of the homotopy classes $p$ and $v_1$ in $\mathbb{S}$ and
  $\mathbb{S}/ p$.
\end{lemma}

\begin{proof}
  By {\cite[Cons~2.58]{Antieau2024}}, we have a fiber sequence
  \[ \Sigma \mathbb{S}_{\tmop{ev}} (1) \otimes_{\mathbb{S}_{\tmop{ev}}}
     N_{h\mathbb{T}_{\tmop{ev}}} \longrightarrow N^{h\mathbb{T}_{\tmop{ev}}}
     \longrightarrow N^{t\mathbb{T}_{\tmop{ev}}} . \]
  The result follows, as the synthetic spectrum $N^{t\mathbb{T}_{\tmop{ev}}}$
  carries a $(\tmop{Fil}_M \tmop{THH} (\mathbb{Z})^{t
  C_p})^{t\mathbb{T}_{\tmop{ev}}}$-module structure, and the
  $\mathbb{E}_{\infty}$-algebra $((\tmop{Fil}_M \tmop{THH} (\mathbb{Z})^{t
  C_p})^{t\mathbb{T}_{\tmop{ev}}})_{(v_0, \tilde{v}_1)}^{\wedge}$ in synthetic
  spectra is zero.
  
  A way to see the later: the $\mathbb{T}_{\tmop{ev}}$-equivariant lift
  $\tilde{\mu}_1$ of the class $\mu_1$ in $\tmop{Fil}_M \tmop{THH}
  (\mathbb{Z}) / (v_0, \tilde{v}_1)$ is send to an invertible
  $\mathbb{T}_{\tmop{ev}}$-equivariant class in $\tmop{Fil}_M \tmop{THH}
  (\mathbb{Z})^{t C_p} / (v_0, \tilde{v}_1)$, thus to an invertible class in
  \[ (\tmop{Fil}_M \tmop{THH} (\mathbb{Z})^{t C_p})^{t\mathbb{T}_{\tmop{ev}}}
     / (v_0, \tilde{v}_1) . \]
  On the other hand, it becomes a zero class in $(\tmop{Fil}_M \tmop{THH}
  (\mathbb{Z}))^{t\mathbb{T}_{\tmop{ev}}}$, which sends to the zero class
  along the map
  \[ (\tmop{Fil}_M \tmop{THH} (\mathbb{Z}))^{t\mathbb{T}_{\tmop{ev}}} / (v_0,
     \tilde{v}_1) \rightarrow (\tmop{Fil}_M \tmop{THH} (\mathbb{Z})^{t
     C_p})^{t\mathbb{T}_{\tmop{ev}}} / (v_0, \tilde{v}_1) . \]
\end{proof}

\begin{corollary}
  \label{cor:exhaustive-mot-fil-TP-poly}Let $R$ be an animated ring. Then the
  motivic filtration on $\tmop{TP}_{/\mathbb{Z}}^{\tmop{poly}} (R)$ is
  exhaustive.
\end{corollary}

\begin{remark}
  Let $R$ be a polynomial ring. Then the map
  \[ \tmop{THH}_{\Prism} (R)^{\Phi C_p} \longrightarrow \tmop{THH} (R)^{t C_p}
  \]
  induced by the cyclotomic Frobenius $\tmop{THH} (R) \rightarrow \tmop{THH}
  (R)^{t C_p}$ is an equivalence after $(p, v_1)$-completion by results in
  {\cite{Mao2024a}}, and passing to associated graded pieces of motivic
  filtrations, it is an equivalence as well by, say,
  {\cite[Prop~2.6]{Bhatt2023a}} and $F$-smoothness of $R$. It follows that the
  motivic filtrations considered in Construction~\ref{cons:mot-fil} coincides
  with the one considered in {\cite{Mao2024a}}.
\end{remark}

We may replace the $(p, v_1)$-completed base $\tmop{THH} (\mathbb{Z})$ (or
more precisely, its reduction modulo $(p, v_1)$) by a Bökstedt ring, a notion
that appears in {\cite{Burklund2023}}.

\begin{definition}
  \label{def:Boekstedt-ring}A {\tmdfn{Bökstedt ring of (height $n \in
  \mathbb{N}$)}} is given by the datum $(R, \mu)$ where
  \begin{enumerate}
    \item $R$ is a $p$-nilpotent $t$-bounded-below $h\mathbb{A}_2$-ring in
    cyclotomic spectra, which satisfies Segal conjecture;
    
    \item $\mu$ is an element of $\pi_{2 p^n} (R)$ which is a Bökstedt class
    {\cite[Def~2.31]{Burklund2023}}.
  \end{enumerate}
\end{definition}

\begin{remark}
  In practice, Bökstedt rings come from modulo out a sequence of homotopy
  classes from an $\mathbb{E}_2$-algebra in cyclotomic spectra, such as
  $\tmop{THH} (\tmop{BP} \langle n \rangle) / (p, v_1, \ldots, v_{n + 1})$.
  Thanks to {\cite{Burklund2022}}, we may raise powers of homotopy classes to
  improve associativity and commutativity. In particular, one might freely
  replace $h\mathbb{A}_2$ by $\mathbb{E}_1$ or even $\mathbb{E}_2$ in
  \Cref{def:Boekstedt-ring}.
\end{remark}

Recall that every Bökstedt ring $(R, \mu)$ has the following properties
{\cite[Lem~2.33~\&~2.36 \& Prop~2.30]{Burklund2023}}:
\begin{itemize}
  \item The spectrum $R / \mu$ is $t$-bounded.
  
  \item The map $R [\mu^{- 1}] \rightarrow R^{t C_p}$ induced by the
  cyclotomic Frobenius $\varphi_p \of R \rightarrow R^{t C_p}$ is an
  equivalence.
  
  \item The cyclotomic spectrum $R$ is cyclotomically $t$-bounded, i.e. the
  spectrum $\tmop{TR} (R)$ is $t$-bounded.
\end{itemize}
And conversely, every $p$-nilpotent $h\mathbb{A}_2$-ring in cyclotomic
spectra, which is cyclotomically $t$-bounded, admits a Bökstedt ring
structure.

\begin{example}
  $(\tmop{THH} (\mathbb{F}_p) / p, (\mu_0 \in \tmop{THH}_2 (\mathbb{F}_p)))$
  is a Bökstedt ring.
\end{example}

\begin{example}
  \label{ex:Boekstedt-THH-BP-n}$(\tmop{THH} (\mathbb{Z}) / (p, v_1), (\mu_1
  \in \tmop{THH}_{2 p} (\mathbb{Z}; \mathbb{F}_p)))$ is a Bökstedt ring, or
  more generally, $(\tmop{THH} (\tmop{BP} \langle n \rangle) / (p, v_1,
  \ldots, v_{n + 1}), (\mu_{n + 1} \in \tmop{THH}_{2 p^{n + 1}} (\tmop{BP}
  \langle n \rangle ; \mathbb{F}_p)))$. This is a consequence of
  {\cite[Thm~G]{Hahn2022}}.
\end{example}

\begin{example}[{\cite[§5]{Mao2024a}}]
  \label{ex:Boekstedt-THH-perfd}Let $S$ be a perfectoid ring. Then
  $(\tmop{THH} (S) / (p, v_1), (\mu_0 \in \pi_2 (\tmop{THH} (S)_p^{\wedge})))$
  is a Bökstedt ring, and the map $\tmop{THH} (\mathbb{Z}) / p \rightarrow
  \tmop{THH} (S) / p$ carries the class $\mu_1$ to $\mu_0^p$ up to a unit.
\end{example}

The ``Nygaard-completeness'' of $\tmop{TP}^{\tmop{poly}}$ is closely related
to the following relative version of Segal conjecture
(\Cref{rem:TP-poly-cpl-Segal-conj}).

\begin{definition}
  \label{def:rel-Segal-conj}Let $(R, \mu)$ be a Bökstedt ring. We say that a
  right $R$-module $M$ in cyclotomic spectra satisfies {\tmdfn{Segal
  conjecture relative to $R$}} if the map
  \[ M [\mu^{- 1}] = M \otimes_R R^{t C_p} \longrightarrow M^{t C_p} \]
  induced by the cyclotomic Frobenius $\varphi_p \of M \rightarrow M^{t C_p}$
  is an equivalence. Equivalently, it is saying that $\tmop{fib} (\varphi_p)
  [\mu^{- 1}] = 0$.
\end{definition}

\begin{remark}
  \label{rem:TP-poly-cpl-Segal-conj}Let $\mathcal{C}$ be a dualizable
  presentable stable $\mathbb{Z}$-linear $\infty$-category. Then by
  definition, the assembly map $\tmop{TP}^{\tmop{poly}} (\mathcal{C})
  \rightarrow \tmop{TP} (\mathcal{C})_{(p, v_1)}^{\wedge}$ is an equivalence
  if the cyclotomic spectrum $\tmop{THH} (\mathcal{C}) / (p, v_1)$ satisfies
  Segal conjecture relative to $\tmop{THH} (\mathbb{Z}) / (p, v_1)$.
\end{remark}

The following statement compares absolute and relative Segal conjecture, which
follows from adapting the proof of {\cite[Lem~2.33]{Burklund2023}}.

\begin{proposition}
  \label{prop:rel-abs-segal-conj}Let $(R, \mu)$ be a Bökstedt ring. Then a
  right $R$-module $M$ in cyclotomic spectrum satisfies Segal conjecture if
  and only the following two conditions hold:
  \begin{enumerate}
    \item The spectrum $M / \mu$ is $t$-bounded above.
    
    \item The $R$-module $M$ in cyclotomic spectrum satisfies Segal conjecture
    relative to $R$.
  \end{enumerate}
\end{proposition}

\begin{proof}
  Let $N \assign \tmop{fib} (\varphi_p)$. The canonical map $N \rightarrow M$
  is an equivalence after modulo $\mu$, as $\mu$ acts invertibly on $M^{t
  C_p}$. Suppose that $M$ satisfies Segal conjecture, then $N$ is $t$-bounded
  above, and in this case, so is $N / \mu$. Thus it suffices to establish the
  equivalence of Segal conjecture and relative Segal conjecture under the
  $t$-bounded-aboveness of $N / \mu$. Indeed, for every $m \in \mathbb{Z}$,
  the sequential tower
  \[ \pi_m (N) \xrightarrow{\mu} \pi_{m + 2 p^n} (N) \xrightarrow{\mu} \pi_{m
     + 4 p^n} (N) \xrightarrow{\mu} \cdots, \]
  which computes $\pi_m (N [\mu^{- 1}])$, stabilizes (uniformly for $m \in
  \mathbb{N}$, say). This implies that, the spectrum $N$ is $t$-bounded above
  (i.e. Segal conjecture) if and only if $N [\mu^{- 1}] = 0$ (i.e. relative
  Segal conjecture).
\end{proof}

We note that, for modules over a Bökstedt ring, it satisfies (absolute) Segal
conjecture if and only if it satisfies it has bounded $\tmop{TR}$.

\begin{lemma}[{\cite[Lem~2.36]{Burklund2023}}]
  \label{lem:SCV-mod-Boekstedt}Let $(R, \mu)$ be a Bökstedt ring, and $M$ a
  right $R$-module in $\tmop{Sp}^{B\mathbb{T}}$. Suppose that the spectrum $M
  / \mu$ is $b$-truncated. Then $M$ satisfies strong canonical vanishing with
  parameter $b$.
\end{lemma}

\begin{proof}
  The proof of {\cite[Lem~2.36]{Burklund2023}} works verbatim. We briefly
  summarize as follows. First pick a lift $\hat{\mu} \in \pi_{2 p^n} (\Sigma
  R_{h\mathbb{T}})$, and there exists $d \gg 0$ such that $a_{\lambda}^d 
  \hat{\mu} = 0$. Then the composite $\mathbb{T}$-equivariant map
  \begin{equation}
    \tau_{> b} M \xrightarrow{c} M \xrightarrow{a_{\lambda}^d} \Sigma^{d
    \lambda} M \label{eq:SCV}
  \end{equation}
  is null homotopic, as the first map $c$ factors as
  \[ \tau_{> b} M \longrightarrow \Sigma^{2 p^n} M \xrightarrow{\tmop{Nm}
     (\hat{\mu})} M \]
  since the composite $\mathbb{T}$-equivariant map $\tau_{> b} M
  \xrightarrow{c} M \rightarrow M / \tmop{Nm} (\hat{\mu})$ is null-homotopic,
  by $b$-truncatedness of $M / \mu$, where $\tmop{Nm}$ is the norm map
  $R_{h\mathbb{T}} [1] \longrightarrow R^{h\mathbb{T}}$, and thus $\tmop{Nm}
  (\hat{\mu})$ is an equivariant lift of $\mu$. But now, the composite
  $\mathbb{T}$-equivariant map
  \[ \Sigma^{2 p^n} M \longrightarrow M \xrightarrow{a_{\lambda}^d} \Sigma^{d
     \lambda} R \]
  is null-homotopic. This gives rise to a null homotopy of the composite
  $\mathbb{T}$-equivariant map \eqref{eq:SCV}.
\end{proof}

\begin{corollary}
  \label{cor:Segal-BP-n-bnd-TR}Let $M$ be a right $\tmop{THH} (\tmop{BP}
  \langle n \rangle)$-module in cyclotomic spectra such that $M / (p, v_1,
  \ldots, v_n)$ satisfies Segal conjecture. Then the spectrum $\tmop{TR} (M) /
  (p, v_1, \ldots, v_{n + 1})$ is $t$-bounded above.
\end{corollary}

\begin{proof}
  Apply \Cref{prop:rel-abs-segal-conj,lem:SCV-mod-Boekstedt} to the Bökstedt
  ring $\tmop{THH} (\tmop{BP} \langle n \rangle) / (p, v_1, \ldots, v_{n +
  1})$ as in \Cref{ex:Boekstedt-THH-BP-n}.
\end{proof}

\begin{example}
  Let $\mathcal{C}$ be a bounded smooth DG-category (over $\mathbb{Z}$). Then
  by \Cref{cor:Segal-BP-n-bnd-TR} and one of main results in
  {\cite{Mao2024a}}, the spectrum $\tmop{TR} (\mathcal{C}) / (p, v_1)$ is
  $t$-bounded.
\end{example}

\begin{example}
  Let $S$ be a perfectoid ring, and $\mathcal{C}$ a bounded smooth $S$-linear
  DG-category. Then by \Cref{cor:Segal-BP-n-bnd-TR} and one of main results in
  {\cite{Mao2024a}}, the spectrum $\tmop{TR} (\mathcal{C}) / (p, v_1)$ is
  $t$-bounded. In particular, when $S$ is a perfect $\mathbb{F}_p$-algebra,
  $\tmop{TR} (\mathcal{C})$ is $p$-completely $t$-bounded above: since
  $\tmop{TR} (\mathbb{F}_p) =\mathbb{Z}_p$, we see that $v_1 = 0$ on
  $\tmop{TR} (\mathcal{C})$.
\end{example}

For modules over a Bökstedt ring, we can simply define the de-completed
$\mathbb{T}$-Tate construction by taking homotopy $\mathbb{T}$-fixed points.
This is justified by the following, an analogue of
\Cref{lem:frob-twist-htpy-fixed-pts-pres-colim}.

\begin{proposition}
  \label{prop:decpl-Tate-as-localization}Let $(R, \mu)$ be a Bökstedt ring.
  Then
  \begin{enumerate}
    \item the functor
    \begin{eqnarray*}
      \tmop{RMod}_R (\tmop{Sp}^{B\mathbb{T}}) & \longrightarrow & \tmop{Sp}\\
      M & \longmapsto & (M \otimes_R^{\mathbb{L}} R^{t C_p})^{h\mathbb{T}}
    \end{eqnarray*}
    preserves sifted colimits; and
    
    \item letting $M$ be an $R$-module in $\tmop{Sp}^{B\mathbb{T}}$, the
    assembly map
    \[ M^{h\mathbb{T}} \otimes_{R^{h\mathbb{T}}, \varphi \assign
       \varphi_p^{h\mathbb{T}}}^{\mathbb{L}} R^{t\mathbb{T}} \longrightarrow
       (M \otimes_R^{\mathbb{L}} R^{t C_p})^{h\mathbb{T}} \]
    is an equivalence.
  \end{enumerate}
\end{proposition}

\begin{proof}
  We prove the both statements in the same time by examining the vertical
  morphism
  \[ \begin{array}{ccccc}
       (\Sigma M_{h\mathbb{T}}) \otimes_{R^{h\mathbb{T}},
       \varphi}^{\mathbb{L}} R^{t\mathbb{T}} & \longrightarrow &
       M^{h\mathbb{T}} \otimes_{R^{h\mathbb{T}}, \varphi} R^{t\mathbb{T}} &
       \longrightarrow & M^{t\mathbb{T}} \otimes_{R^{h\mathbb{T}}, \varphi}
       R^{t\mathbb{T}}\\
       \longdownarrow &  & \longdownarrow &  & \longdownarrow\\
       \Sigma (M \otimes_R^{\mathbb{L}} R^{t C_p})_{h\mathbb{T}} &
       \longrightarrow & (M \otimes_R^{\mathbb{L}} R^{t C_p})^{h\mathbb{T}} &
       \longrightarrow & (M \otimes_R^{\mathbb{L}} R^{t C_p})^{t\mathbb{T}}
     \end{array} \]
  of horizontal fiber sequences, where the right $R^{h\mathbb{T}}$-module
  structure on $M^{t\mathbb{T}}$ is induced by the canonical map $\tmop{can}
  \of R^{h\mathbb{T}} \rightarrow R^{t\mathbb{T}}$, while the base change on
  the first row is along the Frobenius map $\varphi \assign
  \varphi_p^{h\mathbb{T}} \of R^{h\mathbb{T}} \rightarrow R^{t\mathbb{T}}$.
  
  Let $\tilde{\mu} \in \pi_{2 p^n} (R^{h\mathbb{T}})$ be an equivariant lift
  of $\mu$, which further admits a lift to $\pi_{2 p^n} (\Sigma
  R_{h\mathbb{T}})$. Then by assumption, the map $\varphi \of R^{h\mathbb{T}}
  \rightarrow R^{t\mathbb{T}}$ is precisely localization at the class
  $\tilde{\mu}$. However, since $\tilde{\mu}$ lifts along $\tmop{Nm} \of
  \Sigma R_{h\mathbb{T}} \rightarrow R^{h\mathbb{T}}$, we have $\tmop{can}
  (\tilde{\mu}) = 0$. This implies that
  \[ M^{t\mathbb{T}} \otimes_{R^{h\mathbb{T}}, \varphi}^{\mathbb{L}}
     R^{t\mathbb{T}} = 0 \]
  (compare with \Cref{lem:frob-twist-htpy-fixed-pts-pres-colim}). Moreover, by
  Tate orbit and fixed point lemma, we have $(R^{t C_p})^{t\mathbb{T}} \simeq
  0$, and thus
  \[ (M \otimes_R^{\mathbb{L}} R^{t C_p})^{t\mathbb{T}} = 0. \]
  It follows that the two left horizontal maps are equivalences. The
  construction $\Sigma (-)_{h\mathbb{T}}$ preserves sifted colimits, and in
  particular, inverting $\tilde{\mu}$, thus the leftmost vertical map is an
  equivalence. Thus all arrows on the left square are equivalences, and the
  the result follows.
\end{proof}

\begin{corollary}
  \label{cor:TP-poly-base-chg-neg-TC}Let $\mathcal{C}$ be a dualizable
  presentable stable $\mathbb{Z}$-linear $\infty$-category. Then the canonical
  map
  \[ \tmop{TC}^- (\mathcal{C}) \otimes_{\tmop{TC}^- (\mathbb{Z}),
     \varphi}^{\mathbb{L}} \tmop{TP} (\mathbb{Z}) \longrightarrow
     \tmop{TP}_{/\mathbb{Z}}^{\tmop{poly}} (\mathcal{C}) \]
  is an equivalence after $(p, v_1)$-completion\footnote{We do not complete it
  with respect to any ``secondary filtration'' as hinted in
  {\cite[Rem~1.6]{Manam2024}}.}.
\end{corollary}

\begin{proof}
  It follows from applying \Cref{prop:decpl-Tate-as-localization} to the
  Bökstedt ring $\tmop{THH} (\mathbb{Z}) / (p, v_1)$ as in
  \Cref{ex:Boekstedt-THH-BP-n}.
\end{proof}

\begin{corollary}
  \label{cor:TP-poly-perfd-localize-neg-TC}Let $S$ be a perfectoid ring, and
  $\mathcal{C}$ a dualizable presentable stable $S$-linear $\infty$-category.
  Then the canonical map
  \[ \tmop{TC}^- (\mathcal{C}) \otimes_{\tmop{TC}^- (S), \varphi}^{\mathbb{L}}
     \tmop{TP} (S) \longrightarrow \tmop{TP}_{/\mathbb{Z}}^{\tmop{poly}}
     (\mathcal{C}) \]
  is an equivalence after $(p, I)$-completion, where $I \assign \tmop{fib}
  \left( \Prism_S \rightarrow \barPrism_S = S \right)$.
\end{corollary}

\begin{proof}
  Applying \Cref{prop:decpl-Tate-as-localization} to the Bökstedt ring
  $\tmop{THH} (S) / (p, v_1)$ as in \Cref{ex:Boekstedt-THH-perfd}, we get the
  result up to $(p, v_1)$-completion. However, the commutative ring
  $\pi_{\ast} (\tmop{TP} (S)_p^{\wedge})$, we have $v_1 \doteq I^p$ up to a
  unit (cf.~{\cite[§5]{Mao2024a}}), thus $(p, v_1)$-completion is the same as
  $(p, I)$-completion for $\tmop{TP} (S)_p^{\wedge}$-modules.
\end{proof}

\begin{corollary}
  \label{cor:TP-poly-untwist-TP}Let $R$ be a quasiregular semiperfectoid ring.
  Then the canonical map
  \[ \tmop{TC}^- (R) \longrightarrow \tmop{TP}^{\tmop{poly}} (R) \]
  is a localization after $(p, I)$-completion, where $I \assign \tmop{fib}
  \left( \Prism_R \rightarrow \barPrism_R \right)$, which can be identified
  with $\tmop{TP}^{(- 1)} (R)$ in {\cite{Manam2024}}.
\end{corollary}

\begin{proof}
  Pick a surjection $S \twoheadrightarrow R$ for a perfectoid ring $S$. The
  result follows from \Cref{cor:TP-poly-perfd-localize-neg-TC} as the
  Frobenius map $\varphi \of \tmop{TC}^- (S)_p^{\wedge} \rightarrow \tmop{TP}
  (S)_p^{\wedge}$ is already a localization (at a class in $\pi_2 (\tmop{TC}^-
  (S)_p^{\wedge})$). The identification of this localization with
  $\tmop{TP}^{(- 1)} (R)$ is established in {\cite[Rem~1.3]{Gregoric2025}}.
\end{proof}

Finally, we give a conceptual explanation\footnote{We thank Georg
{\tmname{Tamme}} for asking us such a question.} of continuity of $\tmop{TC}$
established in {\cite{Clausen2021}} via de-completed $\tmop{TC}^-$ and
$\tmop{TP}$, which is much more akin to {\cite[Prop~7.12]{Antieau2023}}.
Actually, the argument is parallel, as a consequence of the following simple
category-theoretic observation:

\begin{lemma}
  \label{lem:triangles-lim-Cart-eq}Let $\mathcal{C}$ be an $\infty$-category
  with finite limits. Let
  
  \[\begin{tikzcd} {x_{00}} & {x_{01}} \\ {x_{10}} & {x_{11}} \arrow[from=1-1,
  to=1-2] \arrow[from=1-1, to=2-1] \arrow[from=1-2, to=2-2] \arrow[from=2-1,
  to=1-2] \arrow[from=2-1, to=2-2] \end{tikzcd}\]
  
  {\noindent}be a diagram in $\mathcal{C}$, where the outer square is a
  Cartesian diagram, while the two inner triangles are {\tmstrong{not}}
  commutative.
  
  Then the map from the limit of the upper triangle (viewed as a diagram
  $\partial \Delta^2 \rightarrow \mathcal{C}$) to the limit of the lower
  triangle (viewed as a diagram $\partial \Delta^2 \rightarrow \mathcal{C}$ as
  well) is an equivalence.
\end{lemma}

\begin{proof}
  It follows directly from ``decomposing the diagram'' (e.g. we decompose the
  lower triangle via $\partial \Delta^2 = \Lambda_2^2 \amalg_{\partial
  \Delta^1} \Delta^1$, where the composite functor $\Lambda_2^2 \rightarrow
  \partial \Delta^2 \rightarrow \mathcal{C}$ is the lower cospan $(x_{1
  \nocomma 0} \rightarrow x_{1 \nocomma 1} \leftarrow x_{0 \nocomma 1})$, and
  the composite $\Delta^1 \rightarrow \partial \Delta^2 \rightarrow
  \mathcal{C}$ is the ``diagonal'' arrow $(x_{1 \nocomma 0} \rightarrow x_{0
  \nocomma 1})$, and $x_{0 \nocomma 0}$ is the limit of the first composite
  functor $\Lambda_2^2 \rightarrow \mathcal{C}$). Alternatively, the proof of
  {\cite[Prop~7.12]{Antieau2023}} works in this case as well.
\end{proof}

\begin{proposition}
  \label{prop:TC-Boekstedt-continuity}Let $(R, \mu)$ be a Bökstedt ring. Then
  the composite functor
  \[ \tmop{RMod}_R (\tmop{CycSp}_p) \longrightarrow \tmop{CycSp}_p
     \xrightarrow{\tmop{TC}} \tmop{Sp} \]
  preserves sifted colimits.
\end{proposition}

\begin{proof}
  Let $M$ be a right $R$-module in cyclotomic spectra. Then the Frobenius map
  $\varphi \of M^{h\mathbb{T}} \rightarrow M^{t\mathbb{T}}$ factors as
  $M^{h\mathbb{T}} \xrightarrow{\tilde{\varphi}} (M \otimes_R^{\mathbb{L}}
  R^{t C_p})^{h\mathbb{T}} \rightarrow M^{t\mathbb{T}}$. Now we formulate the
  triangle $\partial \Delta^2 \rightarrow \tmop{Sp}$ given by
  
  \[\begin{tikzcd} & {(M\otimes_R^{\mathbb L}R^{tC_p})^{h\mathbb T}} \\
  {M^{h\mathbb T}} & {M^{t\mathbb T}} \arrow[from=1-2, to=2-2]
  \arrow[from=2-1, to=1-2] \arrow["{\mathrm{can}}", from=2-1, to=2-2]
  \end{tikzcd}\]
  
  {\noindent}where the ``diagonal'' map is $\tilde{\varphi}$. Note that
  $\tmop{TC} (M)$ is precisely the limit of this triangle. Let $P (M)$ denote
  the pullback
  \[ M^{h\mathbb{T}} \times_{\tmop{can}, M^{t\mathbb{T}}} (M
     \otimes_R^{\mathbb{L}} R^{t C_p})^{h\mathbb{T}}, \]
  which gives rise to a functor $P \of \tmop{RMod}_R (\tmop{CycSp}_p)
  \rightarrow \tmop{Sp}$, understood as ``de-completed homotopy fixed
  points''. By \Cref{lem:triangles-lim-Cart-eq}, the spectrum $\tmop{TC} (M)$
  can also be computed by the limit of the upper triangle
  
  \[\begin{tikzcd} {P(M)} & {(M\otimes_R^{\mathbb L}R^{tC_p})^{h\mathbb T}} \\
  {M^{h\mathbb T}} \arrow[from=1-1, to=1-2] \arrow[from=1-1, to=2-1]
  \arrow[from=2-1, to=1-2]
  \end{tikzcd}\]
  
  {\noindent}Since we have a fiber sequence $\Sigma M_{h\mathbb{T}}
  \rightarrow P (M) \rightarrow (M \otimes_R^{\mathbb{L}} R^{t
  C_p})^{h\mathbb{T}}$, it follows from \Cref{prop:decpl-Tate-as-localization}
  that $P$ preserves sifted colimits, and so does $\tmop{TC}$ on
  $\tmop{RMod}_R (\tmop{CycSp}_p)$, as $\tmop{TC} (M)$ is an equalizer of $P
  (M) \rightrightarrows (M \otimes_R^{\mathbb{L}} R^{t C_p})^{h\mathbb{T}}$,
  where we invoke \Cref{prop:decpl-Tate-as-localization} again.
\end{proof}

\begin{corollary}
  \label{cor:cts-TC-p-v-cpl}Let $n \in \mathbb{Z}_{\geqslant - 1}$. Then the
  exact functor
  \[ \tmop{RMod}_{\tmop{THH} (\tmop{BP} \langle n \rangle)} (\tmop{CycSp}_p)
     \rightarrow \tmop{CycSp}_{(p, v_1, \ldots, v_{n + 1})}^{\wedge}
     \xrightarrow{\tmop{TC} (-)^{\wedge}} \tmop{Sp}_{(p, v_1, \ldots, v_{n +
     1})}^{\wedge} \]
  preserves small colimits.
\end{corollary}

\begin{proof}
  It follows from applying \Cref{prop:TC-Boekstedt-continuity} to the
  Bökstedt ring $\tmop{THH} (\tmop{BP} \langle n \rangle) / (p, v_1, \ldots,
  v_{n + 1})$ as in \Cref{ex:Boekstedt-THH-BP-n}.
\end{proof}

\begin{example}
  Let $(M_{\alpha})$ be a uniformly $t$-bounded below filtered diagram of
  $\tmop{THH} (\tmop{BP} \langle n \rangle)$-modules in $\tmop{CycSp}$. Then
  $(\tmop{colim}_{\alpha} \tmop{TC} (M_{\alpha}))_p^{\wedge}$ and $\tmop{TC}
  (\tmop{colim}_{\alpha} M_{\alpha})_p^{\wedge}$ are $t$-bounded below. It
  follows that they are automatically $(p, v_1, \ldots, v_{n + 1})$-complete,
  and hence the canonical map
  \[ \underset{\alpha}{\tmop{colim}} \tmop{TC} (M_{\alpha}) \longrightarrow
     \tmop{TC} \left( \underset{\alpha}{\tmop{colim}} M_{\alpha} \right) \]
  is an equivalence after $p$-completion. This reproduces
  {\cite[Thm~2.7]{Clausen2021}} for $\tmop{THH} (\tmop{BP} \langle n
  \rangle)$-modules in $\tmop{CycSp}$. When the filtered system is not
  uniformly $t$-bounded below, \Cref{cor:cts-TC-p-v-cpl} is potentially more
  informative.
\end{example}

\section{Completeness of motivic filtration}

In this section, we will give a sufficient condition
(\Cref{thm:mot-cpl-TP-poly}) for completeness of motivic filtration on
de-completed topological periodic cyclic homology. We also establishes a Segal
conjecture (\Cref{thm:Segal-conj-F-sm}) when the ring in question is
additionally $F$-smooth. For this purpose, we first introduce a weak finite
generation condition (\Cref{def:weak-ft}) for modules, based on Shannon's
notion of rank of flat modules. Then we deduce a connectivity
(\Cref{cor:Nygaard-gr-unif-bnd}) for associated graded pieces of absolute
prismatic cohomology with respect to Nygaard filtration under this finite
generation of cotangent complex. The main completeness result then follows.

We start with a higher algebra version of Shannon's rank of flat modules:

\begin{definition}[{\cite{Shannon1970}}]
  \label{def:rk-flat-mods}Let $R$ be a connective $\mathbb{E}_1$-ring. We say
  that a flat right $R$-module $P$ is {\tmdfn{of finite rank}} (abbrev.
  $\tmop{rk}_R (P) < \infty$) if there exists a natural number $r \in
  \mathbb{N}$ such that $P$ is a filtered colimit of finite free right
  $R$-modules of rank$\nosymbol \leq r$, and in this case, we say that $P$ is
  of rank$\nosymbol \leq r$ (abbrev. $\tmop{rk}_R (P) \leq r$).
\end{definition}

\begin{example}[Ind-Zariski localizations]
  \label{ex:Zar-loc-rk-1}The flat $\mathbb{Z}$-module $\mathbb{Q}$ is of
  rank$\nosymbol \leq 1$, as we can write $\mathbb{Q}$ as the sequential
  colimit of $\mathbb{Z} \hookrightarrow (2!)^{- 1} \mathbb{Z} \hookrightarrow
  (3!)^{- 1} \mathbb{Z} \hookrightarrow \cdots$. More generally, let $R$ be an
  $\mathbb{E}_1$-ring, $M$ a finite free right $R$-module, and $S \subseteq
  \pi_0 (R)$ a multiplicative subset satisfying right Ore condition in
  $\pi_{\ast} (R)$. Then the localization $S^{- 1} M$, as a flat right
  $R$-module spectrum, is of rank$\nosymbol \leq 1$.
\end{example}

\begin{example}[Standard étale maps]
  \label{ex:et-loc-fin-rk}Let $R$ be a commutative ring, $g \in R [t]$ an
  $R$-polynomial, and $f \in R [t]$ a monic $R$-polynomial. Let $S \assign R
  [t]_f / (g)$. Then $\tmop{rk}_R (S) \leq \deg (g)$. Recall that, when the
  derivative $f'$ is invertible in $S$, such a map $R \rightarrow S$ is called
  a {\tmdfn{standard étale map}}
  {\cite[\href{https://stacks.math.columbia.edu/tag/00UB}{Tag
  00UB}]{stacks-project}}.
\end{example}

\begin{remark}
  \label{rem:rk-sum-tensor-prod}It follows from definition that
  \begin{itemize}
    \item Let $R \rightarrow S$ be a map of connective $\mathbb{E}_1$-rings,
    and $M$ a flat right $R$-module. Then $\tmop{rk}_S (M
    \otimes_R^{\mathbb{L}} S) \leq \tmop{rk}_R (M)$.
    
    \item Let $R$ be a connective $\mathbb{E}_1$-ring, and $M$ and $N$ two
    flat right $R$-modules. Then $\tmop{rk}_R (M \oplus N) \leq \tmop{rk}_R
    (M) + \tmop{rk}_R (N)$.
    
    \item Let $R$ be a connective $\mathbb{E}_{\infty}$-ring, and $M$ and $N$
    two flat $R$-modules. Then $\tmop{rk}_R (M \otimes_R^{\mathbb{L}} N) \leq
    \tmop{rk}_R (M) \tmop{rk}_R (N)$.
    
    \item Let $R$ and $S$ be two connective $\mathbb{E}_1$-rings, $M$ a flat
    right $R$-module, and $N$ a flat right $S$-module. Then $M \times N$ is a
    flat $R \times S$-module, with $\tmop{rk}_{R \times S} (M \times N) \leq
    \max \{ \tmop{rk}_R (M), \tmop{rk}_R (N) \}$.
  \end{itemize}
\end{remark}

We give a criterion for flat modules being of finite rank, as a higher algebra
generalization of {\cite[Thm~2]{Shannon1970}}.

\begin{notation}
  Let $R$ be a connective $\mathbb{E}_1$-ring, and $n \in \mathbb{N}$. We will
  denote by $R^{\oplus n}$, or $R^n$ when there is no ambiguity, the finite
  free (right) $R$-module of rank $n$. More generally, let $M$ be a right
  $R$-module spectrum. We will denote by $M^{\oplus n}$ the direct sum of $n$
  copies of $M$.
\end{notation}

\begin{proposition}
  \label{prop:crit-flat-rank-r}Let $R$ be a connective $\mathbb{E}_1$-ring, $r
  \in \mathbb{N}$, and $P$ a flat right $R$-module. Then $P$ is of
  rank$\nosymbol \leq r$ if and only if, for every finite free right
  $R$-module $L$, every map $L \rightarrow P$ factors through some $R^s$ for
  $s \leq r$.
\end{proposition}

\begin{proof}
  
  \begin{description}
    \item[``Only if'' part] We write $P$ as a filtered colimit
    $\tmop{colim}_{\alpha} L_{\alpha}$ with every $L_{\alpha}$ being finite
    free of rank$\nosymbol \leq r$. By compactness of $L$ in $\tmop{RMod}_R$,
    we have
    \[ \tmop{Hom}_{\tmop{RMod}_R} (L, P) = \underset{\alpha}{\tmop{colim}}
       \tmop{Hom}_{\tmop{RMod}_R} (L, L_{\alpha}) \]
    and the result follows.
    
    \item[``If'' part] It suffices to check that the slice category
    $(\tmop{RMod}_R^{\tmop{free}, \tmop{rk} \leqslant r})_{/ P}$ is filtered.
    By Lazard's theorem {\cite[Thm~7.2.2.15]{Lurie2017}}, the category
    $(\tmop{RMod}_R^{\tmop{free}, \tmop{rk} < \infty})_{/ P}$ is filtered. In
    particular, any diagram $q \of K \rightarrow (\tmop{RMod}_R^{\tmop{free},
    \tmop{rk} \leqslant r})_{/ P}$ extends to a diagram $q^{\vartriangleright}
    \of K^{\vartriangleright} \rightarrow (\tmop{RMod}_R^{\tmop{free},
    \tmop{rk} < \infty})_{/ P}$. By assumption, the map $q^{\vartriangleright}
    (\infty) \rightarrow P$ factors through some $R^s$ for $s \leq r$, where
    $\infty \in K^{\vartriangleright}$ is the cocone point. The result then
    follows.
  \end{description}
\end{proof}

It follows that

\begin{corollary}
  \label{cor:rk-pi-0}Let $R$ be a connective $\mathbb{E}_1$-ring, $r \in
  \mathbb{N}$, and $P$ a flat right $R$-module. Then $P$ is of rank$\nosymbol
  \leq r$ if and only if so is the flat right $\pi_0 (R)$-module $\pi_0 (P)$.
\end{corollary}

As mentioned in the last paragraph of {\cite{Shannon1970}} that, the rank of
flat modules coincides with the usual rank over principal ideal domains.
Actually, the same argument works over Bézout domains, and in particular,
over valuation rings, which we record as follows.

\begin{lemma}
  \label{lem:rk-flat-Bezout}Let $R$ be a Bézout domain (e.g. a valuation
  ring), $K \assign \tmop{Frac} (R)$, and $M$ a flat right $R$-module. Then
  $\tmop{rk}_R (M) = \dim_K (M \otimes_R K)$.
\end{lemma}

\begin{proof}
  First, suppose that $\tmop{rk}_R (M) \leq r$ for some $r \in \mathbb{N}$.
  Then $M$ is a filtered colimit of finite free $R$-modules of rank$\nosymbol
  \leq r$, which implies that $M \otimes_R K$ is a filtered colimit of
  $K$-vector spaces of dimension$\nosymbol \leq r$, thus $\dim_K (M \otimes_R
  K) \leq r$.
  
  On the other hand, suppose that $\dim_K (M \otimes_R K) \leq r$. Since $M$
  is flat, we know that the canonical map $M \rightarrow M \otimes_R K$ is
  injective. Given any $(r + 1)$-tuple $(x_1, \ldots, x_{r + 1}) \in M^{\oplus
  (r + 1)}$, we see that the tuple $(x_1 \otimes 1, \ldots, x_{r + 1} \otimes
  1) \in (M \otimes_R K)^{\oplus (r + 1)}$ is $K$-linear dependent. This
  implies that, there exists a tuple $(a_1, \ldots, a_{r + 1}) \in R^{r + 1}$
  such that $a_1 x_1 + \cdots + a_{r + 1} x_{r + 1} = 0$. Since $R$ is a
  Bézout domain, we may divide $a_i$ by $\gcd (a_1, \ldots, a_{r + 1})$,
  which reduces to the case that $\{ a_1, \ldots, a_{r + 1} \} \subseteq R$
  generates the unit ideal.
  
  We now show that the map $f \of R^{r + 1} \rightarrow M, (b_1, \ldots, b_{r
  + 1}) \mapsto \sum x_i b_i$ factors through $R^r$. Indeed, we note that
  $(a_1, \ldots, a_{r + 1}) \in \ker (f)$, which implies that $f$ factors
  through $Q \assign R^{r + 1} / (a_1, \ldots, a_{r + 1}) R$. However, since
  $\{ a_1, \ldots, a_{r + 1} \} \subseteq R$ generates the unit ideal, the map
  $R \rightarrow R^{n + 1}, r \mapsto (a_1 r, \ldots, a_{r + 1} r)$ of
  $R$-modules admits a retract. It follows that $Q$ is a proper direct summand
  of $R^{r + 1}$. Since $R$ is a Bézout domain, we deduce that $Q$ is a
  finite free $R$-module of rank $r$. By {\cite[Thm~2]{Shannon1970}} (or
  \Cref{prop:crit-flat-rank-r}), we see that $\tmop{rk}_R (M) \leq r$.
\end{proof}

We deduce a descent property of this finite generation condition.

\begin{proposition}
  \label{prop:fin-rk-desc-fin-rk}Let $R \rightarrow S$ be a faithfully flat
  map of connective $\mathbb{E}_{\infty}$-rings with $\tmop{rk}_R (S) <
  \infty$, and $M$ a flat $R$-module. Suppose that the flat $S$-module $M
  \otimes_R^{\mathbb{L}} S$ is of finite rank, then so is the flat $R$-module
  $M$, with inequality
  \[ \tmop{rk}_R (M) \leq \tmop{rk}_R (S) \tmop{rk}_S (M
     \otimes_R^{\mathbb{L}} S) . \]
\end{proposition}

\begin{proof}
  Suppose that $\tmop{rk}_S (M \otimes_R^{\mathbb{L}} S) \leq r$ and
  $\tmop{rk}_R (S) \leq \mu$. By \Cref{cor:rk-pi-0}, we may assume that $R$ is
  concentrated in degree $0$. Let $(x_1, \ldots, x_n) \in M^{\oplus n}$
  corresponding to a map $R^n \rightarrow M$ of $R$-modules. By
  {\cite[Thm~2]{Shannon1970}} (or \Cref{prop:crit-flat-rank-r}), there exists
  $(y_1, \ldots, y_s) \in (M \otimes_R S)^{\oplus s}$ with $s \leq r$, such
  that $x_i \otimes 1$ lies in $\tmop{Span} \{ y_1, \ldots, y_s \}$. Since $S$
  is a filtered colimit of finite free $R$-modules of rank$\nosymbol \leq
  \mu$, we may choose some map $R^t \rightarrow S$ with $t \leq \mu$, and each
  $y_i \in M \otimes_R S$ admits a lift $\tilde{y}_i \in M \otimes_R R^{\oplus
  t} \cong M^{\oplus t}$, and we write $\tilde{y}_i = (z_{i \nocomma 1},
  \ldots, z_{i \nocomma t})$.
  
  We claim that each $x_j$ belongs to $\tmop{Span} \{ z_{1 \nocomma 1}, z_{1
  \nocomma 2}, \ldots, z_{s \nocomma t} \}$. Indeed, by construction, we have
  $\tilde{y}_i$ lies in $\tmop{Span} \{ z_{i \nocomma 1} \otimes 1, \ldots,
  z_{i \nocomma t} \otimes 1 \} \subseteq M \otimes_R R^{\oplus t}$, and thus
  $y_i \in \tmop{Span} \{ z_{i \nocomma 1} \otimes 1, \ldots, z_{i \nocomma t}
  \otimes 1 \} \subseteq M \otimes_R S$. Consequently, each $x_j \otimes 1 \in
  M \otimes_R S$ belongs to $\tmop{Span} \{ z_{1 \nocomma 1} \otimes 1, z_{1
  \nocomma 2} \otimes 1, \ldots, z_{s \nocomma t} \otimes 1 \}$. The result
  then follows from the faithful flatness of $R \rightarrow S$.
\end{proof}

\begin{example}[Zariski descent]
  Let $R$ be a connective $\mathbb{E}_{\infty}$-ring, $d \in \mathbb{N}$, and
  $(f_1, \ldots, f_d) \in \pi_0 (R)^d$ a tuple of elements generating the unit
  ideal of $\pi_0 (R)$. Then the map $R \rightarrow \prod_{i = 1}^d R_{f_i}
  \backassign S$ is a Zariski cover. By
  \Cref{ex:Zar-loc-rk-1,rem:rk-sum-tensor-prod}, the flat $R$-module $S$ is of
  rank$\nosymbol \leq d$.
  
  Let $M$ be a flat $R$-module such that each localization $M_{f_i}$ is a flat
  $R_{f_i}$-module of finite rank. Then by \Cref{rem:rk-sum-tensor-prod}, the
  flat $S$-module $M \otimes_R^{\mathbb{L}} S$ is of finite rank as well, and
  it follows from \Cref{prop:fin-rk-desc-fin-rk} that the flat $R$-module $M$
  is of finite rank.
\end{example}

\begin{example}[Étale descent]
  Let $R \rightarrow S$ be a faithfully flat étale map of connective
  $\mathbb{E}_{\infty}$-rings. By
  {\cite[\href{https://stacks.math.columbia.edu/tag/00UE}{Tag
  00UE}]{stacks-project}}, we may choose a Zariski cover $S \rightarrow
  \prod_{i = 1}^d S_{f_i} \backassign T$ such that each composite $R
  \rightarrow S \rightarrow S_{f_i}$ is a standard étale map on $\pi_0$. By
  \Cref{ex:et-loc-fin-rk,rem:rk-sum-tensor-prod,cor:rk-pi-0} that the flat
  $R$-module $T$ is of finite rank.
  
  Let $M$ be a flat $R$-module such that the flat $S$-module $M
  \otimes_R^{\mathbb{L}} S$ is of finite rank. Then by
  \Cref{rem:rk-sum-tensor-prod}, the flat $T$-module $M \otimes_R^{\mathbb{L}}
  T$ is of finite rank, and it follows from \Cref{prop:fin-rk-desc-fin-rk}
  that the flat $R$-module $M$ is of finite rank as well.
\end{example}

\begin{remark}
  There is a variant of \Cref{def:rk-flat-mods} when $R$ is a connective
  $\mathbb{E}_{\infty}$-rings, where we replace finite free right $R$-modules
  of rank$\nosymbol \leq r$ by finite {\tmem{projective}} $R$-modules of
  rank$\nosymbol \leq r$. For that version, previously statements adapt, and
  one can generalize Bézout domains to Prüfer domains in
  \Cref{lem:rk-flat-Bezout}.
\end{remark}

Now we come back to our weakening of finite generation condition.

\begin{definition}
  \label{def:weak-ft}Let $R$ be a connective $\mathbb{E}_1$-ring, and $r \in
  \mathbb{N}$. We say that a connective right $R$-module $M$ is {\tmdfn{weakly
  finitely generated}}\footnote{As pointed out by ChatGPT-6 Astra, our concept
  is different from the one in {\cite[§8.2]{Bazzoni2022}} under the same
  name. In the special case of $R$ being a commutative chain ring and $M$
  being an $R$-module concentrated in degree $0$, this coincides with
  {\tmdfn{Fleischer rank}} in {\cite[Rem~3.1]{Couchot2013}}.} (of
  rank$\nosymbol \leq r$) if, there exists a flat right $R$-module $P$ of
  finite rank (of rank$\nosymbol \leq r$), along with a surjective\footnote{A
  map of connective $R$-modules is {\tmdfn{surjective}} if it is surjective on
  $\pi_0$.} map $P \rightarrow M$ of $R$-modules.
\end{definition}

\begin{example}
  Let $R$ be a connective $\mathbb{E}_1$-ring, and $M$ a connective right
  $R$-module which is finitely generated (i.e. $\pi_0 (M)$ is a finitely
  generated right $\pi_0 (R)$-module). Then $M$ is weakly finitely generated.
  This justifies the terminology. Slightly more generally, weak finite
  generation is a condition on $\pi_0$ under countable generation hypothesis.
  We collect relevant results in \Cref{app:weak-fg-countable-gen}.
\end{example}

It follows from \Cref{prop:crit-flat-rank-r} that

\begin{corollary}
  \label{cor:flat-mod-weak-fg-fin-rk}Let $R$ be a connective
  $\mathbb{E}_1$-ring, and $M$ a flat right $R$-module. Then $M$ is weakly
  finitely generated if and only if it is of finite rank.
\end{corollary}

A merit of this finite generation condition is an asymptotic connectivity of
exterior powers.

\begin{lemma}
  \label{lem:ft-ext-pow-unif-bnd-below}Let $R$ be an animated ring, and $M$ a
  connective $R$-module weakly finitely generated. Then the $R$-module spectra
  $\left\{ \left( \bigwedgestar_R^j M \right) [- j] \barsuchthat j \in
  \mathbb{N} \right\}$ are uniformly $t$-bounded below, and the bound only
  depends on the rank of $M$.
\end{lemma}

\begin{proof}
  By definition, there exists a natural number $r \in \mathbb{N}$, a flat
  $R$-module $P$ which is a filtered colimit of finite free $R$-modules of
  rank$\nosymbol \leq r$, and a surjective map $P \rightarrow M$ of connective
  $R$-modules. Let $K \assign \tmop{fib} (P \rightarrow M)$, which is
  connective as well, and we get a fiber sequence $P \rightarrow M \rightarrow
  K [1]$. It follows that, for every $j \in \mathbb{N}$, the $R$-module
  spectrum $\left( \bigwedgestar_R^j M \right) [- j]$ admits a finitary
  filtration (cf.~{\cite[§V.4]{Illusie1971}}), whose associated graded pieces
  are given by
  \[ \left( \bigwedgestar_R^i P \right) \otimes_R^{\mathbb{L}} \left(
     \bigwedgestar_R^{j - i} K [1] \right) [- j] \simeq \left(
     \bigwedgestar_R^i P \right) \otimes_R^{\mathbb{L}} (\Gamma_R^{j - i} K)
     [- i], \]
  where we used Quillen and Illusie's décalage
  {\cite[I.4.3.2.1]{Illusie1971}}, see also {\cite[Prop~25.2.4.2]{Lurie2018}}.
  Note that $\bigwedgestar_R^i P = 0$ for $i > r$ by Lazard's theorem
  {\cite[Thm~7.2.2.15]{Lurie2017}}, as $\bigwedgestar_R^i (-)$ preserves
  filtered colimits of $R$-module spectra, and it vanishes on finite free
  $R$-modules of rank$\nosymbol \leq r$. It follows that $\left(
  \bigwedgestar_R^j M \right) [- j]$ is $(- r)$-connective.
\end{proof}

Let $R$ be an animated ring. {\cite[Rem~5.5.9]{Bhatt2022}} showed that
$\tmop{gr}_N^m  \Prism_R \{ n \}$ is bounded below. Combining with
\Cref{lem:ft-ext-pow-unif-bnd-below}, we adapt the same argument to show that,
under such a weak finite generation, they are uniformly bounded below. Recall
that the ($p$-completed) diffracted Hodge cohomology $\Omega_R^{\DHod}$ admits
an increasing filtration $\tmop{Fil}_{\ast}^{\tmop{conj}} \Omega_R^{\DHod}$,
called the {\tmdfn{conjugate filtration}}, {\cite[Cons~4.7.1]{Bhatt2022}}.
Connectivity of exterior powers gives rise to that of several cohomology
theories:

\begin{lemma}
  \label{lem:DHod-conn}Let $R$ be an animated ring such that the spectra
  \[ \left\{ \left( \bigwedgestar_R^n L_{R /\mathbb{Z}} \right) [- n]
     \otimes_{\mathbb{Z}}^{\mathbb{L}} \mathbb{F}_p \barsuchthat n \in
     \mathbb{Z} \right\} \]
  is uniformly $t$-bounded below. Then the $R$-module spectra
  \[ \left\{ \tmop{Fil}_s^{\tmop{conj}} \Omega_R^{\DHod}
     \otimes_{\mathbb{Z}}^{\mathbb{L}} \mathbb{F}_p \barsuchthat s \in
     \mathbb{N} \right\} \]
  are uniformly $t$-bounded below, and the bound is the same as the previous
  bound.
\end{lemma}

\begin{proof}
  {\cite[Rem~4.5.3]{Bhatt2022}} identifies the associated graded pieces
  $\tmop{gr}_{\ast}^{\tmop{conj}} \Omega_R^{\DHod}$ with ($p$-completed)
  shifted exterior powers $\left( \bigwedgestar_R^{\ast} L_{R /\mathbb{Z}}
  \right) \left[ - \mathord{\ast} \right]$, and the result then follows.
\end{proof}

\begin{corollary}
  \label{cor:Nygaard-gr-unif-bnd}Let $R$ be an animated ring such that the
  spectra
  \[ \left\{ \left( \bigwedgestar_R^n L_{R /\mathbb{Z}} \right) [- n]
     \otimes_{\mathbb{Z}}^{\mathbb{L}} \mathbb{F}_p \barsuchthat n \in
     \mathbb{Z} \right\} \]
  is uniformly $t$-bounded below. Then the $R$-module spectra
  \[ \left\{ \tmop{gr}_N^m  \Prism_R \{ n \} \otimes_{\mathbb{Z}}^{\mathbb{L}}
     \mathbb{F}_p \barsuchthat (m, n) \in \mathbb{Z} \times \mathbb{Z}
     \right\} \]
  and
  \[ \left\{ \tmop{Fil}_m^{\tmop{conj}}  \barPrism_R \{ n \}
     \otimes_{\mathbb{Z}}^{\mathbb{L}} \mathbb{F}_p \barsuchthat (m, n) \in
     \mathbb{Z} \times \mathbb{Z} \right\} \]
  are uniformly $t$-bounded below, and the bound only depends on the previous
  uniform $t$-bound. In particular, the spectra $\left\{ \hatPrism_R \{ n \}
  \barsuchthat n \in \mathbb{Z} \right\}$ and $\left\{ \Prism_R \{ n \}
  \barsuchthat n \in \mathbb{Z} \right\}$ are $p$-completely\footnote{That is
  to say, they are uniformly $t$-bounded below after $(-)
  \otimes_{\mathbb{Z}}^{\mathbb{L}} \mathbb{F}_p$.} uniformly $t$-bounded
  below.
\end{corollary}

\begin{proof}
  It follows from \Cref{lem:DHod-conn} and the fiber sequences
  \[ \tmop{gr}_N^m  \Prism_R \{ n \} \longrightarrow
     \tmop{Fil}_m^{\tmop{conj}} \Omega_R^{\DHod} \xrightarrow{\Theta + m}
     \tmop{Fil}_{m - 1}^{\tmop{conj}} \Omega_R^{\DHod} \]
  and
  \[ \tmop{Fil}_m^{\tmop{conj}}  \barPrism_R \{ n \} \longrightarrow
     \tmop{Fil}_m^{\tmop{conj}} \Omega_R^{\DHod} \xrightarrow{\Theta + n}
     \tmop{Fil}_m^{\tmop{conj}} \Omega_R^{\DHod} \]
  as established in {\cite[Rem~4.7.5~\&~5.5.8]{Bhatt2022}}.
\end{proof}

These connectivity estimates are sufficient for establishing the completeness
of motivic filtration on $\tmop{TP}_{/\mathbb{Z}}^{\tmop{poly}}$. Recall that
it is exhaustive (\Cref{cor:exhaustive-mot-fil-TP-poly}).

\begin{theorem}
  \label{thm:mot-cpl-TP-poly}Let $R$ be an animated ring such that such that
  the spectra
  \[ \left\{ \left( \bigwedgestar_R^n L_{R /\mathbb{Z}} \right) [- n]
     \otimes_{\mathbb{Z}}^{\mathbb{L}} \mathbb{F}_p \barsuchthat n \in
     \mathbb{Z} \right\} \]
  is uniformly $t$-bounded below. Then $\tmop{Fil}_M^n  \left(
  \tmop{THH}_{\Prism} (R)^{\Phi C_p} \right)_p^{\wedge}$ is linearly
  increasingly connective as $n \rightarrow \infty$, and the motivic
  filtration on $\tmop{TP}_{/\mathbb{Z}}^{\tmop{poly}} (R)$ is complete.
\end{theorem}

\begin{proof}
  Recall that the motivic filtration on $\tmop{THH} (R)$ is complete
  {\cite[Thm~6.2.4 \& Prop~6.2.12]{Bhatt2022}}. By
  \Cref{rem:Frob-THH-Z-inv-mu1}, we are reduced to show that, after inverting
  $\mu_1$ in $\tmop{Fil}_M (\tmop{THH} (R) \otimes_{\mathbb{Z}}^{\mathbb{L}}
  \mathbb{F}_p)$, the motivic filtration remains complete. We prove this by
  connectivity. Indeed, the $r$-th filtration piece of $\tmop{Fil}_M
  (\tmop{THH} (R) \otimes_{\mathbb{Z}}^{\mathbb{L}} \mathbb{F}_p) [\mu_1^{-
  1}]$ is the sequential colimit of
  \[ \tmop{Fil}_M^r (\tmop{THH} (R) \otimes_{\mathbb{Z}}^{\mathbb{L}}
     \mathbb{F}_p) \xrightarrow{\mu_1} \tmop{Fil}_M^{r + p} (\tmop{THH} (R)
     \otimes_{\mathbb{Z}}^{\mathbb{L}} \mathbb{F}_p) [- 2 p]
     \xrightarrow{\mu_1} \cdots \]
  We note that $\{ \tmop{Fil}_M^n \tmop{THH} (R) [- 2 n] \barsuchthat n \in
  \mathbb{Z} \}$ is $p$-completely $t$-bounded below by motivic completeness
  on $\tmop{THH} (R)$ and \Cref{cor:Nygaard-gr-unif-bnd}. It follows that
  \[ \{ \tmop{Fil}_M^{r + n} \tmop{THH} (R) [- 2 n] \}_{n \in \mathbb{Z}} = \{
     \Sigma^{2 r} \tmop{Fil}_M^n \tmop{THH} (R) [- 2 n] \}_{n \in \mathbb{N}}
  \]
  is $p$-completely uniformly $t$-bounded below, and the given bound tends to
  $\infty$ as $r \rightarrow \infty$ with linear growth. In particular, the
  motivic filtration on $\left( \tmop{THH}_{\Prism} (R)^{\Phi C_p}
  \right)_p^{\wedge}$ is complete, and it remains complete after taking
  $(-)^{h\mathbb{T}_{\tmop{ev}}}$.
\end{proof}

This also implies that, when $R$ is additionally $F$-smooth à la
Bhatt--Mathew {\cite{Bhatt2023a}}, the cyclotomic spectrum $\tmop{THH} (R) /
(p, v_1)$ satisfies Segal conjecture relative to $\tmop{THH} (\mathbb{Z}) /
(p, v_1)$, by passing to associated graded pieces of motivic filtrations.
Actually, the same argument implies that absolute Segal conjecture holds,
which is slightly stronger by \Cref{prop:rel-abs-segal-conj}.

\begin{theorem}
  \label{thm:Segal-conj-F-sm}Let $R$ be an $F$-smooth $p$-quasisyntomic
  commutative ring. Suppose such that the spectra
  \[ \left\{ \left( \bigwedgestar_R^n L_{R /\mathbb{Z}} \right) [- n]
     \otimes_{\mathbb{Z}}^{\mathbb{L}} \mathbb{F}_p \barsuchthat n \in
     \mathbb{Z} \right\} \]
  is uniformly $t$-bounded below. Then the fiber $\tmop{fib} (\tmop{THH} (R)
  \rightarrow \tmop{THH} (R)^{t C_p})$ of the cyclotomic Frobenius map is
  $p$-completely $t$-bounded above.
\end{theorem}

\begin{proof}
  {\cite[Rem~4.11]{Bhatt2023a}} showed\footnote{We thank Hyungseop
  {\tmname{Kim}} for pointing this out to us.} that the filtered Frobenius map
  \[ \tmop{Fil}_N^{\ast}  \Prism_R \longrightarrow \Prism_R^{[\ast]} \]
  induces an equivalence
  \[ \tmop{gr}_N^i  \Prism_R \longrightarrow \tau_{\geqslant - i} \tmop{gr}^i 
     \Prism_R^{[\ast]} \simeq \tau_{\geqslant - i}  \barPrism_R \{ i \} . \]
  By \Cref{cor:Nygaard-gr-unif-bnd}, the Hodge--Tate cohomology $\barPrism_R
  \{ i \}$ for all $i$'s are uniformly $t$-bounded below. It follows that, the
  map
  \[ \tmop{gr}_M^{\ast} \tmop{THH} (R) \longrightarrow \tmop{gr}_M^{\ast}
     \tmop{THH} (R)^{t C_p} \]
  is an equivalence after $p$-completion as $\mathord{\ast} \gg 0$. The result
  then follows from the completeness of motivic filtrations on $\tmop{THH}$
  and $\tmop{THH}^{t C_p}$, {\cite[Thm~6.2.4 \& Prop~6.2.12]{Bhatt2022}} along
  with the trick of tensoring with $\mathbb{T}/ C_p$ as in
  {\cite[Lem~3.8]{Riggenbach2022}}.
\end{proof}

We give some special cases of \Cref{thm:mot-cpl-TP-poly,thm:Segal-conj-F-sm}.
We learned the following result from Daniel {\tmname{Fink}}.

\begin{lemma}
  \label{lem:abs-cot-cplx-retr-rel-frob}Let $k \rightarrow R$ be a map of
  animated $\mathbb{F}_p$-algebra. Then the $R$-module spectrum $L_{R / k}$
  can be realized, functorially in maps $(k \rightarrow R) \in \tmop{Fun}
  (\Delta^1, \tmop{CAlg}_{\mathbb{F}_p}^{\tmop{an}})$ of animated
  $\mathbb{F}_p$-algebras, as a retract of the cotangent complex $L_{R /
  \varphi_k^{\ast} R}$ of the relative Frobenius map $\varphi_{R / k} \of
  \varphi_k^{\ast} R \rightarrow R$, viewed as a module over the target.
\end{lemma}

\begin{proof}
  This is implicit in the proof of
  {\cite[Lem~2.19]{Blickle2025}}\footnote{This turns out to be implicit in
  {\cite[Lem~6.5.13(i)]{Gabber2003}} as well.}. We briefly review their
  argument. We have the transitivity sequence
  \[ L_{\varphi_k^{\ast} R / k} \otimes_{\varphi_k^{\ast} R, \varphi_{R /
     k}}^{\mathbb{L}} R \longrightarrow L_{R / k} \longrightarrow L_{R /
     \varphi_k^{\ast} R} \]
  associated to the composite map $k \rightarrow \varphi_k^{\ast} R
  \xrightarrow{\varphi_{R / k}} R$ of animated rings. It remains to specify a
  null-homotopy of the first map. Indeed, the first map is zero on standard
  maps $k \rightarrow R$ of the form $\mathbb{F}_p [x_1, \ldots, x_m]
  \rightarrow \mathbb{F}_p [x_1, \ldots, x_m, y_1, \ldots, y_n]$ for $(m, n)
  \in \mathbb{N}^2$, which gives rise to a null-homotopy of the first map
  functorially in such standard maps $k \rightarrow R$. We then left Kan
  extend this null homotopy to all maps $k \rightarrow R$ of animated
  $\mathbb{F}_p$-algebras.
\end{proof}

We introduce the following weakened version of $F$-finiteness.

\begin{definition}
  We say that an animated $\mathbb{F}_p$-algebra $R$ is {\tmdfn{weakly
  $F$-finite}} if the Frobenius map $\varphi \of R \rightarrow R$ realizes the
  target as a weakly finitely generated module over the source.
\end{definition}

This weakened $F$-finiteness gives rise to a weakened finite generation, just
as usual $F$-finiteness:

\begin{lemma}
  \label{lem:weak-ft-weak-fg-cot-cplx}Let $R \rightarrow S$ be a map of
  animated rings. Suppose that, as an $R$-modules, $S$ is weakly finitely
  generated. Then the cotangent complex $L_{S / R}$ is weakly finitely
  generated.
\end{lemma}

\begin{proof}
  Pick a surjection $P \rightarrow S$ where $P$ is a flat $R$-module of
  rank$\nosymbol \leq r$. We write $P$ as a filtered colimit $\tmop{colim}_{i
  \in \mathcal{I}} F_i$ of finite free $R$-modules of rank$\nosymbol \leq r$.
  Then we get a surjection
  \[ \underset{i \in \mathcal{I}}{\tmop{colim}} \tmop{LSym}_R (F_i)
     \longrightarrow S \]
  of animated $R$-algebras, which gives rise to a surjection
  \[ \underset{i \in \mathcal{I}}{\tmop{colim}} L_{\tmop{LSym}_R (F_i) / R}
     \otimes_{\tmop{LSym}_R (F_i)}^{\mathbb{L}} S \longrightarrow L_{S / R} \]
  of $S$-module spectra, where the source is a filtered colimit of finite free
  $S$-modules of rank$\nosymbol \leq r$. The result then follows.
\end{proof}

\begin{corollary}
  \label{cor:weak-F-fin-weak-fg-cot-cplx}Let $A$ be a weakly $F$-finite
  animated $\mathbb{F}_p$-algebra. Denote by $\varphi \of A \rightarrow A$ its
  Frobenius map. Then the cotangent complex $L_{A /\mathbb{F}_p}$ is weakly
  finitely generated.
\end{corollary}

\begin{proof}
  It follows from
  \Cref{lem:abs-cot-cplx-retr-rel-frob,lem:weak-ft-weak-fg-cot-cplx}.
\end{proof}

\begin{example}
  \label{ex:F-fin-mot-cpl}Let $R$ be an animated ring such that the animated
  $\mathbb{F}_p$-algebra $R \otimes_{\mathbb{Z}}^{\mathbb{L}} \mathbb{F}_p$ is
  weakly $F$-finite. It then follows from
  \Cref{cor:weak-F-fin-weak-fg-cot-cplx,lem:ft-ext-pow-unif-bnd-below,thm:mot-cpl-TP-poly}
  that the motivic filtration on $\tmop{TP}_{/\mathbb{Z}}^{\tmop{poly}} (R)$
  is complete.
\end{example}

\begin{example}
  \label{ex:F-fin-F-sm-Segal-conj}Let $R$ be an $F$-smooth $p$-quasisyntomic
  ring such that the animated ring $\mathbb{F}_p$-algebra $R
  \otimes_{\mathbb{Z}}^{\mathbb{L}} \mathbb{F}_p$ is weakly $F$-finite. By
  \Cref{cor:weak-F-fin-weak-fg-cot-cplx,lem:ft-ext-pow-unif-bnd-below,thm:Segal-conj-F-sm},
  we see that $\tmop{THH} (R)_p^{\wedge}$ satisfies Segal conjecture. This
  includes qualitative versions for both smooth algebras over a perfectoid
  ring {\cite[Prop~5.10]{Mathew2021}} and regular rings whose reduction
  $\tmop{mod} p$ is $F$-finite ({\cite[Cor~5.3]{Mathew2021}} in the excellent
  case, and {\cite[Cor~4.20]{Bhatt2023a}} in general).
\end{example}

\begin{example}
  We give an example of a Cartier smooth $\mathbb{F}_p$-algebra which is
  weakly $F$-finite but not $F$-finite. Let $A$ be a smooth
  $\mathbb{F}_p$-algebra, and $S \subseteq A$ a multiplicative subset. Then
  the subring $R \assign A + t (S^{- 1} A) [t] \subseteq (S^{- 1} A) [t]$ is
  Cartier smooth, as it is the filtered colimit $\tmop{colim}_{s \in S} A
  [s^{- 1} t]$. The $R$-module $L_{R /\mathbb{F}_p}$ is not finitely generated
  in general (say, $A =\mathbb{F}_p [x]$ and $S = \{ x^n \barsuchthat n \in
  \mathbb{N} \}$), but $R$ is weakly $F$-finite.
\end{example}

\begin{example}[non $F$-finite DVR]
  \label{ex:non-F-fin-DVR}Let $k$ be a perfect field of $\tmop{char} p$, and
  let $f = \sum_{n = 1}^{\infty} a_n t^n \in k \llbracket t \rrbracket$ be a
  formal power series with zero constant coefficient which is transcendental
  over $k [t]$. Let $f_j \assign \sum_{n = 1}^j a_n t^n$, and
  \begin{eqnarray*}
    u_j & \assign & \frac{(f - f_j)^p}{t^j}\\
    & = & \sum_{n = j + 1}^{\infty} a_n^p t^{np - j}
  \end{eqnarray*}
  for $j \in \mathbb{N}$. Then the localization $k [t, u_0, \ldots, u_j,
  \ldots]_{(t, u_0, \ldots, u_j, \ldots)} \subseteq k \llbracket t \rrbracket$
  is a DVR, but it is not excellent, as established in {\cite{Valabrega1973}}.
  However, it is weakly $F$-finite, being a filtered union of the subrings $k
  [t, u_0]_{(t, u_0)} \subseteq k [t, u_1]_{(t, u_1)} \subseteq k [t,
  u_2]_{(t, u_2)} \subseteq \cdots$ (note that we have $tu_{j + 1} - u_j \in
  tk \llbracket t \rrbracket$).
\end{example}

\begin{example}[non weakly $F$-finite DVR]
  We point out that there exists a DVR of char $p$ with perfect residue field,
  which is not weakly $F$-finite over $\mathbb{F}_p$. This is the sort of DVR
  constructed in {\cite[§4.1]{Datta2018}} (cf. {\cite[§2]{Li2026}}). Indeed,
  let $k \assign \overline{\mathbb{F}}_p$ and $K \assign k (x_1, \ldots, x_n,
  \ldots)$. There an embedding $\imath \of K \rightarrow k \pparen{t}$ of
  $k$-algebras obtained by assigning the image of $x_n$ recursively via a
  cardinality argument in {\cite[Rem~4.3]{Datta2018}}. Let $V_p \assign k
  \llbracket t \rrbracket \cap \iota (K)$ with $\mathfrak{m}_{V_p} \assign tk
  \llbracket t \rrbracket \cap \iota (K)$. Then $(V_p, \mathfrak{m}_{V_p})$ is
  a DVR over $k$ with residue field being $V_p /\mathfrak{m}_{V_p} = k$ and
  fractional field $\tmop{Frac} (V_p)$ being $\iota (K) \cong K$, as explained
  in {\cite[§4.1]{Datta2018}}. It follows that
  \[ L_{V_p /\mathbb{F}_p} \otimes_{V_p}^{\mathbb{L}} K = L_{K /\mathbb{F}_p}
     = \bigoplus_{n = 1}^{\infty} K \mathd x_i, \]
  and thus by \Cref{cor:flat-mod-weak-fg-fin-rk,lem:rk-flat-Bezout} and the
  flatness of the $V_p$-module $L_{V_p /\mathbb{F}_p}$ established in
  {\cite[Prop~9.5.1]{Bhatt2018a}}, the $V_p$-module $L_{V_p /\mathbb{F}_p}$ is
  not weakly finitely generated, and by
  \Cref{cor:weak-F-fin-weak-fg-cot-cplx}, the $\mathbb{F}_p$-algebra $V_p$ is
  not weakly $F$-finite over $\mathbb{F}_p$.
\end{example}

\section{Crystalline degeneration of Segal conjecture}

In {\cite{Hahn2022}}, the authors proved Segal conjecture for $\tmop{THH}
(\tmop{BP} \langle n \rangle)$. In this section, we show that their argument
implies a general statement (\Cref{prop:crys-deg-Segal-conj}) which allows us
to deduce Segal conjecture of a ``formal $\tmop{BP} \langle n
\rangle$-scheme'' from Segal conjecture for its special fiber over
$\mathbb{F}_p$. We refer to this argument as {\tmdfn{crystalline
degeneration}} by virtue of its resemblance to {\cite[§2.2]{Antieau2024}}.

First, the proof of {\cite[Prop~C.5.4]{Hahn2022a}} actually shows that

\begin{lemma}
  \label{lem:conv-ANSS}Let $A$ be a connective $\mathbb{E}_1$-ring. Then for
  every $n \in \mathbb{Z}$, the filtered object
  \[ \pi_n \left( \tmop{THH} \left( \tmop{desc}_{\mathbb{F}_p}^{\geqslant
     \mathord{\ast}} (A) \right) / v_0 \right) \]
  is eventually constant, converging to $\tmop{THH} (A)_p^{\wedge}$.
\end{lemma}

\begin{proposition}
  \label{prop:crys-deg-Segal-conj}Let $M$ be a right $\tmop{THH} (\tmop{BP}
  \langle n \rangle)$-module in genuine $C_p$-spectra, such that $M / (p, v_1,
  \ldots, v_n)$ is $t$-bounded below. Suppose that the canonical map
  \[ M^{\Phi C_p} \otimes_{\tmop{THH} (\tmop{BP} \langle n
     \rangle)}^{\mathbb{L}} \tmop{THH} (\mathbb{F}_p) \longrightarrow (M
     \otimes_{\tmop{THH} (\tmop{BP} \langle n \rangle)}^{\mathbb{L}}
     \tmop{THH} (\mathbb{F}_p))^{t C_p} \]
  is an equivalence on $\pi_{\ast}$ for $\mathord{\ast} \gg 0$. Then so is the
  map
  \[ M^{\Phi C_p} \otimes_{\tmop{BP} \langle n \rangle}^{\mathbb{L}}
     \mathbb{F}_p \longrightarrow M^{t C_p} \otimes_{\tmop{BP} \langle n
     \rangle}^{\mathbb{L}} \mathbb{F}_p \]
  induced by the canonical map $M^{\Phi C_p} \rightarrow M^{t C_p}$.
\end{proposition}

\begin{proof}
  We endow $\tmop{BP} \langle n \rangle$ the $\mathbb{N}$-graded decreasing
  filtration $\tmop{desc}_{\mathbb{F}_p}^{\geqslant \mathord{\ast}} (\tmop{BP}
  \langle n \rangle)$ coming from the Adams spectral sequence as in
  {\cite[Rec~4.3.2]{Hahn2022}}, with
  \[ \pi_{\ast} \left( \tmop{gr} \left( \tmop{desc}_{\mathbb{F}_p}^{\geqslant
     \mathord{\ast}} (\tmop{BP} \langle n \rangle) \right) \right)
     =\mathbb{F}_p [v_0, \ldots, v_n], \]
  where $\tmop{wt} (v_i) = 2 p^i - 1$. Then the cyclotomic spectrum
  $\tmop{THH} (\tmop{BP} \langle n \rangle)$ is equipped with an
  $\mathbb{N}$-graded decreasing filtration $\tmop{Fil}_{\tmop{cris}}^{\ast}$
  as well, concretely given by
  \begin{eqnarray*}
    \tmop{Fil}_{\tmop{cris}}^{\ast} \tmop{THH} (\tmop{BP} \langle n \rangle) &
    = & \tmop{THH} \left( \tmop{desc}_{\mathbb{F}_p}^{\geqslant
    \mathord{\ast}} (\tmop{BP} \langle n \rangle) \right) .
  \end{eqnarray*}
  This induces an $\mathbb{N}$-graded decreasing filtration
  $\tmop{Fil}_{\tmop{cris}}^{\ast}$ on the canonical map $M^{\Phi C_p}
  \rightarrow M^{t C_p}$, and it suffices to see that the filtered map (for
  convenience, we omit the grading shear $L_p$ on the source)
  \begin{equation}
    \tmop{Fil}_{\tmop{cris}}^{\ast} M^{\Phi C_p}
    \otimes_{\tmop{desc}_{\mathbb{F}_p}^{\geqslant \mathord{\ast}} (\tmop{BP}
    \langle n \rangle)}^{\mathbb{L}} \mathbb{F}_p \longrightarrow
    \tmop{Fil}_{\tmop{cris}}^{\ast} M^{t C_p}
    \otimes_{\tmop{desc}_{\mathbb{F}_p}^{\geqslant \mathord{\ast}} (\tmop{BP}
    \langle n \rangle)}^{\mathbb{L}} \mathbb{F}_p \label{eq:fil-frob-mod-v}
  \end{equation}
  is an equivalence in high enough degrees.
  
  We first show this after passing to associated graded pieces. For this, we
  do note need the bounded-belowness of the genuine $C_p$-spectrum $M$.
  Indeed, let $M_{\mathbb{F}_p} \assign M \otimes_{\tmop{THH} (\tmop{BP}
  \langle n \rangle)}^{\mathbb{L}} \tmop{THH} (\mathbb{F}_p)$. The associated
  graded pieces of the map \eqref{eq:fil-frob-mod-v} is the composite graded
  map
  \begin{eqnarray*}
    &  & M_{\mathbb{F}_p}^{\Phi C_p} \otimes_{\tmop{THH}
    (\mathbb{F}_p)}^{\mathbb{L}} \tmop{THH} (\mathbb{F}_p [v_0, \ldots, v_n])
    / (v_0, \ldots, v_n)\\
    & \longrightarrow & M_{\mathbb{F}_p}^{t C_p} \otimes_{\tmop{THH}
    (\mathbb{F}_p)}^{\mathbb{L}} \tmop{THH} (\mathbb{F}_p [v_0, \ldots, v_n])
    / (v_0, \ldots, v_n)\\
    & \longrightarrow & (M_{\mathbb{F}_p} \otimes_{\tmop{THH}
    (\mathbb{F}_p)}^{\mathbb{L}} \tmop{THH} (\mathbb{F}_p [v_0, \ldots,
    v_n]))^{t C_p} / (v_0, \ldots, v_n) .
  \end{eqnarray*}
  The first map is an equivalence in high enough degrees, as it is tensoring
  the canonical map $M_{\mathbb{F}_p}^{\Phi C_p} \rightarrow
  M_{\mathbb{F}_p}^{t C_p}$, which is an equivalence in high enough degrees by
  assumption, with $\tmop{THH} (\mathbb{S} [v_0, \ldots, v_n]) / (v_0, \ldots,
  v_n)$. The second map is an equivalence by {\cite[Lem~4.1.3]{Hahn2022}} (see
  the last paragraph of the proof of {\cite[Prop~4.2.2]{Hahn2022}}).
  
  We now show the completeness of filtrations on the source and the target of
  the filtered map \eqref{eq:fil-frob-mod-v}. By \Cref{lem:conv-ANSS} and that
  the genuine $C_p$-spectrum $M$ is bounded below, the filtration
  $\tmop{Fil}_{\tmop{cris}}^{\ast}$ is $(p, v_1, \ldots, v_n)$-completely
  eventually constant on $M^{\Phi C_p}$ and $M$, thus the filtration
  $\tmop{Fil}_{\tmop{cris}}^{\ast}$ is complete on $M^{t C_p}$ by
  {\cite[Cor~C.4.3]{Hahn2022}}. The result then follows from the equivalence
  \[ \tmop{desc}_{\mathbb{F}_p}^{\geqslant \mathord{\ast}} (\tmop{BP} \langle
     n \rangle) / (v_0, \tilde{v}_1, \ldots, \tilde{v}_n) \simeq \mathbb{F}_p
  \]
  as in the proof of {\cite[Thm~4.3.1]{Hahn2022}}.
\end{proof}

\begin{example}
  Let $R$ be a $(- 1)$-connective $\mathbb{E}_1$-$\mathbb{Z}$-algebra. Suppose
  that the cyclotomic spectrum $\tmop{THH} (R
  \otimes_{\mathbb{Z}}^{\mathbb{L}} \mathbb{F}_p)$ satisfies Segal conjecture.
  Then \Cref{prop:crys-deg-Segal-conj} implies that, the cyclotomic spectrum
  $\tmop{THH} (R)$ satisfies Segal conjecture $p$-completely. This can be
  viewed as a shadow of the fact that, a $p$-quasisyntomic $p$-torsion-free
  commutative ring $R$ is $F$-smooth if the commutative $\mathbb{F}_p$-algebra
  $R / p$ is Cartier smooth, {\cite[Cor~4.17]{Bhatt2023a}}.
\end{example}

\section{Relative conjugate filtration}\label{sec:rel-conj-fil}

In this section, we will present a general procedure
(Construction~\ref{cons:KO-THH-decpl-tCp}), which allows us to produce a
relative conjugate filtration (\Cref{ex:KO-Hdg-Tate}) on Hodge--Tate
cohomology $\barPrism_R$ for any map $k \rightarrow R$ of animated rings, a
mixed characteristic analogue of {\cite[Prop~3.22]{Bhatt2012a}}. In the same
time, it also allows us to produce a relative conjugate filtration
(\Cref{ex:KO-THH-Phi-Cp}) on both $\tmop{THH}_{\Prism} (\mathcal{C})^{\Phi
C_p}$ and $\tmop{THH} (\mathcal{C}) \otimes_{\tmop{THH} (k)}^{\mathbb{L}}
\tmop{THH} (k)^{t C_p}$ after $p$-completion, for every animated ring $k$ and
every $k$-linear DG-category $\mathcal{C}$, mixed characteristic and relative
versions of Kaledin's conjugate filtration. We also discuss the completeness
of these filtrations (\Cref{prop:KO-cpl-SynSp-hat,prop:KO-cpl-Frob-twist}).

{\construction{\label{cons:KO-THH-decpl-tCp}To simplify notation, everything
is $p$-completed. Let $k$ be an animated ring. The theory of synthetic
$\tmop{THH}$ as in {\cite{Antieau2024a}} gives rise to an
$\mathbb{E}_{\infty}$-algebra
\[ \tmop{Fil}_M^{\geqslant 0} (\tmop{THH} (k) \otimes_{\tmop{THH}
   (\mathbb{Z})}^{\mathbb{L}} \tmop{THH} (\mathbb{Z})^{t C_p}) =
   \tmop{Fil}_M^{\geqslant 0} \tmop{THH}_{\Prism} (k)^{\Phi C_p} \]
in $\tmop{SynSp}$. This gives rise to a $\mathbb{N}$-filtered
$\mathbb{E}_{\infty}$-algebra\footnote{Along with the ``internal'' filtration
of synthetic spectra, this can be seen as a bi-filtered object, but we prefer
to view it as a filtered synthetic spectrum. In addition, the subscript
$\tmop{KO}$ is an abbreviation for Katz--Oda, as the procedure here is the
same as the procedure of producing Katz--Oda filtration on derived de Rham
cohomology.} $\tmop{Fil}_{\tmop{KO}}^{\ast} \tmop{Fil}_M^{\geqslant 0}
\tmop{THH}_{\Prism} (k)^{\Phi C_p}$ in $\tmop{SynSp}$, concretely given by
\[ \tmop{Fil}_{\tmop{KO}}^j \tmop{Fil}_M^{\geqslant 0} \tmop{THH}_{\Prism}
   (k)^{\Phi C_p} \assign \tmop{Fil}_M^{\geqslant j} \tmop{THH}_{\Prism}
   (k)^{\Phi C_p} \in \tmop{SynSp} \]
for $j \in \mathbb{N}$, where given a filtered object $\tmop{Fil} W$, the
filtered object of $\tmop{Fil}^{\geqslant j} W$ is
\[ \cdots \xleftarrow{\simeq} \tmop{Fil}_M^j W \xleftarrow{\simeq}
   \tmop{Fil}_M^j W \leftarrow \tmop{Fil}_M^{j + 1} W \leftarrow \cdots, \]
that is, first restricting to $\mathbb{N}_{\geqslant j}$ then left Kan
extending along $(\mathbb{N}_{\geqslant j}, \geq) \rightarrow (\mathbb{N},
\geq)$.

Let $N$ be a $\tmop{Fil}_M^{\geqslant 0} \tmop{THH}_{\Prism} (k)^{\Phi
C_p}$-module in $\tmop{SynSp}$. The $\mathbb{N}$-graded filtration on
$\tmop{Fil}_M^{\geqslant 0} \tmop{THH}_{\Prism} (k)^{\Phi C_p}$ gives rise to
an $\mathbb{Z}$-graded filtration on
\[ N' \assign N \otimes_{\tmop{Fil}_M^{\geqslant 0} \tmop{THH}_{\Prism}
   (k)^{\Phi C_p}}^{\mathbb{L}} \tmop{Fil}_M \tmop{THH}_{\Prism} (k)^{\Phi
   C_p}, \]
given by
\[ \tmop{Fil}_{\tmop{KO}}^j N' \assign N \otimes_{\tmop{Fil}_M^{\geqslant 0}
   \tmop{THH}_{\Prism} (k)^{\Phi C_p}}^{\mathbb{L}} \tmop{Fil}_M^{\geqslant j}
   \tmop{THH}_{\Prism} (k)^{\Phi C_p} \in \tmop{SynSp} \]
for $j \in \mathbb{Z}$. We see that
\begin{eqnarray*}
  \tmop{gr}_{\tmop{KO}}^j N' & = & N \otimes_{\tmop{Fil}_M^{\geqslant 0}
  \tmop{THH}_{\Prism} (k)^{\Phi C_p}}^{\mathbb{L}} \tmop{ins}^j \tmop{gr}_M^j
  \tmop{THH} (k)^{t C_p}\\
  & = & N \otimes_{\tmop{Fil}_M^{\geqslant 0} \tmop{THH}_{\Prism} (k)^{\Phi
  C_p}}^{\mathbb{L}} \tmop{ins}^j  \barPrism_k \{ j \} [2 j] \in \tmop{SynSp},
\end{eqnarray*}
where $\tmop{ins}^j$ is the left Kan extension along $\{ j \} \hookrightarrow
(\mathbb{N}, \geq)$, see {\cite[Not~3.1.1]{Raksit2020}}. Moreover, unwinding
the construction, we see that the $\tmop{Fil}_M^{\geqslant 0}
\tmop{THH}_{\Prism} (k)^{\Phi C_p}$-module structure on $\tmop{ins}^j 
\barPrism_k \{ j \} [2 j]$ is induced by the $\tmop{gr}_M^0
\tmop{THH}_{\Prism} (k)^{\Phi C_p} = \barPrism_k$-module structure on
$\barPrism_k \{ j \} [2 j]$.

We may replace $\tmop{THH}_{\Prism} (k)^{\Phi C_p}$ by $\tmop{THH} (k)^{t
C_p}$ throughout, and for every $\tmop{Fil}_M^{\geqslant 0} \tmop{THH} (k)^{t
C_p}$-module $N$ in $\tmop{SynSp}$, we get an $\mathbb{Z}$-graded filtration
on
\[ N' \assign N \otimes_{\tmop{Fil}_M^{\geqslant 0} \tmop{THH} (k)^{t
   C_p}}^{\mathbb{L}} \tmop{Fil}_M \tmop{THH} (k)^{t C_p}, \]
with associated graded pieces
\[ \tmop{gr}_{\tmop{KO}}^j N' = N \otimes_{\tmop{Fil}_M^{\geqslant 0}
   \tmop{THH} (k)^{t C_p}}^{\mathbb{L}} \tmop{ins}^j  \hatbarPrism_k \{ j \}
   [2 j] \in \tmop{SynSp} \]
for $j \in \mathbb{Z}$.}}

\begin{example}
  \label{ex:KO-Frob-twist-THH-k-mod}To simplify notation, everything is
  $p$-completed. Let $k$ be an animated ring, and $M$ a $\tmop{Fil}_M
  \tmop{THH} (k)$-module in $\tmop{SynSp}$. Apply
  Construction~\ref{cons:KO-THH-decpl-tCp} to the $\tmop{Fil}_M
  \tmop{THH}_{\Prism} (k)^{\Phi C_p}$-module
  \[ M \otimes_{\tmop{Fil}_M \tmop{THH} (k)}^{\mathbb{L}}
     \tmop{Fil}_M^{\geqslant 0} \tmop{THH}_{\Prism} (k)^{\Phi C_p} \]
  in $\tmop{SynSp}$, we get a filtered synthetic spectrum
  \[ \tmop{Fil}_{\tmop{KO}}^{\ast} \left( M \otimes_{\tmop{Fil}_M \tmop{THH}
     (k)}^{\mathbb{L}} \tmop{Fil}_M \tmop{THH}_{\Prism} (k)^{\Phi C_p} \right)
  \]
  with $j$-th associated graded pieces given by
  \begin{eqnarray*}
    &  & M \otimes_{\tmop{Fil}_M \tmop{THH} (k)}^{\mathbb{L}} \tmop{ins}^j 
    \barPrism_k \{ j \} [2 j]\\
    & = & (M \otimes_{\tmop{Fil}_M \tmop{THH} (k)}^{\mathbb{L}} \tmop{ins}^0
    \tmop{gr}_M^0 \tmop{THH} (k)) \otimes_{\tmop{ins}^0 \tmop{gr}_M^0
    \tmop{THH} (k)}^{\mathbb{L}} \tmop{ins}^j  \barPrism_k \{ j \} [2 j]\\
    & = & (M \otimes_{\tmop{Fil}_M \tmop{THH} (k)}^{\mathbb{L}} \tmop{ins}^0
    k) \otimes_{\tmop{ins}^0 k}^{\mathbb{L}} \tmop{ins}^j  \barPrism_k \{ j \}
    [2 j],
  \end{eqnarray*}
  for $j \in \mathbb{Z}$, where in the second equality, we take advantage of
  the commutative diagram
  \[ \begin{array}{ccc}
       \tmop{Fil}_M \tmop{THH} (k) & \longrightarrow & \tmop{Fil}_M
       \tmop{THH}_{\Prism} (k)^{\Phi C_p}\\
       \longdownarrow &  & \longdownarrow\\
       \tmop{ins}^0 \tmop{gr}_M^0 \tmop{THH} (k) & \longrightarrow &
       \tmop{ins}^0 \tmop{gr}_M^0 \tmop{THH}_{\Prism} (k)^{\Phi C_p}
     \end{array} \]
  of $\mathbb{E}_{\infty}$-algebras of synthetic spectra. Unwinding the
  construction, we see that the $\tmop{ins}^0 k$-module structure on
  $\tmop{ins}^j  \barPrism_k \{ j \} [2 j]$ is given by the canonical map
  $\varphi_0 \of k \rightarrow \barPrism_k$. Similarly, we get a filtered
  synthetic spectrum
  \[ \tmop{Fil}_{\tmop{KO}}^{\ast} (M \otimes_{\tmop{Fil}_M \tmop{THH}
     (k)}^{\mathbb{L}} \tmop{Fil}_M \tmop{THH} (k)^{t C_p}) \]
  with $j$-th associated graded piece given by
  \[ (M \otimes_{\tmop{Fil}_M \tmop{THH} (k)}^{\mathbb{L}} \tmop{ins}^0 k)
     \otimes_{\tmop{ins}^0 k}^{\mathbb{L}} \tmop{ins}^j  \hatbarPrism_k \{ j
     \} [2 j] \]
  for $j \in \mathbb{Z}$.
\end{example}

\begin{notation}
  Let $k \rightarrow R$ be a map of animated ring, and $\mathcal{C}$ a
  dualizable presentable stable $k$-linear $\infty$-category. We will denote
  by $\tmop{THH}_{\Prism / k} (\mathcal{C})^{\Phi C_p}$ the
  $\mathbb{T}$-equivariant spectrum
  \[ \tmop{THH} (\mathcal{C}) \otimes_{\tmop{THH} (k)}^{\mathbb{L}} \tmop{THH}
     (k)^{t C_p}, \]
  and by $\tmop{Fil}_M \tmop{THH}_{\Prism / k} (R)^{\Phi C_p}$ the
  $\mathbb{T}_{\tmop{ev}}$-equivariant synthetic spectrum
  \[ \tmop{Fil}_M \tmop{THH} (R) \otimes_{\tmop{Fil}_M \tmop{THH}
     (k)}^{\mathbb{L}} \tmop{Fil}_M \tmop{THH} (k)^{t C_p} . \]
\end{notation}

\begin{example}
  \label{ex:KO-Hdg-Tate}To simplify notation, everything is $p$-completed. Let
  $k \rightarrow R$ be a map of animated rings. Take $M = \tmop{Fil}_M
  \tmop{THH} (R)$ in \Cref{ex:KO-Frob-twist-THH-k-mod}, we get a filtered
  synthetic spectrum\footnote{We replace the subscript $\tmop{KO}$ by
  $\tmop{KO} / k$ to stress that the base is $k$.} $\tmop{Fil}_{\tmop{KO} / k}
  \tmop{Fil}_M \tmop{THH}_{\Prism} (R)^{\Phi C_p}$, with
  \begin{eqnarray*}
    \tmop{gr}_{\tmop{KO} / k}^j \tmop{Fil}_M \tmop{THH}_{\Prism} (R)^{\Phi
    C_p} & = & (\tmop{Fil}_M \tmop{THH} (R) \otimes_{\tmop{Fil}_M \tmop{THH}
    (k)}^{\mathbb{L}} \tmop{ins}^0 k) \otimes_{\tmop{ins}^0 k}^{\mathbb{L}}
    \tmop{ins}^j  \barPrism_k \{ j \} [2 j]\\
    & = & \tmop{Fil}_M \tmop{HH} (R / k) \otimes_{\tmop{ins}^0
    k}^{\mathbb{L}} \tmop{ins}^j  \barPrism_k \{ j \} [2 j] \in \tmop{SynSp},
  \end{eqnarray*}
  and a filtered synthetic spectrum $\tmop{Fil}_{\tmop{KO} / k} \tmop{Fil}_M
  \tmop{THH}_{\Prism / k} (R)^{\Phi C_p}$, with
  \[ \tmop{gr}_{\tmop{KO} / k}^j \tmop{Fil}_M \tmop{THH}_{\Prism / k}
     (R)^{\Phi C_p} = \tmop{Fil}_M \tmop{HH} (R / k) \otimes_{\tmop{ins}^0
     k}^{\mathbb{L}} \tmop{ins}^j  \hatbarPrism_k \{ j \} [2 j] \in
     \tmop{SynSp} . \]
  Passing to the $i$th associated graded piece of the synthetic spectra, we
  get filtered objects $\tmop{Fil}_{\tmop{KO} / k}  \barPrism_R \{ i \} [2
  i]$, with
  \begin{eqnarray*}
    \tmop{gr}_{\tmop{KO} / k}^j  \barPrism_R \{ i \} [2 i] & = & \left(
    \bigwedgestar_R^{i - j} L_{R / k} \right) [i - j] \otimes_k^{\mathbb{L}}
    \barPrism_k \{ j \} [2 j]\\
    & = & \left( \bigwedgestar_R^{i - j} L_{R / k} \right)
    \otimes_k^{\mathbb{L}} \barPrism_k \{ j \} [i + j],\\
    \tmop{gr}_{\tmop{KO} / k}^j  \barPrism_R \{ i \} & = & \left(
    \bigwedgestar_R^{i - j} L_{R / k} \right) \otimes_k^{\mathbb{L}}
    \barPrism_k \{ j \} [j - i] .
  \end{eqnarray*}
  We note that this filtration is concentrated in degrees $j \leq i$, and in
  the usual convention, one might invert the sign of the weight to make it an
  increasing filtration. One also sees that this (inverted) increasing
  filtration starts with zero, which is a special case of
  \Cref{prop:KO-cpl-SynSp-hat}. When $k$ is an animated
  $\mathbb{F}_p$-algebra, this recovers {\cite[Prop~3.22]{Bhatt2012a}} for the
  composite map $\mathbb{F}_p \rightarrow k \rightarrow R$ of animated rings,
  thus we view our result as a mixed characteristic analogue of that.
\end{example}

\begin{example}
  \label{ex:KO-THH-Phi-Cp}To simplify notation, everything is $p$-completed.
  Let $k$ be an animated ring, and $\mathcal{C}$ a dualizable presentable
  $k$-linear stable $\infty$-category. Take $M$ to be $\tmop{THH}
  (\mathcal{C})$, with induced synthetic spectrum structure from $\tmop{Fil}_M
  \tmop{THH} (k)$, in \Cref{ex:KO-Frob-twist-THH-k-mod}, we get a filtered
  synthetic spectrum, whose underlying filtered spectrum under forgetful
  functor $\tmop{SynSp} \rightarrow \tmop{Sp}$ is denoted by
  $\tmop{Fil}_{\tmop{KO} / k} \tmop{THH}_{\Prism} (\mathcal{C})^{\Phi C_p}$,
  with
  \[ \tmop{gr}_{\tmop{KO} / k}^j \tmop{THH}_{\Prism} (\mathcal{C})^{\Phi C_p}
     = \tmop{HH} (\mathcal{C}/ k) \otimes_k^{\mathbb{L}} \barPrism_k \{ j \}
     [2 j], \]
  and another filtered spectrum $\tmop{Fil}_{\tmop{KO} / k} \tmop{THH}_{\Prism
  / k} (\mathcal{C})^{\Phi C_p}$, with
  \[ \tmop{gr}_{\tmop{KO} / k}^j \tmop{THH}_{\Prism / k} (\mathcal{C})^{\Phi
     C_p} = \tmop{HH} (\mathcal{C}/ k) \otimes_k^{\mathbb{L}} \hatbarPrism_k
     \{ j \} [2 j] . \]
\end{example}

We make some final remarks on completeness of Katz--Oda filtration, which is,
by construction, exhaustive. The case is direct if we consider $\tau$-complete
synthetic spectra instead.

\begin{notation}
  Let $\widehat{\tmop{SynSp}} \subseteq \tmop{SynSp}$ denote the full
  subcategory of $\tau$-complete synthetic spectra, which admits a left
  adjoint $\widehat{(-)} \of \tmop{SynSp} \rightarrow \widehat{\tmop{SynSp}}$.
\end{notation}

\begin{proposition}
  \label{prop:KO-cpl-SynSp-hat}Let $k$ be an animated ring, and $N$ a
  $\tmop{Fil}_M^{\geqslant 0} \tmop{THH}_{\Prism} (k)^{\Phi C_p}$-module in
  $\tmop{SynSp}$ which is weight-bounded below. Then the filtered synthetic
  spectra $\tmop{Fil}_{\tmop{KO}} N'$ in both cases of
  Construction~\ref{cons:KO-THH-decpl-tCp} is complete after $\tau$-completion
  to $\widehat{\tmop{SynSp}}$.
\end{proposition}

\begin{proof}
  Since the functor $\tmop{gr} \of \widehat{\tmop{SynSp}} \rightarrow
  \tmop{Mod}_{\tmop{gr} (\mathbb{S}_{\tmop{ev}})} (\tmop{Gr} (\tmop{Sp}))$
  preserves small colimits and is conservative, it suffice to in
  $\tmop{Mod}_{\tmop{gr} (\mathbb{S}_{\tmop{ev}})} (\tmop{Gr} (\tmop{Sp}))$.
  In this case, since $N$ is weight-bounded below, one can see that the limit
  tower $(\tmop{Fil}_{\tmop{KO}}^n N')_n$, after passing to $\tmop{Gr}
  (\tmop{Sp})$, is degreewise eventually zero.
\end{proof}

Now we address the motivic completeness on the nose, which is equivalent to
the motivic completeness of the underlying spectrum if the motivic
completeness holds after $\tau$-completion.

\begin{proposition}
  \label{prop:KO-cpl-Frob-twist}Let $k$ be an animated ring, and $M$ a
  $\tmop{Fil}_M \tmop{THH} (k)$-module in $\tmop{SynSp}$ which is
  weight-bounded below. Suppose that $M$ is uniformly $t$-bounded below.
  \begin{enumerate}
    \item If the spectra $\left\{ \tmop{gr}_N^m  \Prism_k \{ n \} \barsuchthat
    (m, n) \in \mathbb{Z} \times \mathbb{Z} \right\}$ is uniformly $t$-bounded
    below (this is the case when the $k$-module $L_{k /\mathbb{Z}}$ is
    $p$-completely weakly finitely generated by
    \Cref{cor:Nygaard-gr-unif-bnd}), then the Katz--Oda filtration on
    \[ M \otimes_{\tmop{Fil}_M \tmop{THH} (k)}^{\mathbb{L}} \tmop{Fil}_M
       \tmop{THH}_{\Prism} (k)^{\Phi C_p} \]
    in \Cref{ex:KO-Frob-twist-THH-k-mod} is complete;
    
    \item If the Nygaard-completed Hodge--Tate $\hatbarPrism \{ n \}$ is
    $t$-bounded below uniformly in $n \in \mathbb{Z}$ (this is the case when
    the $k$-module $L_{k /\mathbb{Z}}$ is $p$-completely weakly finitely
    generated by \Cref{cor:Nygaard-gr-unif-bnd}), then the Katz--Oda
    filtration on
    \[ M \otimes_{\tmop{Fil}_M \tmop{THH} (k)}^{\mathbb{L}} \tmop{Fil}_M
       \tmop{THH} (k)^{t C_p} \]
    in \Cref{ex:KO-Frob-twist-THH-k-mod} is complete.
  \end{enumerate}
\end{proposition}

\begin{proof}
  Both cases follows from connectivity estimates of $\tmop{Fil}_M
  \tmop{THH}_{\Prism} (k)^{\Phi C_p}$ and $\tmop{Fil}_M \tmop{THH} (k)^{t
  C_p}$. Indeed, in the first case, connectivity estimates are established in
  the proof of \Cref{thm:mot-cpl-TP-poly}, while the second connectivity
  estimates follow from the motivic completeness of $\tmop{THH} (k)^{t C_p}$
  for animated rings $k$, cf.~{\cite[Cor~6.2.15]{Bhatt2022}} along with the
  trick of tensoring with $\mathbb{T}/ C_p$ as in
  {\cite[Lem~3.8]{Riggenbach2022}}. 
\end{proof}

\section{Transitivity of $F$-smoothness}

The goal of this section is to establish transitivity of $F$-smoothness. For
this purpose, we first introduce a variant of conjugate filtration in
\Cref{sec:rel-conj-fil} on associate graded pieces of Nygaard-filtered
absolute prismatic cohomology, which leads to a transitivity of weak
$F$-smoothness (\Cref{thm:weak-F-sm-quasism}). Then we discuss a relative
Nygaard-completeness (\Cref{def:rel-Nyg-cpl}). We give an animated version of
transitivity of $p$-Cartier smoothness (the first case of
\Cref{thm:transitivity-Cart-sm}), and a transitivity of $F$-smoothness
(\Cref{thm:dR-Hdg-cpl-rel-Nyg-cpl,cor:transitivity-F-sm}). The strategy of two
proofs are parallel:
\begin{enumerate}
  \item Establish a flatness of the map between cohomology theories via
  conjugate filtrations
  (\Cref{lem:dR-rel-conj-fil,lem:Tor-ampl-prism-coh-quasi-syn}).
  
  \item Using this flatness to deduce completeness of Katz--Oda filtrations
  (\Cref{cor:quasi-lci-bnd-Tor-ampl-KO-cpl,cor:KO-tw-Nyg-cpl}).
  
  \item Analyze associated graded pieces of Katz--Oda filtrations.
\end{enumerate}
Moreover, we define {\tmdfn{universal convergence}} of Hodge and Nygaard
filtration (\Cref{def:Hdg-filt-univ-conv,def:Nyg-filt-univ-conv}) generalizing
``having a universal Cartier isomorphism'' {\cite[Def~9.5.9]{Bhatt2018a}}.
This notion gives rise to variants of transitivity of $p$-Cartier smoothness
and $F$-smoothness without bounded $\tmop{Tor}$-amplitude (the first case of
\Cref{thm:transitivity-Cart-sm} and
\Cref{thm:transitivity-F-sm-base-conv-univ}). We also establish transitivity
of universal convergence
(\Cref{thm:transitivity-univ-conv-Hdg-filt,thm:transitivity-univ-conv-Nyg-filt}).
Over a perfectoid base, two universal convergences could also be compared
(\Cref{thm:Hdg-Nyg-conv-univ,prop:Hdg-Nyg-conv-univ}).

{\construction{\label{cons:KO-THH}To simplify notation, everything is
$p$-completed. Let $k$ be an animated ring, and $M$ a $\tmop{Fil}_M \tmop{THH}
(k)$-module in $\tmop{SynSp}$. As in Construction~\ref{cons:KO-THH-decpl-tCp}
and \Cref{ex:KO-Frob-twist-THH-k-mod}, the filtration given by
$\tmop{Fil}_{\tmop{KO}}^j \tmop{Fil}_M \tmop{THH} (k) \assign
\tmop{Fil}_M^{\geqslant j} \tmop{THH} (k)$ on $\tmop{Fil}_M \tmop{THH} (k)$
gives rise to a filtration $\tmop{Fil}_{\tmop{KO}} M$ on $M$, with
\begin{eqnarray*}
  \tmop{gr}_{\tmop{KO}}^j M & = & M \otimes_{\tmop{Fil}_M \tmop{THH}
  (k)}^{\mathbb{L}} \tmop{ins}^j \tmop{gr}_M^j \tmop{THH} (k)\\
  & = & (M \otimes_{\tmop{Fil}_M \tmop{THH} (k)}^{\mathbb{L}} \tmop{ins}^0 k)
  \otimes_{\tmop{ins}^0 k}^{\mathbb{L}} \tmop{ins}^j \tmop{gr}_N^j  \Prism_k
  \{ j \} [2 j] \in \tmop{Mod}_{\tmop{Fil}_M \tmop{THH} (k)} (\tmop{SynSp}) .
\end{eqnarray*}}}

\begin{example}
  \label{ex:KO-THH}To simplify notation, everything is $p$-completed. Let $k$
  be an animated ring, and $\mathcal{C}$ a dualizable presentable $k$-linear
  stable $\infty$-category. Take $M$ to be $\tmop{THH} (\mathcal{C})$, with
  induced synthetic spectrum structure from $\tmop{Fil}_M \tmop{THH} (k)$, in
  \Cref{ex:KO-Frob-twist-THH-k-mod}, we get a filtered synthetic spectrum,
  whose underlying filtered spectrum under forgetful functor $\tmop{SynSp}
  \rightarrow \tmop{Sp}$ is denoted by $\tmop{Fil}_{\tmop{KO} / k} \tmop{THH}
  (\mathcal{C})$, with
  \[ \tmop{gr}_{\tmop{KO} / k}^j \tmop{THH} (\mathcal{C}) = \tmop{HH}
     (\mathcal{C}/ k) \otimes_k^{\mathbb{L}} \tmop{gr}_N^j  \Prism_k \{ j \}
     [2 j] . \]
\end{example}

\begin{example}
  \label{ex:KO-Nygaard-gr}To simplify notation, everything is $p$-completed.
  Let $k \rightarrow R$ be a map of animated rings. Then as in
  \Cref{ex:KO-Hdg-Tate}, applying Construction~\ref{cons:KO-THH} to $M =
  \tmop{Fil}_M \tmop{THH} (R)$, we get a filtered synthetic spectrum
  $\tmop{Fil}_{\tmop{KO} / k} \tmop{Fil}_M \tmop{THH} (R)$ which admits a
  commutative ring structure, with
  \begin{eqnarray*}
    \tmop{gr}_{\tmop{KO} / k}^j \tmop{Fil}_M \tmop{THH} (R) & = &
    (\tmop{Fil}_M \tmop{THH} (R) \otimes_{\tmop{Fil}_M \tmop{THH}
    (k)}^{\mathbb{L}} \tmop{ins}^0 k) \otimes_{\tmop{ins}^0 k}^{\mathbb{L}}
    \tmop{ins}^j \tmop{gr}_N^j  \Prism_k \{ j \} [2 j]\\
    & = & \tmop{Fil}_M \tmop{HH} (R / k) \otimes_{\tmop{ins}^0
    k}^{\mathbb{L}} \tmop{ins}^j \tmop{gr}_N^j  \Prism_k \{ j \} [2 j] .
  \end{eqnarray*}
  Passing to the $i$th associated graded pieces of the synthetic spectra, we
  get filtered objects $\tmop{Fil}_{\tmop{KO} / k} \tmop{gr}_N^i  \Prism_k \{
  i \} [2 i]$, with
  \begin{eqnarray*}
    \tmop{gr}_{\tmop{KO} / k}^j \tmop{gr}_N^i  \Prism_k \{ i \} [2 i] & = &
    \left( \bigwedgestar_R^{i - j} L_{R / k} \right) [i - j]
    \otimes_k^{\mathbb{L}} \tmop{gr}_N^j  \Prism_k \{ j \} [2 j]\\
    & = & \left( \bigwedgestar_R^{i - j} L_{R / k} \right)
    \otimes_k^{\mathbb{L}} \tmop{gr}_N^j  \Prism_k \{ j \} [i + j] \in D
    (R),\\
    \tmop{gr}_{\tmop{KO} / k}^j \tmop{gr}_N^i  \Prism_k \{ i \} & = & \left(
    \bigwedgestar_R^{i - j} L_{R / k} \right) \otimes_k^{\mathbb{L}}
    \tmop{gr}_N^j  \Prism_k \{ j \} [j - i] \in D (R),
  \end{eqnarray*}
  which is concentrated in weights $j \leq i$, where the $R$-module structure
  is induced by that on the first tensor factor, thus the filtration is
  complete and exhaustive, a special case of \Cref{prop:KO-cpl-THH}.
\end{example}

Again, by construction, Katz--Oda filtration on $\tmop{THH}$ is always
exhaustive. Completeness of Katz--Oda filtration on $\tmop{THH}$ is much
better than that on $\tmop{THH}_{\Prism}^{\Phi C_p}$ and $\tmop{THH}_{\Prism /
k}^{\Phi C_p}$:

\begin{proposition}
  \label{prop:KO-cpl-THH}Let $k$ be an animated ring, and $M$ a $\tmop{Fil}_M
  \tmop{THH} (k)$-module in $\tmop{SynSp}$ which is weight-bounded below, and
  $\tau$-complete. Then the filtered synthetic spectra $\tmop{Fil}_{\tmop{KO}}
  M$ in Construction~\ref{cons:KO-THH} is complete.
\end{proposition}

\begin{proof}
  First, each synthetic spectrum $\tmop{Fil}_{\tmop{KO}}^j \tmop{Fil}_M
  \tmop{THH} (k)$ is complete, as the maps $\tmop{Fil}_{\tmop{KO}}^{j + 1} M
  \rightarrow \tmop{Fil}_{\tmop{KO}}^j M$ induces an equivalence on the limit
  of the tower of multiplication by $\tau$, thus the limits of these towers
  are the same, being zero. Thus it suffices to check the completeness of
  $\tmop{Fil}_{\tmop{KO}} M$ after $\tau$-completion to
  $\widehat{\tmop{SynSp}}$, and the proof of \Cref{prop:KO-cpl-SynSp-hat}
  adapts.
\end{proof}

\begin{example}
  \label{ex:KO-Nygaard-gr-monogen-poly}Let $k$ be an animated ring. Take $R$
  to be the monogenic polynomial $k$-algebra $k [t]$. Then
  \Cref{ex:KO-Nygaard-gr} gives rise to a 2-step increasing filtration, which
  corresponds to a fiber sequence
  \[ L_{k [t] / k}^1 \otimes_k^{\mathbb{L}} \tmop{gr}_N^{i - 1}  \Prism_k [-
     1] \longrightarrow \tmop{gr}_N^i  \Prism_{k [t]} \longrightarrow k [t]
     \otimes_k^{\mathbb{L}} \tmop{gr}_N^i  \Prism_k, \]
  which is the fiber sequence in {\cite[Prop~4.8]{Bhatt2023a}}, used to
  establish the weak $F$-smoothness of $k [t]$ from that of $k$.
\end{example}

\Cref{ex:KO-Nygaard-gr-monogen-poly} hints a possibility of applying
\Cref{ex:KO-Nygaard-gr} to understand weak $F$-smoothness. This is indeed the
case (\Cref{thm:weak-F-sm-quasism}). We first review weak $F$-smoothness in
{\cite[Def~4.2]{Bhatt2023a}}. We give a slightly different definition, which
is equivalent to theirs by \Cref{rem:compare-weak-F-sm-BM}.

\begin{definition}
  \label{def:weak-F-sm}We say that an animated ring $A$ is {\tmdfn{weakly
  $F$-smooth}} if it is $p$-quasi-lci over $\mathbb{Z}$, and for each $i \in
  \mathbb{Z}$, the $(A \otimes_{\mathbb{Z}}^{\mathbb{L}} \mathbb{F}_p)$-module
  spectrum
  \[ \tmop{fib} \left( \mu_1 \of \tmop{gr}_N^i  \Prism_A \{ j \}
     \otimes_{\mathbb{Z}}^{\mathbb{L}} \mathbb{F}_p \longrightarrow
     \tmop{gr}_N^{i + p}  \Prism_A \{ j + p \}
     \otimes_{\mathbb{Z}}^{\mathbb{L}} \mathbb{F}_p \right) \]
  has $\tmop{Tor}$-amplitude in degrees$\nosymbol \leq - (i + 2)$ for
  every\footnote{This is equivalent to requiring the same for a single $j \in
  \mathbb{Z}$, as Breuil--Kisin twists trivialize after passing to associated
  graded pieces of the Nygaard filtration by {\cite[Rem~5.5.15]{Bhatt2022}}.}
  $j \in \mathbb{Z}$, where $\mu_1$ is induced by the class $\mu_1$ in
  $\tmop{Fil}_M (\tmop{THH} (\mathbb{Z}) \otimes_{\mathbb{Z}}^{\mathbb{L}}
  \mathbb{F}_p)$, which gives rise to a map
  \[ \tmop{Fil}_M (\tmop{THH} (A) \otimes_{\mathbb{Z}}^{\mathbb{L}}
     \mathbb{F}_p) \longrightarrow \tmop{Fil}_M (\tmop{THH} (A)
     \otimes_A^{\mathbb{L}} \mathbb{F}_p) (- p) [- 2 p] \]
  of synthetic spectra.
\end{definition}

\begin{remark}
  Let $A$ be an animated ring which is $p$-quasi-lci over $\mathbb{Z}$. Then
  for every pair $(i, j) \in \mathbb{Z}^2$ with $i < p$, by
  {\cite[Cor~5.5.18]{Bhatt2022}}, we have $\tmop{gr}_N^i  \Prism_A \{ j \}
  \otimes_{\mathbb{Z}}^{\mathbb{L}} \mathbb{F}_p \simeq \left(
  \bigwedgestar_A^i L_{A /\mathbb{Z}} \right)
  \otimes_{\mathbb{Z}}^{\mathbb{L}} \mathbb{F}_p [- i]$, which is zero when $i
  < 0$. Since $A$ is $p$-quasi-lci over $\mathbb{Z}$, the $A
  \otimes_{\mathbb{Z}}^{\mathbb{L}} \mathbb{F}_p$-module $\tmop{gr}_N^i 
  \Prism_A \{ j \} \otimes_{\mathbb{Z}}^{\mathbb{L}} \mathbb{F}_p$ has
  $\tmop{Tor}$-amplitude$\nosymbol \leq 0$. It follows that the condition on
  $\tmop{fib} (\mu_1)$ in \Cref{def:weak-F-sm} is automatic when $i < 0$, and
  moreover, the $A \otimes_{\mathbb{Z}}^{\mathbb{L}} \mathbb{F}_p$-module
  $\tmop{fib} (\mu_1)$ has $\tmop{Tor}$-amplitude$\nosymbol < 0$ for every
  pair $(i, j) \in \mathbb{Z}^2$.
\end{remark}

In {\cite[Def~4.2]{Bhatt2023a}}, the authors considered the $p$-complete
$\tmop{Tor}$-amplitude of $\tmop{fib} (\mu_1)$ as an $A$-module spectrum. To
see that it is equivalent\footnote{In an earlier version of this article, we
mistakenly deduced this equivalence from the converse direction of
\Cref{lem:Tor-ampl-mod-I} without assuming that the $A / I$-module $M$ is
extended from an $A$-module $M_0$, cf.~\Cref{ex:mod-I-tor-ampl-infty}.} to
ours, we need a simple lemma for $\tmop{Tor}$-amplitude of $A / I$-modules as
$A$-modules for an invertible ideal $I \subseteq A$.

\begin{lemma}
  \label{lem:Tor-ampl-mod-I}Let $A$ be an animated ring, and $I \rightarrow A$
  a generalized Cartier divisor, and $M$ an $A / I$-module. Let $n \in
  \mathbb{Z}$ be an integer. If the $A / I$-module $M$
  $\tmop{Tor}$-amplitude$\nosymbol \leq n$, then the $A$-module $M$ has
  $\tmop{Tor}$-amplitude$\nosymbol \leq n + 1$. The converse is true if there
  exists an $A$-module $M_0$ such that $M = M_0 \otimes_A^{\mathbb{L}} (A /
  I)$.
\end{lemma}

\begin{proof}
  Suppose that the $A / I$-module $M$ has $\tmop{Tor}$-amplitude$\nosymbol
  \leq n$. For every coconnective $A$-module $N$, as the spectrum $(A / I)
  \otimes_A^{\mathbb{L}} N$ is $1$-truncated, the spectrum
  \[ M \otimes_A^{\mathbb{L}} N \simeq M \otimes_{A / I}^{\mathbb{L}} ((A / I)
     \otimes_A^{\mathbb{L}} N) \]
  is $(n + 1)$-truncated.
  
  On the other hand, suppose that $M = M_0 \otimes_A^{\mathbb{L}} (A / I)$ and
  the $A$-module $M$ has $\tmop{Tor}$-amplitude$\nosymbol \leq n + 1$. Then
  the multiplication map $m \of M \otimes_A^{\mathbb{L}} (A / I) \rightarrow
  M$ admits an $A / I$-linear retract. This implies that the $A / I$-module
  $\tmop{fib} (m) = M \otimes_{A / I}^{\mathbb{L}} (I / I^2) [1]$ has
  $\tmop{Tor}$-amplitude$\nosymbol \leq n + 1$. As the $A / I$-module $I /
  I^2$ is connectively invertible, the result follows.
\end{proof}

\begin{example}
  \label{ex:mod-I-tor-ampl-infty}The $\mathbb{Z}$-module $\mathbb{F}_p$ has
  $\tmop{Tor}$-amplitude$\nosymbol \leq 1$, while as a $\mathbb{Z}/
  p^2$-module, it has unbounded $\tmop{Tor}$-amplitude, since $\mathbb{F}_p
  \otimes_{\mathbb{Z}/ p^2}^{\mathbb{L}} \mathbb{F}_p \simeq \bigoplus_{n \in
  \mathbb{N}} \mathbb{F}_p [n]$. This shows that the converse of
  \Cref{lem:Tor-ampl-mod-I} is false in general.
\end{example}

\begin{corollary}
  \label{cor:A-mod-I-mods-I-cpl-Tor-ampl}Let $A$ be an animated ring, and $I
  \rightarrow A$ a generalized Cartier divisor, and $M$ an $A / I$-module. Let
  $n \in \mathbb{Z}$ be an integer. Then the $A$-module $M$ has $I$-complete
  $\tmop{Tor}$-amplitude$\nosymbol \leq n$ if and only if the $A$-module $M$
  has $\tmop{Tor}$-amplitude$\nosymbol \leq n$.
\end{corollary}

\begin{proof}
  ``if'' part follows from definitions. We suppose that the $A$-module $M$ has
  $I$-complete $\tmop{Tor}$-amplitude$\nosymbol \leq n$. By definition, the $A
  / I$-module $M \otimes_A^{\mathbb{L}} (A / I)$ has
  $\tmop{Tor}$-amplitude$\nosymbol \leq n$. It follows from
  \Cref{lem:Tor-ampl-mod-I} that the $A$-module $M \otimes_A^{\mathbb{L}} (A /
  I)$ has $\tmop{Tor}$-amplitude$\nosymbol \leq n + 1$. As in the proof of the
  converse direction of \Cref{lem:Tor-ampl-mod-I}, the multiplication map $m
  \of M \otimes_A^{\mathbb{L}} (A / I) \rightarrow M$ admits an $A$-linear
  retract $M \rightarrow M \otimes_A^{\mathbb{L}} (A / I)$ induced by $A
  \rightarrow A / I$, thus the $A$-module $\tmop{fib} (m) = M \otimes_A I [1]$
  has $\tmop{Tor}$-amplitude$\nosymbol \leq n + 1$. As the $A$-module $I$ is
  connectively invertible, the result follows.
\end{proof}

\begin{remark}
  \label{rem:compare-weak-F-sm-BM}Let $A$ be an animated ring which is
  $p$-quasi-lci over $\mathbb{Z}$. Then the proof of
  {\cite[Prop~4.3]{Bhatt2023a}} shows that $A$ is weak $F$-smooth in our sense
  if and only if, for every $i \in \mathbb{Z}$, the $A
  \otimes_{\mathbb{Z}}^{\mathbb{L}} \mathbb{F}_p$-module spectrum
  \[ \tmop{fib} \left( \varphi_i \of \tmop{gr}_N^i  \Prism_A \rightarrow
     \barPrism_A \{ i \} \right) \otimes_{\mathbb{Z}}^{\mathbb{L}}
     \mathbb{F}_p \]
  has $\tmop{Tor}$-amplitude$\nosymbol \leq - (i + 2)$, which is equivalent to
  that the $A$-module spectrum $\tmop{fib} (\varphi_i)
  \otimes_{\mathbb{Z}}^{\mathbb{L}} \mathbb{F}_p$ has
  $\tmop{Tor}$-amplitude$\nosymbol \leq - (i + 1)$ by
  \Cref{lem:Tor-ampl-mod-I}, which is equivalent to that the $A$-module
  spectrum $\tmop{fib} (\varphi_i) \otimes_{\mathbb{Z}}^{\mathbb{L}}
  \mathbb{F}_p$ has $p$-completed $\tmop{Tor}$-amplitude$\nosymbol \leq - (i +
  1)$ by \Cref{cor:A-mod-I-mods-I-cpl-Tor-ampl}. This is equivalent to weak
  $F$-smoothness in {\cite[Def~4.2]{Bhatt2023a}}.
\end{remark}

\begin{theorem}
  \label{thm:weak-F-sm-quasism}Let $k \rightarrow R$ be a map of animated
  rings. Suppose that
  \begin{itemize}
    \item the animated ring $k$ is weakly $F$-smooth; and that
    
    \item the $R$-module spectrum $L_{R / k}$ is $p$-completely flat.
  \end{itemize}
  Then the animated ring $R$ is weakly $F$-smooth as well.
\end{theorem}

\begin{proof}
  First, as the animated ring $k$ is $p$-quasi-lci, and the $R$-module $L_{R /
  k}$ is $p$-completely flat, the animated ring $R$ is $p$-quasi-lci as well,
  by transitivity of cotangent complex. For animated rings $A$, let $F^i (A)$
  denote the fiber
  \[ \tmop{fib} \left( \mu_1 \of \tmop{gr}_N^i  \Prism_A \{ i \}
     \otimes_{\mathbb{Z}}^{\mathbb{L}} \mathbb{F}_p \longrightarrow
     \tmop{gr}_N^{i + p}  \Prism_A \{ i + p \}
     \otimes_{\mathbb{Z}}^{\mathbb{L}} \mathbb{F}_p \right) \]
  being an $A$-module. We want to invoke \Cref{ex:KO-Nygaard-gr}, but we have
  to be a bit more careful, as the map $\mu_1$ involves a shift of Katz--Oda
  weight, internal filtration weight, and homotopy degree. Indeed, $\mu_1$
  comes from a class in $\tmop{Fil}_M (\tmop{THH} (\mathbb{Z})
  \otimes_{\mathbb{Z}}^{\mathbb{L}} \mathbb{F}_p)$, which has filtration
  weight $p$ and homotopy degree $2 p$. This gives rise to a class in
  $\tmop{Fil}_{\tmop{KO}} \tmop{Fil}_M (\tmop{THH} (k)
  \otimes_{\mathbb{Z}}^{\mathbb{L}} \mathbb{F}_p)$, which has Katz--Oda weight
  $p$, internal filtration weight $p$, and homotopy degree $2 p$. In the end,
  we indeed get a Katz--Oda filtration $\tmop{Fil}_{\tmop{KO}} F^i (R)$ on the
  $R$-module $F^i (R)$, with
  \[ \tmop{gr}_{\tmop{KO}}^j F^i (R) = \left( \bigwedgestar_R^{i - j} L_{R /
     k} \right) \otimes_k^{\mathbb{L}} F^j (k) [j - i] \in D (R), \]
  which corresponds to a finite-step exhaustive increasing filtration. Now for
  every $j$, the $R$-module spectrum
  \[ \left( \bigwedgestar_R^{i - j} L_{R / k} \right) \otimes_k^{\mathbb{L}}
     F^j (k) [j - i] \]
  has $p$-complete $\tmop{Tor}$-amplitude$\nosymbol \leq - (j + 2) + (j - i) =
  - (i + 2)$, as for every $(R \otimes_{\mathbb{Z}}^{\mathbb{L}}
  \mathbb{F}_p)$-module $N$ which is concentrated in degree $0$, the
  $k$-module $N \otimes_R^{\mathbb{L}} \left( \bigwedgestar_R^{i - j} L_{R /
  k} \right)$ is concentrated in degree $0$ by $p$-complete flatness of the
  $R$-module $L_{R / k}$. The result then follows.
\end{proof}

Now we prove transitivity of $F$-smoothness. For this, it seems reasonable to
introduce a relative version of Nygaard-completeness, which is precisely
relative Segal conjecture up to completions with respect to motivic
filtrations.

\begin{definition}
  \label{def:rel-Nyg-cpl}Let $k \rightarrow R$ be a map of animated rings. We
  say that the absolute prismatic cohomology $\Prism_R \{ \ast \}$ is
  {\tmdfn{Nygaard-complete relative to $\Prism_k \{ \ast \}$}} if the
  canonical map
  \[ \tmop{gr}_M \tmop{THH} (R) \otimes_{\tmop{gr}_M \tmop{THH}
     (k)}^{\mathbb{L}} \tmop{gr}_M \tmop{THH} (k)^{t C_p} \longrightarrow
     \tmop{gr}_M \tmop{THH} (R)^{t C_p} \]
  of graded $\mathbb{E}_{\infty}$-rings is an equivalence after
  $p$-completion, or equivalently\footnote{We use the symmetric monoidal
  structure on the shearing, which is explained in
  {\cite[Prop~3.3.4]{Raksit2020}}.}, the canonical map
  \[ \left( \bigoplus_{n \in \mathbb{Z}} \barPrism_R \{ n \} \right)
     \otimes_{\bigoplus_{n \in \mathbb{Z}} \barPrism_k \{ n \}}^{\mathbb{L}}
     \left( \bigoplus_{n \in \mathbb{Z}} \hatbarPrism_k \{ n \} \right)
     \longrightarrow \bigoplus_{n \in \mathbb{Z}} \hatbarPrism_R \{ n \} \]
  of graded $\mathbb{E}_{\infty}$-rings is an equivalence after
  $p$-completion.
\end{definition}

\begin{remark}
  \label{rem:Nyg-rel-cpl}Let $k \rightarrow R$ be a map of animated rings.
  Then by definition, the absolute prismatic cohomology $\Prism_R \{ \ast \}$
  is Nygaard-complete relative to $\Prism_k \{ \ast \}$ if and only if the
  filtered graded $\mathbb{E}_{\infty}$-ring
  \[ \left( \bigoplus_{n \in \mathbb{Z}} \tmop{Fil}_N  \barPrism_R \{ n \}
     \right) \otimes_{\bigoplus_{n \in \mathbb{Z}} \tmop{Fil}_N  \barPrism_k
     \{ n \}}^{\mathbb{L}} \left( \bigoplus_{n \in \mathbb{Z}} \tmop{Fil}_N 
     \hatbarPrism_k \{ n \} \right) \]
  is complete up to $p$-completion.
\end{remark}

The rest of this section is devoted to the following theorem.

\begin{theorem}
  \label{thm:dR-Hdg-cpl-rel-Nyg-cpl}Let $k \rightarrow R$ be a map of animated
  rings such that $k$ is $p$-completely $t$-bounded\footnote{We say that a
  spectrum $X \in \tmop{Sp}$ is {\tmdfn{$p$-completely $t$-bounded}} if the
  spectrum $X / p$ is $t$-bounded.} and $R$ has bounded $p$-complete
  $\tmop{Tor}$-amplitude over $k$. Suppose that
  \begin{itemize}
    \item the animated ring $k$ is weakly $F$-smooth; and that
    
    \item the $R$-module spectrum $L_{R / k}$ is $p$-completely flat.
  \end{itemize}
  If the $p$-completed derived de Rham cohomology $\tmop{dR}_{R / k}$ is
  $p$-completely Hodge-complete, then the absolute prismatic cohomology
  $\Prism_R \{ \ast \}$ is Nygaard-complete relative to $\Prism_k \{ \ast \}$.
\end{theorem}

Before its proof, we mention that it implies transitivity of $F$-smoothness:

\begin{corollary}
  \label{cor:transitivity-F-sm}Let $k$ be an $F$-smooth $p$-quasisyntomic
  ring, and $k \rightarrow R$ a $p$-Cartier smooth map {\cite{Bouis2023}} such
  that $R$ has bounded $p^{\infty}$-torsion\footnote{We do not know whether
  this assumption is redundant.}, and that $R$ has bounded $p$-complete
  $\tmop{Tor}$-amplitude over $k$. Then $R$ is $F$-smooth $p$-quasisyntomic as
  well.
\end{corollary}

\begin{proof}
  We examine the transitivity sequence
  \[ L_{k /\mathbb{Z}_p} \otimes_k^{\mathbb{L}} R \longrightarrow L_{R
     /\mathbb{Z}_p} \longrightarrow L_{R / k} . \]
  Since $k$ is $p$-quasisyntomic, the $k$-module $L_{k /\mathbb{Z}_p}$ has
  $p$-complete $\tmop{Tor}$-amplitude in $[0, 1]$. By assumption, the
  $R$-module $L_{R / k}$ is $p$-completely flat\footnote{The $p$-discreteness
  of $k \rightarrow R$ implies that the canonical map $(R
  \otimes_{\mathbb{Z}}^{\mathbb{L}} \mathbb{F}_p) \otimes_{k
  \otimes_{\mathbb{Z}}^{\mathbb{L}} \mathbb{F}_p}^{\mathbb{L}} \pi_0 (k
  \otimes_{\mathbb{Z}}^{\mathbb{L}} \mathbb{F}_p) \rightarrow \pi_0 (R
  \otimes_{\mathbb{Z}}^{\mathbb{L}} \mathbb{F}_p)$ is an equivalence.}, thus
  $L_{R /\mathbb{Z}_p}$ has $p$-complete $\tmop{Tor}$-amplitude in $[0, 1]$ as
  well, thus $R$ is $p$-quasisyntomic. Since the ring $k$ is $F$-smooth, it is
  weakly $F$-smooth {\cite[Prop~4.3]{Bhatt2023a}}. By
  \Cref{thm:weak-F-sm-quasism}, the ring $R$ is weakly $F$-smooth as well. Now
  the result follows from \Cref{thm:dR-Hdg-cpl-rel-Nyg-cpl} and the
  Nygaard-completeness of the absolute prismatic cohomology $\Prism_k$.
\end{proof}

Now we come to the proof of \Cref{thm:dR-Hdg-cpl-rel-Nyg-cpl}. Let $k
\rightarrow R$ be a map of animated rings. The strategy is similar to
\Cref{sec:rel-conj-fil}, namely introducing a filtration on the graded
$\mathbb{E}_{\infty}$-ring $\bigoplus_{n \in \mathbb{Z}} \hatbarPrism_k \{
\ast \}$, showing that it induces a complete ``Katz--Oda'' filtration on the
map
\[ \left( \bigoplus_{n \in \mathbb{Z}} \barPrism_R \{ n \} \right)
   \otimes_{\bigoplus_{n \in \mathbb{Z}} \barPrism_k \{ n \}}^{\mathbb{L}}
   \left( \bigoplus_{n \in \mathbb{Z}} \hatbarPrism_k \{ n \} \right)
   \longrightarrow \bigoplus_{n \in \mathbb{Z}} \hatbarPrism_R \{ n \} \]
and then show that, on each ``Katz--Oda'' associated graded pieces, the
filtered graded spectrum is complete. A naive guess is to take the Nygaard
filtration on $\bigoplus_{n \in \mathbb{Z}} \hatbarPrism_k \{ n \}$ as in
Construction~\ref{cons:KO-THH-decpl-tCp} etc. Tensor products do not preserve
completeness, and such a ``Katz--Oda'' filtration is a priori not necessarily
complete. Fortunately, our assumptions allow us to control the convergence.

We start with a closer analysis of convergence of the Hodge filtration on de
Rham cohomology, and to see how to deduce transitivity of $p$-Cartier
smoothness via the usual Katz--Oda filtration. As a toy example, for a map $k
\rightarrow R$ of commutative rings with the $R$-module $L_{R / k}$ being
flat, the $r$-th filtered piece of the Hodge-completed derived de Rham
cohomology $\widehat{\tmop{dR}}_{R / k}$ is computed by the cochain complex
\[ \cdots \longrightarrow 0 \longrightarrow \Omega_{R / k}^r
   \xrightarrow{\mathd} \Omega_{R / k}^{r + 1} \xrightarrow{\mathd} \cdots, \]
which is $(- r)$-truncated. That is to say, the convergence of the
Hodge-completed derived de Rham cohomology $\widehat{\tmop{dR}}_{R / k}$ has a
linearly increasing coconnectivity estimate. This extends to the animated
setting:

\begin{notation}
  Let $k$ be an animated $\mathbb{F}_p$-algebra. We denote by $\varphi_k$ the
  Frobenius map $k \rightarrow k$, and by $\varphi_k^{\ast}$ the (derived)
  base change $(-) \otimes_{k, \varphi_k}^{\mathbb{L}} k \of D (k) \rightarrow
  D (k)$.
\end{notation}

\begin{lemma}
  \label{lem:lin-conv-Hdg-fil-quasi-sm}Let $k \rightarrow R$ be a map of
  animated rings with $R$ being $n$-truncated for some $n \in \mathbb{N}$.
  Suppose that the $R$-module $L_{R / k}$ is flat. Then for every $r \in
  \mathbb{N}$, the $r$-th Hodge-filtered piece $\tmop{Fil}_H^r 
  \widehat{\tmop{dR}}_{R / k}$ is $(n - r)$-truncated.
\end{lemma}

\begin{proof}
  Note that, for every $s \in \mathbb{N}_{> r}$, the spectrum
  $\tmop{gr}_H^{[r, s)} \tmop{dR}_{R / k}$ is a finite extension of
  \[ \tmop{gr}_H^i \tmop{dR}_{R / k} = \left( \bigwedgestar_R^i L_{R / k}
     \right) [- i] \]
  for $i = r, r + 1, \ldots, s - 1$, and for such $i$'s, since the $R$-module
  $L_{R / k}$ is flat, so is the exterior power $\bigwedgestar_R^i L_{R / k}$,
  and thus $\tmop{gr}_H^i \tmop{dR}_{R / k}$ is $(n - i)$-truncated, thus the
  spectrum $\tmop{gr}_H^{[r, s)} \tmop{dR}_{R / k}$ is $(n - r)$-truncated.
  Taking $s \rightarrow \infty$, we see that the spectrum $\tmop{Fil}_H^r 
  \widehat{\tmop{dR}}_{R / k}$ is $(n - r)$-truncated as well.
\end{proof}

As a consequence of \Cref{lem:lin-conv-Hdg-fil-quasi-sm} and definition of
$\tmop{Tor}$-amplitude, we have

\begin{corollary}
  \label{cor:lin-conv-bnd-Tor-ampl-Hdg-fil-quasi-sm}Let $k \rightarrow R$ be a
  map of animated $\mathbb{F}_p$-algebras with $R$ being $n$-truncated for
  some $n \in \mathbb{N}$, and $m \in \mathbb{Z}$ an integer. Suppose that the
  $R$-module $L_{R / k}$ is flat. Then for every $\varphi_k^{\ast} R$-module
  $M$ of $\tmop{Tor}$-amplitude$\nosymbol \leq m$, and every $r \in
  \mathbb{N}$, the $r$-th Hodge-filtered piece
  \[ M \otimes_{\varphi_k^{\ast} R}^{\mathbb{L}} \tmop{Fil}_H^r 
     \widehat{\tmop{dR}}_{R / k} \]
  is $(m + n - r)$-truncated. In particular, the filtered spectrum $M
  \otimes_{\varphi_k^{\ast} R}^{\mathbb{L}} \tmop{Fil}_H 
  \widehat{\tmop{dR}}_{R / k}$ is complete.
\end{corollary}

Sometimes, we have an efficient bound of the $\tmop{Tor}$-amplitude of
$\tmop{Fil}_H^r  \widehat{\tmop{dR}}_{R / k}$, in which case a variant of
\Cref{cor:lin-conv-bnd-Tor-ampl-Hdg-fil-quasi-sm} holds. This notion is
partially inspired by universal Cartier isomorphisms in
{\cite[Def~9.5.14]{Bhatt2018a}} (cf.~\Cref{lem:univ-Cart-iso-Hdg-conv-univ}).

\begin{definition}
  \label{def:Hdg-filt-univ-conv}Let $k \rightarrow R$ be a map of animated
  $\mathbb{F}_p$-algebras. We say that a filtered $\tmop{Fil}_H \tmop{dR}_{R /
  k}$-module $\tmop{Fil} M$ {\tmdfn{converges universally (of index $\alpha
  \in \mathbb{Z}$)}} if, for every $r \in \mathbb{N}$, the $\varphi_k^{\ast}
  R$-module $\tmop{Fil}^r M$ has $\tmop{Tor}$-amplitude$\nosymbol \leq \alpha
  - r$, and it converges universally if such an integer $\alpha$ exists. In
  particular, we say that the Hodge filtration on $\widehat{\tmop{dR}}_{R /
  k}$ {\tmdfn{converges universally (resp. of index $\alpha \in \mathbb{Z}$)}}
  if the filtered $\tmop{Fil}_H \tmop{dR}_{R / k}$-module $\tmop{Fil}_H 
  \widehat{\tmop{dR}}_{R / k}$ converges universally (resp. of index $\alpha
  \in \mathbb{Z}$).
\end{definition}

It follows directly from definitions that

\begin{corollary}
  \label{cor:Hdg-conv-univ-trunc-Tor-ampl}Let $k \rightarrow R$ be a map of
  animated $\mathbb{F}_p$-algebras, and $N$ a filtered $\tmop{Fil}_H
  \tmop{dR}_{R / k}$-module which converges universally of index $\alpha \in
  \mathbb{Z}$. Then for $\varphi_k^{\ast} R$-module $M$, and every $r \in
  \mathbb{N}$,
  \begin{itemize}
    \item if $M$ is $m$-truncated for some $m \in \mathbb{Z}$, then the $r$-th
    Hodge-filtered piece
    \[ M \otimes_{\varphi_k^{\ast} R}^{\mathbb{L}} \tmop{Fil}^r N \]
    is $(\alpha + m - r)$-truncated. In particular, the filtered spectrum $M
    \otimes_{\varphi_k^{\ast} R}^{\mathbb{L}} \tmop{Fil} N$ is complete; and
    
    \item if the $\varphi_k^{\ast} R$-module $M$ has
    $\tmop{Tor}$-amplitude$\nosymbol \leq m$ for some $m \in \mathbb{Z}$, then
    the $r$-th Hodge-filtered piece
    \[ M \otimes_{\varphi_k^{\ast} R}^{\mathbb{L}} \tmop{Fil}^r N \]
    has $\tmop{Tor}$-amplitude$\nosymbol \leq \alpha + m - r$.
  \end{itemize}
\end{corollary}

Universal convergence of Hodge filtration implies flatness of Frobenius.

\begin{lemma}
  \label{lem:conv-univ-frob-flat}Let $k \rightarrow R$ be a map of animated
  $\mathbb{F}_p$-algebras, and $\tmop{Fil} M$ a filtered $\tmop{Fil}_H
  \tmop{dR}_{R / k}$-module equipped with a map $\tmop{Fil}_H \tmop{dR}_{R /
  k} \rightarrow \tmop{Fil} M$ which induces an equivalence on associated
  graded pieces. Suppose that the $\tmop{Fil}_H \tmop{dR}_{R / k}$-module
  $\tmop{Fil} M$ converges universally of index $\alpha \in \mathbb{Z}$. Then
  the Frobenius map $\varphi_{R / k} \of \varphi_k^{\ast} R \rightarrow R$ has
  $\tmop{Tor}$-amplitude$\nosymbol \leq \alpha$, and for every $r \in
  \mathbb{N}$, the $\varphi_k^{\ast} R$-module $\bigwedgestar_R^r L_{R / k}$
  (via forgetting along $\varphi_{R / k}$) has
  $\tmop{Tor}$-amplitude$\nosymbol \leq \alpha$.
\end{lemma}

\begin{proof}
  By assumption, for every $r \in \mathbb{N}$, the $\varphi_k^{\ast}
  R$-modules $\tmop{Fil}^r M [r]$ and $\tmop{Fil}^{r + 1} M [r + 1]$ have
  $\tmop{Tor}$-amplitude$\nosymbol \leq \alpha$. Now the result follows from
  the fiber sequence
  \[ \tmop{Fil}^r M [r] \longrightarrow \tmop{gr}^r M [r] \longrightarrow
     \tmop{Fil}^{r + 1} M [r + 1] \]
  of $\varphi_k^{\ast} R$-modules, and the equivalence $\bigwedgestar_R^r L_{R
  / k} = \tmop{gr}_H^r \tmop{dR}_{R / k} [r] \xrightarrow{\simeq} \tmop{gr}^r
  M [r]$ of $R$-modules (thus as $\varphi_k^{\ast} R$-modules).
\end{proof}

The Hodge-filtration on the Hodge-completed derived de Rham cohomology
converges universally for Cartier smooth $\mathbb{F}_p$-algebras with
universal Cartier isomorphism {\cite[Def~9.5.14]{Bhatt2018a}}.

\begin{lemma}
  \label{lem:univ-Cart-iso-Hdg-conv-univ}Let $R$ be a Cartier smooth
  $k$-algebra where $k \assign \mathbb{F}_p$. Suppose that, for every $r \in
  \mathbb{N}$, the $\varphi_k^{\ast} R$\footnote{This is just $R$. It is
  written as $\varphi_k^{\ast} R$ to stress that the $\varphi_k^{\ast}
  R$-module structures on $\Omega_{R / k}^{\ast}$ is Frobenius twisted, so
  that the differentials $\mathd \of \Omega_{R / k}^{\ast} \rightarrow
  \Omega_{R / k}^{\mathord{\ast} + 1}$ are $\varphi_k^{\ast}
  R$-linear.}-module $B \Omega_{R / k}^r \assign \tmop{Im} \left( \Omega_{R /
  k}^{r - 1} \xrightarrow{\mathd} \Omega_{R / k}^r \right)$ is flat. Then the
  Hodge filtration on $\widehat{\tmop{dR}}_{R / k}$ converges universally of
  index $0$.
\end{lemma}

\begin{proof}
  Indeed, for every $r \in \mathbb{N}$, the $\varphi_k^{\ast} R$-module
  $\tmop{Fil}_H^r  \widehat{\tmop{dR}}_{R / k}$ is represented by the cochain
  complex
  \[ \cdots \longrightarrow 0 \longrightarrow \Omega_{R / k}^r
     \xrightarrow{\mathd} \Omega_{R / k}^{r + 1} \xrightarrow{\mathd} \cdots
  \]
  of $\varphi_k^{\ast} R$-modules. We get a short exact sequence
  \[ \begin{array}{ccccccccc}
       \cdots & \longrightarrow & 0 & \longrightarrow & B \Omega_{R / k}^r &
       \longrightarrow & 0 & \longrightarrow & \cdots\\
       &  & \longdownequal &  & \longdownarrow &  & \longdownarrow &  & \\
       \cdots & \longrightarrow & 0 & \longrightarrow & \Omega_{R / k}^r &
       \xrightarrow{\mathd} & \Omega_{R / k}^{r + 1} & \xrightarrow{\mathd} &
       \cdots\\
       &  & \longdownequal &  & \longdownarrow &  & \longdownequal &  & \\
       \cdots & \longrightarrow & 0 & \longrightarrow & \Omega_{R / k}^r / B
       \Omega_{R / k}^r & \xrightarrow{\mathd} & \Omega_{R / k}^{r + 1} &
       \xrightarrow{\mathd} & \cdots
     \end{array} \]
  of cochain complexes of $\varphi_k^{\ast} R$-modules. As the $R$-module
  $L_{R / k}$ is flat, by Cartier isomorphism, we see that $\varphi_k^{\ast}
  R$-module represented by the bottom cochain complex has
  $\tmop{Tor}$-amplitude$\nosymbol \leq - r$. It follows from flatness of the
  $\varphi_k^{\ast} R$-module $B \Omega_{R / k}^r$ that the $\varphi_k^{\ast}
  R$-module $\tmop{Fil}_H^r  \widehat{\tmop{dR}}_{R / k}$ has
  $\tmop{Tor}$-amplitude$\nosymbol \leq - r$ as well.
\end{proof}

\begin{example}
  Let $R$ be a regular Noetherian $k$-algebra where $k \assign \mathbb{F}_p$.
  By {\cite[Thm~9.5.1~\&~9.5.18]{Bhatt2018a}}, the $k$-algebra $R$ is Cartier
  smooth, and has a universal Cartier isomorphism. It follows from
  \Cref{lem:univ-Cart-iso-Hdg-conv-univ} that the Hodge filtration on
  $\widehat{\tmop{dR}}_{R / k}$ converges universally of index $0$.
\end{example}

\begin{example}
  Let $R$ be a perfect $k$-algebra where $k \assign \mathbb{F}_p$. By
  {\cite[Prop~9.5.11]{Bhatt2018a}}, the $k$-algebra $R$ is Cartier smooth, and
  has a universal Cartier isomorphism. It follows from
  \Cref{lem:univ-Cart-iso-Hdg-conv-univ} that the Hodge filtration on
  $\widehat{\tmop{dR}}_{R / k}$ converges universally of index $0$.
\end{example}

This actually admits an animated generalization:

\begin{lemma}
  \label{lem:rel-perf-Hdg-filt-cpl-conv-univ}Let $k \rightarrow R$ be a
  relative perfect map of animated $\mathbb{F}_p$-algebras. Then the
  Hodge-filtration on $\tmop{dR}_{R / k}$ is complete, and it converges
  universally of index $0$.
\end{lemma}

\begin{proof}
  First, by \Cref{lem:abs-cot-cplx-retr-rel-frob}, the $R$-module $L_{R / k}$
  is a retract of $L_{R / \varphi_k^{\ast} R}$, which is zero as the relative
  Frobenius map $\varphi_{R / k} \of \varphi_k^{\ast} R \rightarrow R$ is an
  equivalence\footnote{This argument is actually implicit in
  {\cite[Lem~6.5.13(i)]{Gabber2003}} as well.}. It follows that the Hodge
  filtration on $\widehat{\tmop{dR}}_{R / k}$ is concentrated in weight $0$,
  and so is the conjugate filtration on $\tmop{dR}_{R / k}$, and thus the
  canonical map $\tmop{dR}_{R / k} \rightarrow \widehat{\tmop{dR}}_{R / k}$ is
  identified with the relative Frobenius map $\varphi_{R / k}$, and the
  Hodge-completeness follows. The universal convergence follows from
  definitions as well.
\end{proof}

\begin{example}
  Let $R$ be a valuation $k$-algebra where $k \assign \mathbb{F}_p$. Gabber
  showed that the $k$-algebra $R$ is Cartier smooth,
  {\cite[Cor~A.4]{Kerz2021}}. In particular, for every $r \in \mathbb{N}$, the
  $R$-module $\Omega_{R / k}^r$ is flat thus torsion-free, thus it is also
  torsion-free viewed as a $\varphi_k^{\ast} R$-module. It follows that the
  sub-$\varphi_k^{\ast} R$-module $B \Omega_{R / k}^r \subseteq \Omega_{R /
  k}^r$ is torsion-free as well, thus flat, as $\varphi_k^{\ast} R \cong R$ is
  a valuation ring. It then follows from
  \Cref{lem:univ-Cart-iso-Hdg-conv-univ} that the Hodge filtration on
  $\widehat{\tmop{dR}}_{R / k}$ converges universally of index $0$.
\end{example}

We now explain how to deduce an animated version of transitivity of
$p$-Cartier smoothness established in {\cite[Lem~2.10]{Bouis2023}} from
\Cref{lem:lin-conv-Hdg-fil-quasi-sm}, and the proof of
\Cref{thm:dR-Hdg-cpl-rel-Nyg-cpl} shares the same strategy. The key is to
bound the $\tmop{Tor}$-amplitude of the map $\tmop{dR}_{R / k} \rightarrow
\tmop{dR}_{A / k}$ for quasi-lci maps $R \rightarrow A$:

\begin{lemma}
  \label{lem:dR-rel-conj-fil}Let $k$ be an animated $\mathbb{F}_p$-algebra,
  and $R \rightarrow A$ be a quasi-lci map of animated $k$-algebras. Then the
  de Rham cohomology $\tmop{dR}_{A / k}$, being a $(\varphi_k^{\ast} A
  \otimes_{\varphi_k^{\ast} R}^{\mathbb{L}} \tmop{dR}_{R / k})$-module, admits
  an exhaustive increasing filtration, with associated graded pieces of the
  form
  \[ \varphi_k^{\ast} N \otimes_{\varphi_k^{\ast} R}^{\mathbb{L}} \tmop{dR}_{R
     / k} \]
  for some $A$-module $N$ of $\tmop{Tor}$-amplitude$\nosymbol \leq 0$.
\end{lemma}

\begin{proof}
  We simply take the relative conjugate filtration on $\tmop{dR}_{A / k}$
  constructed in {\cite[Prop~3.22]{Bhatt2012a}}. The associated graded pieces
  are given by
  \[ \varphi_k^{\ast}  \left( \bigwedgestar_A^s L_{A / R} \right) [- s]
     \otimes_{\varphi_k^{\ast} R}^{\mathbb{L}} \tmop{dR}_{R / k} . \]
  Now as $R \rightarrow A$ is quasi-lci, for every $s \in \mathbb{N}$, the
  $A$-module $\left( \bigwedgestar_A^s L_{A / R} \right) [- s]$ has
  $\tmop{Tor}$-amplitude$\nosymbol \leq 0$ (cf.~{\cite[Lem~4.34]{Bhatt2018}}).
  The result then follows.
\end{proof}

This bound of $\tmop{Tor}$-amplitude leads to completeness of the Katz--Oda
filtration, under either a boundedness of $\tmop{Tor}$-amplitude, or a
universal convergence of Hodge filtration.

\begin{corollary}
  \label{cor:quasi-lci-bnd-Tor-ampl-KO-cpl}Let $k$ be an animated ring, and $R
  \rightarrow A$ a $p$-quasi-lci map of animated $k$-algebras, and $\tmop{Fil}
  M$ a filtered $\tmop{Fil}_H \tmop{dR}_{R / k}$-module. Suppose that
  {\tmstrong{either}}
  \begin{itemize}
    \item the $R$-module $A$ has $p$-complete $\tmop{Tor}$-amplitude$\nosymbol
    \leq m$ for some $m \in \mathbb{Z}$, that the $R$-module $L_{R / k}$ is
    $p$-completely flat, and that $R$ is $p$-completely $n$-truncated for some
    $n \in \mathbb{N}$, with $\tmop{Fil} M = \tmop{Fil}_H 
    \widehat{\tmop{dR}}_{R / k}$; {\tmstrong{or}}
    
    \item the $\tmop{Fil}_H \tmop{dR}_{R \otimes_{\mathbb{Z}}^{\mathbb{L}}
    \mathbb{F}_p / k \otimes_{\mathbb{Z}}^{\mathbb{L}} \mathbb{F}_p}$-module
    $\tmop{Fil} M \otimes_{\mathbb{Z}}^{\mathbb{L}} \mathbb{F}_p$ converges
    universally, and that the spectrum $\varphi_{k
    \otimes_{\mathbb{Z}}^{\mathbb{L}} \mathbb{F}_p}^{\ast} (A
    \otimes_{\mathbb{Z}}^{\mathbb{L}} \mathbb{F}_p)$ is $m$-truncated for some
    $m \in \mathbb{N}$.
  \end{itemize}
  Then the Katz--Oda filtration
  \[ \tmop{Fil}_{\tmop{KO} / R} (\tmop{dR}_{A / k} \otimes_{\tmop{dR}_{R /
     k}}^{\mathbb{L}} M) \assign \tmop{dR}_{A / k} \otimes_{\tmop{dR}_{R /
     k}}^{\mathbb{L}} \tmop{Fil} M \]
  is $p$-completely complete.
\end{corollary}

\begin{proof}
  First, we may replace $k$ by $k \otimes_{\mathbb{Z}}^{\mathbb{L}}
  \mathbb{F}_p$, so without loss of generality, we may equip $k$ an animated
  $\mathbb{F}_p$-algebra structure.
  
  We fix $r \in \mathbb{N}$. Then by \Cref{lem:dR-rel-conj-fil}, the $r$-th
  Katz--Oda filtered piece
  \[ \tmop{dR}_{A / k} \otimes_{\tmop{dR}_{R / k}}^{\mathbb{L}} \tmop{Fil}^r M
  \]
  admits an exhaustive increasing filtration, with associated graded pieces of
  the form
  \begin{equation}
    \varphi_k^{\ast} N \otimes_{\varphi_k^{\ast} R}^{\mathbb{L}} \tmop{Fil}^r
    M \label{eq:gr-KO-conj-fil}
  \end{equation}
  for some $A$-module $N$ of $\tmop{Tor}$-amplitude$\nosymbol \leq 0$. We
  consider the two cases separately:
  \begin{itemize}
    \item Suppose that the $R$-module $A$ has $\tmop{Tor}$-amplitude$\nosymbol
    \leq m$, and that $R$ is $n$-truncated, with $\tmop{Fil} M = \tmop{Fil}_H 
    \widehat{\tmop{dR}}_{R / k}$. Then every $\varphi_k^{\ast} R$-module
    $\varphi_k^{\ast} N$ in \eqref{eq:gr-KO-conj-fil} has
    $\tmop{Tor}$-amplitude$\nosymbol \leq m$. It follows from
    \Cref{cor:lin-conv-bnd-Tor-ampl-Hdg-fil-quasi-sm} and the flatness of
    $L_{R / k}$ that the spectrum $\varphi_k^{\ast} N
    \otimes_{\varphi_k^{\ast} R}^{\mathbb{L}} \tmop{Fil}^r M$ is $(m + n -
    r)$-truncated, thus so is the $r$-th filtered piece $\tmop{dR}_{A / k}
    \otimes_{\tmop{dR}_{R / k}}^{\mathbb{L}} \tmop{Fil}^r M$. The result then
    follows.
    
    \item Suppose that $\tmop{Fil}_H \tmop{dR}_{R / k}$-module $\tmop{Fil} M$
    converges universally of index $\alpha \in \mathbb{Z}$, and that the
    spectrum $\varphi_k^{\ast} A$ is $m$-truncated. Then every spectrum
    $\varphi_k^{\ast} N$ in \eqref{eq:gr-KO-conj-fil} is $m$-truncated. It
    follows from \Cref{cor:Hdg-conv-univ-trunc-Tor-ampl} that the spectrum
    $\varphi_k^{\ast} N \otimes_{\varphi_k^{\ast} R}^{\mathbb{L}} \tmop{Fil}^r
    M$ is $(\alpha + m - r)$-truncated, and so is the $r$-th filtered piece
    $\tmop{dR}_{A / k} \otimes_{\tmop{dR}_{R / k}}^{\mathbb{L}} \tmop{Fil}^r
    M$. The result then follows.
  \end{itemize}
\end{proof}

Now we deduce the following animated version of transitivity of $p$-Cartier
smoothness, under either a boundedness of $\tmop{Tor}$-amplitude, or a
universal convergence of Hodge filtration.

\begin{theorem}
  \label{thm:transitivity-Cart-sm}Let $k$ be an animated ring, and $R
  \rightarrow A$ a map of animated $k$-algebras with the $A$-module $L_{A /
  R}$ being $p$-completely flat, such that {\tmstrong{either}}
  \begin{itemize}
    \item the $R$-module $A$ has bounded $p$-complete $\tmop{Tor}$-amplitude,
    the $R$-module $L_{R / k}$ is $p$-completely flat, and $R$ is
    $p$-completely $t$-bounded; {\tmstrong{or}}
    
    \item the Hodge-filtration on $\widehat{\tmop{dR}}_{R
    \otimes_{\mathbb{Z}}^{\mathbb{L}} \mathbb{F}_p / k
    \otimes_{\mathbb{Z}}^{\mathbb{L}} \mathbb{F}_p}$ converges universally,
    the spectrum $A$ is $p$-completely $t$-bounded, the $R$-module $L_{R / k}$
    has $p$-complete $\tmop{Tor}$-amplitude$\nosymbol \leq 1$, and the
    spectrum $\varphi_{k \otimes_{\mathbb{Z}}^{\mathbb{L}}
    \mathbb{F}_p}^{\ast} (A \otimes_{\mathbb{Z}}^{\mathbb{L}} \mathbb{F}_p)$
    is $t$-bounded\footnote{This is satisfied if the Frobenius map $\varphi_{k
    \otimes_{\mathbb{Z}}^{\mathbb{L}} \mathbb{F}_p}$ has bounded $p$-complete
    $\tmop{Tor}$-amplitude, as we have already assumed that $A$ is
    $p$-completely $t$-bounded.}.
  \end{itemize}
  If the $p$-completed derived de Rham cohomology $\tmop{dR}_{A / R}$ is
  $p$-completely Hodge-complete, then the filtered spectrum
  \[ \tmop{Fil}_H \tmop{dR}_{A / k} \otimes_{\tmop{Fil}_H \tmop{dR}_{R /
     k}}^{\mathbb{L}} \tmop{Fil}_H  \widehat{\tmop{dR}}_{R / k} \]
  is $p$-completely complete. In particular, if the $p$-completed derived de
  Rham cohomology $\tmop{dR}_{R / k}$ is additionally $p$-completely Hodge
  complete, then so is $\tmop{dR}_{A / k}$.
\end{theorem}

\begin{proof}
  We could replace $k$ by $k \otimes_{\mathbb{Z}}^{\mathbb{L}} \mathbb{F}_p$,
  and without loss of generality, we may endow $k$ an animated
  $\mathbb{F}_p$-algebra structure. By
  \Cref{cor:quasi-lci-bnd-Tor-ampl-KO-cpl}, it suffices to show that, for
  every $i \in \mathbb{N}$, the Hodge-filtration on
  \begin{eqnarray*}
    \tmop{gr}_{\tmop{KO} / R}^i (\tmop{dR}_{A / k} \otimes_{\tmop{dR}_{R /
    k}}^{\mathbb{L}} \widehat{\tmop{dR}}_{R / k}) & = & \tmop{dR}_{A / k}
    \otimes_{\tmop{dR}_{R / k}}^{\mathbb{L}} \tmop{gr}_H^i \tmop{dR}_{R / k}\\
    & = & \tmop{dR}_{A / R} \otimes_R^{\mathbb{L}} \left( \bigwedgestar_R^i
    L_{R / k} \right) [- i]\\
    & = & \widehat{\tmop{dR}}_{A / R} \otimes_R^{\mathbb{L}} \left(
    \bigwedgestar_R^i L_{R / k} \right) [- i]
  \end{eqnarray*}
  is complete, where the last equivalence is Hodge-completeness of
  $\tmop{dR}_{A / R}$.
  \begin{itemize}
    \item In the first case, as $R$ is $t$-bounded above, and $A$ has bounded
    $\tmop{Tor}$-amplitude over $R$, the spectrum $A$ is $t$-bounded above as
    well. The result follows from \Cref{lem:lin-conv-Hdg-fil-quasi-sm} and the
    flatness of the $A$-module $L_{A / R}$.
    
    \item In the second case, the $R$-module $\left( \bigwedgestar_R^i L_{R /
    k} \right) [- i]$ has $\tmop{Tor}$-amplitude$\nosymbol \leq 0$
    (cf.~{\cite[Lem~4.34]{Bhatt2018}}), and the result follows from
    \Cref{lem:lin-conv-Hdg-fil-quasi-sm} and the flatness of the $A$-module
    $L_{A / R}$ as well.
  \end{itemize}
\end{proof}

We can also establish a version of transitivity of universal convergence of
Hodge filtrations via the Katz--Oda filtration:

\begin{theorem}
  \label{thm:transitivity-univ-conv-Hdg-filt}Let $k$ be an animated
  $\mathbb{F}_p$-algebra, and $R \rightarrow A$ be a quasi-lci map of animated
  $k$-algebras, such that
  \begin{itemize}
    \item the Hodge-filtration on $\widehat{\tmop{dR}}_{R / k}$ converges
    universally of index $\alpha \in \mathbb{N}$, with $R$-module $L_{R / k}$
    being flat; and
    
    \item the Hodge-filtration on $\widehat{\tmop{dR}}_{A / R}$ converges
    universally of index $\beta \in \mathbb{N}$, and the Hodge-filtration on
    $\tmop{dR}_{A / R}$ is complete.
  \end{itemize}
  Then the Hodge filtration on
  \[ \tmop{dR}_{A / k} \otimes_{\tmop{dR}_{R / k}}^{\mathbb{L}}
     \widehat{\tmop{dR}}_{R / k}, \]
  viewed as a filtered $\tmop{dR}_{A / k}$-module, converges universally of
  index $\alpha + \beta$.
\end{theorem}

\begin{proof}
  Let $r \in \mathbb{N}$. By definition, we have to control the
  $\tmop{Tor}$-amplitude of the $\varphi_k^{\ast} A$-module
  \[ E_r \assign \tmop{Fil}_H^r (\tmop{dR}_{A / k} \otimes_{\tmop{dR}_{R /
     k}}^{\mathbb{L}} \widehat{\tmop{dR}}_{R / k}) . \]
  We filter this $\varphi_k^{\ast} A$-module by the Katz--Oda filtration, and
  the key observation is that, the $s$-th Katz--Oda filtration piece of
  $\tmop{Fil}_H^r (\tmop{dR}_{A / k} \otimes_{\tmop{dR}_{R / k}}^{\mathbb{L}}
  \widehat{\tmop{dR}}_{R / k})$ is the same as
  \[ \tmop{dR}_{A / k} \otimes_{\tmop{dR}_{R / k}}^{\mathbb{L}} \tmop{Fil}_H^s
     \widehat{\tmop{dR}}_{R / k} \]
  for $s \geq r$ (which is in general false for $s < r$)\footnote{More
  generally, let $\tmop{Fil} R$ be an $(\mathbb{N}, \geq)$-filtered
  $\mathbb{E}_1$-ring, and $\tmop{Fil} M$ an $(\mathbb{N}, \geq)$-filtered
  right $\tmop{Fil} R$-module, and $\tmop{Fil} N$ a $\mathbb{Z}$-filtered left
  $\tmop{Fil} R$-module concentrated in weights$\nosymbol \geq s$. Then the
  filtered map
  \[ \tmop{Fil} M \otimes_{\tmop{Fil} R}^{\mathbb{L}} \tmop{Fil} N \rightarrow
     \tmop{ins}^0 \tmop{Fil}^0 M \otimes_{\tmop{ins}^0 \tmop{Fil}^0
     R}^{\mathbb{L}} \tmop{Fil} N \]
  is an equivalence. Indeed, the description of the relative tensor product in
  terms of bar construction allows us to reduce to the special case that
  $\tmop{Fil} R$ is the unit. Then it follows directly from cofinality.}. Let
  $T_r \assign \tmop{dR}_{A / k} \otimes_{\tmop{dR}_{R / k}}^{\mathbb{L}}
  \tmop{Fil}_H^r  \widehat{\tmop{dR}}_{R / k}$. By universal convergence of
  the Hodge filtration on $\widehat{\tmop{dR}}_{R / k}$, it follows from
  \Cref{lem:dR-rel-conj-fil} that the $\varphi_k^{\ast} A$-module $T_r$ has
  $\tmop{Tor}$-amplitude$\nosymbol \leq \alpha - r \leq \alpha + \beta - r$.
  
  Now it remains to see that, the $s$-th associated graded pieces of Katz--Oda
  filtration on
  \[ \tmop{Fil}_H^r (\tmop{dR}_{A / k} \otimes_{\tmop{dR}_{R /
     k}}^{\mathbb{L}} \widehat{\tmop{dR}}_{R / k}) \]
  has $\tmop{Tor}$-amplitude$\nosymbol \leq - r$ over $\varphi_k^{\ast} A$ for
  $s < r$. This is precisely given by
  \begin{equation}
    \tmop{Fil}_H^{r - s}  \widehat{\tmop{dR}}_{A / R} \otimes_R^{\mathbb{L}}
    \left( \bigwedgestar_R^s L_{R / k} \right) [- s] \label{eq:gr-KO-dR}
  \end{equation}
  as in the proof of \Cref{thm:transitivity-Cart-sm}, which has
  $\tmop{Tor}$-amplitude$\nosymbol \leq \beta + (s - r) - s = \beta - r$ over
  \[ \varphi_R^{\ast} A = \varphi_k^{\ast} (A) \otimes_{\varphi_k^{\ast} R,
     \varphi_{R / k}}^{\mathbb{L}} R. \]
  By universal convergence of the Hodge-filtration on $\widehat{\tmop{dR}}_{R
  / k}$ and \Cref{lem:conv-univ-frob-flat}, the relative Frobenius map
  $\varphi_{R / k} \of \varphi_k^{\ast} R \rightarrow R$ has
  $\tmop{Tor}$-amplitude$\nosymbol \leq \alpha$. It follows that
  \eqref{eq:gr-KO-dR} has $\tmop{Tor}$-amplitude$\nosymbol \leq \alpha + \beta
  - r$ over $\varphi_k^{\ast} A$. To conclude, the $\varphi_{\ast} A$-module
  $E_r$ has $\tmop{Tor}$-amplitude$\nosymbol \leq \alpha + \beta - r$.
\end{proof}

\begin{remark}
  \label{rem:quant-Katz-Oda-conv}We may prove \Cref{thm:transitivity-Cart-sm}
  in the same line as in the proof of
  \Cref{thm:transitivity-univ-conv-Hdg-filt}, which does not invoke
  \Cref{cor:quasi-lci-bnd-Tor-ampl-KO-cpl} explicitly, but a more
  ``quantitative'' version of that, i.e. the $t$-boundedness of ``Katz--Oda
  tails'' $T_r$.
\end{remark}

The Nygaard filtration on $\hatPrism_R$ might be more complicated. However, if
we take its reduction $\tmop{mod} p$, we can take the {\tmdfn{twisted Nygaard
filtration}} instead.

\begin{remark}[{\cite[Cons~3.3, Rem~3.4 \& Rem~4.5]{Bhatt2023a}}]
  \label{rem:tw-Nyg-fil}Let $R$ be an animated ring. Then there exists a
  functorial $\mathbb{N}$-graded filtration $\tmop{Fil}_{\tilde{N}} \left(
  \bigoplus_{n \in \mathbb{Z}} \barPrism_R \{ n \}
  \otimes_{\mathbb{Z}}^{\mathbb{L}} \mathbb{F}_p \right)$, called the
  {\tmdfn{twisted Nygaard filtration}}, on the graded
  $\mathbb{E}_{\infty}$-ring $\bigoplus_{n \in \mathbb{Z}} \barPrism_R \{ n \}
  \otimes_{\mathbb{Z}}^{\mathbb{L}} \mathbb{F}_p$, with associated graded
  pieces given by
  \[ \tmop{cofib} \left( \mu_1 \of \bigoplus_{i \in \mathbb{N}} \bigoplus_{n
     \in \mathbb{Z}} \tmop{gr}_N^{i - p}  \Prism_R \{ n \}
     \otimes_{\mathbb{Z}}^{\mathbb{L}} \mathbb{F}_p \longrightarrow
     \bigoplus_{i \in \mathbb{N}} \bigoplus_{n \in \mathbb{Z}} \tmop{gr}_N^i 
     \Prism_R \{ n \} \otimes_{\mathbb{Z}}^{\mathbb{L}} \mathbb{F}_p \right) .
  \]
  Then
  \begin{enumerate}
    \item \label{pt:tw-Nyg-perf}the $\tmop{Fil}_N \left( \bigoplus_{n \in
    \mathbb{Z}} \barPrism_R \{ n \} \right)$-module $\tmop{Fil}_{\tilde{N}}
    \left( \bigoplus_{n \in \mathbb{Z}} \barPrism_R \{ n \}
    \otimes_{\mathbb{Z}}^{\mathbb{L}} \mathbb{F}_p \right)$ is perfect, as it
    is the cofiber of a twisted endomorphism given by a class in
    {\cite[Ex~3.2]{Bhatt2023a}};
    
    \item the $\tmop{mod} p$-reduction $\barPrism_R \{ \ast \}
    \otimes_{\mathbb{Z}}^{\mathbb{L}} \mathbb{F}_p$ of Nygaard-completed
    absolute Hodge--Tate cohomology coincides with the completion of the
    twisted Nygaard completion of $\barPrism_R \{ \ast \}
    \otimes_{\mathbb{Z}}^{\mathbb{L}} \mathbb{F}_p$; and
    
    \item if the animated ring $R$ is weakly $F$-smooth and $p$-completely
    $t$-bounded above, then then the convergence $\tmop{Fil}_{\tilde{N}}
    \left( \bigoplus_{n \in \mathbb{N}} \hatbarPrism_R \{ n \}
    \otimes_{\mathbb{Z}}^{\mathbb{L}} \mathbb{F}_p \right)$ has a linearly
    increasing coconnectivity estimate (the neutral $t$-structure on graded
    spectra implies uniformity in $n \in \mathbb{Z}$).
  \end{enumerate}
\end{remark}

\begin{warning}
  \label{warn:gr-Nyg-R-mod-struct}Let $R$ be an animated ring. The $R
  \otimes_{\mathbb{Z}}^{\mathbb{L}} \mathbb{F}_p$-module structure on
  $\tmop{gr}_N^{\ast}  \barPrism_R$ inherited from that on $\tmop{gr}_N^{\ast}
  \Prism_R$ is incompatible with the one inherited from $\barPrism_R$. They
  differ up to a Frobenius, as explained in {\cite[Cons~3.1]{Bhatt2023a}}. The
  same for the twisted Nygaard filtration. To fix the terminology, we mean the
  $R \otimes_{\mathbb{Z}}^{\mathbb{L}} \mathbb{F}_p$-module structure on
  $\tmop{gr}_N^{\ast}  \barPrism_R$ inherited from that on $\tmop{gr}_N^{\ast}
  \Prism_R$ by default, unless otherwise explicitly said.
\end{warning}

\begin{remark}
  \label{rem:tw-Nyg-in-lieu-of-Nyg}By \Cref{rem:tw-Nyg-fil}, we may replace
  $\barPrism$ by $\barPrism \otimes_{\mathbb{Z}}^{\mathbb{L}} \mathbb{F}_p$,
  and $\tmop{Fil}_N$ by $\tmop{Fil}_{\tilde{N}}$ throughout in, say,
  \Cref{rem:Nyg-rel-cpl}.
\end{remark}

\begin{example}[cf.~{\cite[Ex~3.2]{Bhatt2023a}}]
  \label{ex:perfd-twi-Nyg-fil}Let $R$ be a $p$-torsion-free perfectoid ring,
  and let $(A, I)$ denote its associated perfect (transversal) prism. Then the
  prismatic cohomology $\Prism_R$ is $A$, with Nygaard filtration given by the
  $\varphi^{- 1} (I)$-adic filtration. As the Cartier divisor $I \subseteq A$
  is distinguished, we know that $\varphi^{- 1} (I)^p$ becomes $(p)$ in
  $\barPrism_R = R$, and the filtered ring $\tmop{Fil}_{\tilde{N}} \left(
  \bigoplus_{n \in \mathbb{N}} \hatbarPrism_R \{ n \}
  \otimes_{\mathbb{Z}}^{\mathbb{L}} \mathbb{F}_p \right)$ is concretely given
  by
  \[ R \{ \ast \} / p \longleftarrow \varphi^{- 1} (I) R \{ \ast \} / p
     \longleftarrow \cdots \longleftarrow \varphi^{- 1} (I)^{p - 1} R \{ \ast
     \} / p \longleftarrow 0 \longleftarrow \cdots \]
  which is concentrated in weights $[0, p - 1]$.
\end{example}

The following lemma controls ``$\tmop{Tor}$-amplitude'' of the map
$\barPrism_k \rightarrow \barPrism_R$ for quasi-lci maps $k \rightarrow R$ of
bounded $\tmop{Tor}$-amplitude, an analogue of \Cref{lem:dR-rel-conj-fil}.

\begin{lemma}
  \label{lem:Tor-ampl-prism-coh-quasi-syn}Let $k \rightarrow R$ be a
  $p$-quasi-lci map of animated rings. Then the graded $R
  \otimes_k^{\mathbb{L}} \left( \bigoplus_{n \in \mathbb{Z}} \barPrism_k \{ n
  \} \right)$-module spectrum $\bigoplus_{n \in \mathbb{Z}} \barPrism_R \{ n
  \}$ admits an exhaustive increasing filtration, with associated graded
  pieces being retracts of graded $R \otimes_k^{\mathbb{L}} \left(
  \bigoplus_{n \in \mathbb{Z}} \barPrism_k \{ n \} \right)$-modules of the
  form
  \[ N \otimes_k^{\mathbb{L}} \left( \bigoplus_{n \in \mathbb{Z}} \barPrism_k
     \{ n \} \right) \]
  with $N$ being a $R$-module spectrum with $p$-complete
  $\tmop{Tor}$-amplitude$\nosymbol \leq 0$.
\end{lemma}

\begin{proof}
  We simply take the Katz--Oda filtration in \Cref{ex:KO-Hdg-Tate}, and the
  associated graded pieces are direct summands (thus retracts) of the graded
  graded $R \otimes_k^{\mathbb{L}} \left( \bigoplus_{n \in \mathbb{Z}}
  \barPrism_k \{ n \} \right)$-module
  \[ \left( \bigoplus_{s \in \mathbb{N}} \left( \bigwedgestar_R^s L_{R / k}
     \right) [- s] \right) \otimes_k^{\mathbb{L}} \left( \bigoplus_{n \in
     \mathbb{Z}} \barPrism_k \{ n \} \right) . \]
  Now as $k \rightarrow R$ is $p$-quasi-lci, for every $s \in \mathbb{N}$, the
  $R$-module $\left( \bigwedgestar_R^s L_{R / k} \right) [- s]$ has
  $p$-complete $\tmop{Tor}$-amplitude$\nosymbol \leq 0$
  (cf.~{\cite[Lem~4.34]{Bhatt2018}}), thus so does the $R$-module
  \[ N \assign \bigoplus_{s \in \mathbb{N}} \left( \bigwedgestar_R^s L_{R / k}
     \right) [- s] . \]
\end{proof}

With this $\tmop{Tor}$-amplitude of absolute prismatic cohomology, we deduce
the completeness of ``Katz--Oda filtration'' out of twisted Nygaard
filtration, just as the first case of
\Cref{cor:quasi-lci-bnd-Tor-ampl-KO-cpl}:

\begin{corollary}
  \label{cor:KO-tw-Nyg-cpl}Let $k \rightarrow R$ be a $p$-quasi-lci map of
  animated rings such that $k$ is $p$-completely $t$-bounded and $R$ has
  bounded $p$-complete $\tmop{Tor}$-amplitude over $k$. Suppose that the
  animated ring $k$ is weakly $F$-smooth. Then the filtered graded spectrum
  \[ \left( \bigoplus_{n \in \mathbb{Z}} \barPrism_R \{ n \}
     \otimes_{\mathbb{Z}}^{\mathbb{L}} \mathbb{F}_p \right)
     \otimes_{\bigoplus_{n \in \mathbb{Z}} \barPrism_k \{ n \}
     \otimes_{\mathbb{Z}}^{\mathbb{L}} \mathbb{F}_p}^{\mathbb{L}} \left(
     \bigoplus_{n \in \mathbb{Z}} \tmop{Fil}_{\tilde{N}} \left( \hatbarPrism_k
     \{ n \} \otimes_{\mathbb{Z}}^{\mathbb{L}} \mathbb{F}_p \right) \right) \]
  is complete, and its convergence has a linearly increasing $p$-complete
  coconnectivity estimate.
\end{corollary}

\begin{proof}
  The result follows from
  \Cref{rem:tw-Nyg-fil,lem:Tor-ampl-prism-coh-quasi-syn} (cf. the proof of
  \Cref{cor:quasi-lci-bnd-Tor-ampl-KO-cpl}).
\end{proof}

Now we follow the same line as in the proof of
\Cref{thm:transitivity-Cart-sm}.

\begin{proof*}{Proof of \Cref{thm:dR-Hdg-cpl-rel-Nyg-cpl}}
  By \Cref{rem:Nyg-rel-cpl,rem:tw-Nyg-in-lieu-of-Nyg,cor:KO-tw-Nyg-cpl}, it
  suffices to check that, for every $j \in \mathbb{N}$, the filtered graded
  spectrum
  \[ F_j \assign \left( \bigoplus_{n \in \mathbb{Z}} \tmop{Fil}_{\tilde{N}}
     \left( \barPrism_R \{ n \} \otimes_{\mathbb{Z}}^{\mathbb{L}} \mathbb{F}_p
     \right) \right) \otimes_{\bigoplus_{n \in \mathbb{Z}}
     \tmop{Fil}_{\tilde{N}} \left( \barPrism_k \{ n \}
     \otimes_{\mathbb{Z}}^{\mathbb{L}} \mathbb{F}_p \right)}^{\mathbb{L}}
     \left( \bigoplus_{n \in \mathbb{Z}} \tmop{ins}^j \tmop{gr}_{\tilde{N}}^j 
     \barPrism_k \{ n \} \otimes_{\mathbb{Z}}^{\mathbb{L}} \mathbb{F}_p
     \right) \]
  is complete. We start with $j = 0$. By point~\ref{pt:tw-Nyg-perf} of
  \Cref{rem:tw-Nyg-fil}, it suffices to check the completeness of the filtered
  graded spectrum
  \[ \left( \bigoplus_{n \in \mathbb{Z}} \tmop{Fil}_N  \barPrism_R \{ n \}
     \right) \otimes_{\bigoplus_{n \in \mathbb{Z}} \tmop{Fil}_N  \barPrism_k
     \{ n \}}^{\mathbb{L}} \left( \bigoplus_{n \in \mathbb{Z}} \tmop{ins}^0
     \tmop{gr}_N^0  \barPrism_k \{ n \} \right) \]
  after $p$-completion. By de Rham comparison, this filtered spectrum can be
  identified with the filtered graded spectrum $\bigoplus_{n \in \mathbb{Z}}
  \tmop{Fil}_H \tmop{dR}_{R / k} \otimes_{\mathbb{Z}}^{\mathbb{L}}
  \mathbb{F}_p$, which is complete.
  
  Now we fix $j \in \mathbb{N}$. Note that we have
  \[ F_j = F_0 \otimes_{\bigoplus_{n \in \mathbb{Z}} \tmop{ins}^0
     \tmop{gr}_{\tilde{N}}^0 \left( \barPrism_k \{ n \}
     \otimes_{\mathbb{Z}}^{\mathbb{L}} \mathbb{F}_p \right)}^{\mathbb{L}}
     \bigoplus_{n \in \mathbb{Z}} \tmop{ins}^j \tmop{gr}_{\tilde{N}}^j \left(
     \barPrism_k \{ n \} \otimes_{\mathbb{Z}}^{\mathbb{L}} \mathbb{F}_p
     \right) . \]
  By weak $F$-smoothness of $k$ and \Cref{rem:tw-Nyg-fil} that the filtered
  $\left( \bigoplus_{n \in \mathbb{Z}} \tmop{ins}^0 \tmop{gr}_{\tilde{N}}^0
  \left( \barPrism_k \{ n \} \otimes_{\mathbb{Z}}^{\mathbb{L}} \mathbb{F}_p
  \right) \right)$-module
  \[ \bigoplus_{n \in \mathbb{Z}} \tmop{ins}^j \tmop{gr}_{\tilde{N}}^j \left(
     \barPrism_k \{ n \} \otimes_{\mathbb{Z}}^{\mathbb{L}} \mathbb{F}_p
     \right) \]
  has bounded $\tmop{Tor}$-amplitude, and the completeness of $F_j$ then
  follows from the completeness of $F_0$ and the convergence estimate by
  \Cref{lem:lin-conv-Hdg-fil-quasi-sm}.
\end{proof*}

Finally, inspired by \Cref{def:Hdg-filt-univ-conv}, we introduce universal
convergence for Nygaard filtration:

\begin{definition}
  \label{def:Nyg-filt-univ-conv}Let $A$ be an animated ring. We say that a
  filtered graded $\tmop{Fil}_{\tilde{N}} \left( \barPrism_A \{ \ast \}
  \otimes_{\mathbb{Z}}^{\mathbb{L}} \mathbb{F}_p \right)$-module $\tmop{Fil}
  M$ {\tmdfn{converges universally of index $\alpha \in \mathbb{Z}$}} if, for
  every $r \in \mathbb{N}$, the graded $A$-module $\tmop{Fil}^r M$, inherited
  from the $A$-module structure on $\barPrism_A$
  (cf.~Warning~\ref{warn:gr-Nyg-R-mod-struct}) has $p$-complete
  $\tmop{Tor}$-amplitude\footnote{By \Cref{cor:A-mod-I-mods-I-cpl-Tor-ampl},
  it is equivalent to non-$p$-complete $\tmop{Tor}$-amplitude.}$\nosymbol \leq
  \alpha - r$, and it {\tmdfn{converges universally}} if such an $\alpha \in
  \mathbb{Z}$ exists. In particular, we say that the Nygaard
  filtration\footnote{Here we are actually considering the twisted Nygaard
  filtration is actually on $\hatbarPrism_A \{ \ast \}
  \otimes_{\mathbb{Z}}^{\mathbb{L}} \mathbb{F}_p$, and thus this is an abuse
  of language.} on $\hatbarPrism_A \{ \ast \}$ {\tmdfn{converges universally
  (resp. of index $\alpha \in \mathbb{Z}$)}} if the filtered
  $\tmop{Fil}_{\tilde{N}} \left( \barPrism_A \{ \ast \}
  \otimes_{\mathbb{Z}}^{\mathbb{L}} \mathbb{F}_p \right)$-module
  $\tmop{Fil}_{\tilde{N}} \left( \hatbarPrism_A \{ \ast \}
  \otimes_{\mathbb{Z}}^{\mathbb{L}} \mathbb{F}_p \right)$ converges
  universally (resp. of index $\alpha \in \mathbb{Z}$).
\end{definition}

\begin{remark}
  \label{rem:Fil-tilde-N-A-mod-struct}Let $A$ be an animated ring. Unlike in
  Warning~\ref{warn:gr-Nyg-R-mod-struct}, when we talk about universal
  convergence, we mean the $A$-module structure on $\tmop{Fil}_{\tilde{N}}
  \left( \barPrism_A \{ \ast \} \otimes_{\mathbb{Z}}^{\mathbb{L}} \mathbb{F}_p
  \right)$ inherited from the $A$-module structure on $\barPrism_A$, while
  when we talk about associated graded pieces $\tmop{gr}_{\tilde{N}} \left(
  \barPrism_A \{ \ast \} \otimes_{\mathbb{Z}}^{\mathbb{L}} \mathbb{F}_p
  \right)$, we refer to that inherited from $\tmop{gr}_N  \Prism_A$.
\end{remark}

Similar to \Cref{lem:conv-univ-frob-flat}, the universal convergence of
Nygaard filtration gives rise to a bound of $\tmop{Tor}$-amplitude of
Frobenius.

\begin{lemma}
  \label{lem:Nyg-conv-univ-Tor-Frob}Let $R$ be an animated ring, and
  $\tmop{Fil} M$ a filtered graded $\tmop{Fil}_{\tilde{N}} \left( \barPrism_R
  \{ \ast \} \otimes_{\mathbb{Z}}^{\mathbb{L}} \mathbb{F}_p \right)$-module
  equipped with a map $\tmop{Fil}_{\tilde{N}} \left( \barPrism_R \{ \ast \}
  \otimes_{\mathbb{Z}}^{\mathbb{L}} \mathbb{F}_p \right) \rightarrow
  \tmop{Fil} M$ which induces an equivalence on associated graded pieces.
  Suppose that the $\tmop{Fil}_{\tilde{N}} \left( \barPrism_R \{ \ast \}
  \otimes_{\mathbb{Z}}^{\mathbb{L}} \mathbb{F}_p \right)$-module $\tmop{Fil}
  M$ converges universally of index $\alpha \in \mathbb{Z}$. Then the
  Frobenius map $\psi_R \of R \rightarrow R \otimes_{\mathbb{Z}}^{\mathbb{L}}
  \mathbb{F}_p$ has $p$-complete $\tmop{Tor}$-amplitude$\nosymbol \leq
  \alpha$.
\end{lemma}

\begin{proof}
  We examine the fiber sequence
  \[ \begin{array}{ccccc}
       &  & \tmop{gr}_{\tilde{N}}^0 \left( \barPrism_R \{ \ast \}
       \otimes_{\mathbb{Z}}^{\mathbb{L}} \mathbb{F}_p \right) &  & \\
       &  & \longdownarrow \nocomma \simeq &  & \\
       \tmop{Fil}^0 M & \longrightarrow & \tmop{gr}^0 M & \longrightarrow &
       \tmop{Fil}^1 M [1]
     \end{array} \]
  of $R$-modules, where unlike in \Cref{rem:Fil-tilde-N-A-mod-struct}, we view
  the top $\tmop{gr}_{\tilde{N}}^0 \left( \barPrism_R \{ \ast \}
  \otimes_{\mathbb{Z}}^{\mathbb{L}} \mathbb{F}_p \right)$ as an $R$-module
  inherited from the map $R \rightarrow \barPrism_R$. By assumption, we know
  that the map $\psi_R \of R \rightarrow R \otimes_{\mathbb{Z}}^{\mathbb{L}}
  \mathbb{F}_p = \tmop{gr}_{\tilde{N}}^0 \left( \barPrism_R \{ \ast \}
  \otimes_{\mathbb{Z}}^{\mathbb{L}} \mathbb{F}_p \right)$ has
  $\tmop{Tor}$-amplitude$\nosymbol \leq \alpha$.
\end{proof}

Similarly to the proof of the second case of \Cref{thm:transitivity-Cart-sm},
it follows from the proof of \Cref{thm:dR-Hdg-cpl-rel-Nyg-cpl} that the
following variant of transitivity without boundedness of
$\tmop{Tor}$-amplitude when the Nygaard filtration on the base converges
universally.

\begin{theorem}
  \label{thm:transitivity-F-sm-base-conv-univ}Let $k \rightarrow R$ be a
  quasi-lci map of animated rings such that $R$ is $p$-completely $t$-bounded.
  Suppose that
  \begin{itemize}
    \item the animated ring $k$ is weakly $F$-smooth and the Nygaard
    filtration on $\hatbarPrism_k \{ \ast \}$ converges universally; and that
    
    \item the $R$-module spectrum $L_{R / k}$ is $p$-completely flat.
  \end{itemize}
  If the $p$-completed derived de Rham cohomology $\tmop{dR}_{R / k}$ is
  $p$-completely Hodge-complete, then the absolute prismatic cohomology
  $\Prism_R \{ \ast \}$ is Nygaard-complete relative to $\Prism_k \{ \ast \}$.
\end{theorem}

Adapting the proof of \Cref{thm:dR-Hdg-cpl-rel-Nyg-cpl} in the line the proof
of \Cref{thm:transitivity-univ-conv-Hdg-filt}
(cf.~\Cref{rem:quant-Katz-Oda-conv}), we get a transitivity of universal
convergence of Nygaard filtration.

\begin{theorem}
  \label{thm:transitivity-univ-conv-Nyg-filt}Let $k \rightarrow R$ be a
  quasi-lci map of animated rings. Suppose that
  \begin{itemize}
    \item the animated ring $k$ is weakly $F$-smooth and the Nygaard
    filtration on $\hatbarPrism_k \{ \ast \}$ converges universally of index
    $\alpha \in \mathbb{N}$; and that
    
    \item the Hodge filtration on $\tmop{dR}_{R / k}$ is $p$-completely
    complete, and the Hodge filtration on $\widehat{\tmop{dR}}_{R
    \otimes_{\mathbb{Z}}^{\mathbb{L}} \mathbb{F}_p / k
    \otimes_{\mathbb{Z}}^{\mathbb{L}} \mathbb{F}_p}$ converges universally of
    index $\beta \in \mathbb{N}$.
  \end{itemize}
  Then the filtered graded $\tmop{Fil}_{\tilde{N}} \left( \barPrism_R \{ \ast
  \} \otimes_{\mathbb{Z}}^{\mathbb{L}} \mathbb{F}_p \right)$-module
  \[ \tmop{Fil}_{\tilde{N}} \left( \barPrism_R \{ \ast \}
     \otimes_{\mathbb{Z}}^{\mathbb{L}} \mathbb{F}_p \right)
     \otimes_{\tmop{Fil}_{\tilde{N}} \left( \barPrism_k \{ \ast \}
     \otimes_{\mathbb{Z}}^{\mathbb{L}} \mathbb{F}_p \right)}^{\mathbb{L}}
     \tmop{Fil}_{\tilde{N}} \left( \hatbarPrism_k \{ \ast \}
     \otimes_{\mathbb{Z}}^{\mathbb{L}} \mathbb{F}_p \right) \]
  converges universally.
\end{theorem}

\begin{proof}
  Let $E_r$ denote the $r$-th filtered piece of the filtered graded
  $\tmop{Fil}_{\tilde{N}} \left( \barPrism_R \{ \ast \}
  \otimes_{\mathbb{Z}}^{\mathbb{L}} \mathbb{F}_p \right)$-module above. By
  definition, we have to control its $\tmop{Tor}$-amplitude as a graded
  $R$-module, inherited from the graded $R$-module structure on $\barPrism_R
  \{ \ast \}$. As in the proof of \Cref{thm:transitivity-univ-conv-Hdg-filt},
  we endow $E_r$ with Katz--Oda filtration as in the proof of
  \Cref{thm:dR-Hdg-cpl-rel-Nyg-cpl}, induced by the twisted Nygaard filtration
  on $\hatbarPrism_k \{ \ast \} \otimes_{\mathbb{Z}}^{\mathbb{L}}
  \mathbb{F}_p$, and we have to control the tail $\tmop{Fil}_{\tmop{KO}}^r
  E_r$ and $\tmop{gr}_{\tmop{KO}}^s E_r$ for $0 \leq s < r$.
  
  As in the proof of \Cref{thm:transitivity-univ-conv-Hdg-filt}, the tail
  $\tmop{Fil}_{\tmop{KO}}^r E_r$ is simply given by
  \[ \left( \bigoplus_{n \in \mathbb{Z}} \barPrism_R \{ n \}
     \otimes_{\mathbb{Z}}^{\mathbb{L}} \mathbb{F}_p \right)
     \otimes_{\bigoplus_{n \in \mathbb{Z}} \barPrism_k \{ n \}
     \otimes_{\mathbb{Z}}^{\mathbb{L}} \mathbb{F}_p}^{\mathbb{L}} \left(
     \bigoplus_{n \in \mathbb{Z}} \tmop{Fil}_{\tilde{N}}^r \left(
     \hatbarPrism_k \{ n \} \otimes_{\mathbb{Z}}^{\mathbb{L}} \mathbb{F}_p
     \right) \right) . \]
  Now the base change along $\barPrism_k \{ \ast \}
  \otimes_{\mathbb{Z}}^{\mathbb{L}} \mathbb{F}_p \rightarrow \barPrism_R \{
  \ast \} \otimes_{\mathbb{Z}}^{\mathbb{L}} \mathbb{F}_p$ above can be
  rewritten as the ($p$-completed) base change along $\barPrism_k \{ \ast \}
  \rightarrow \barPrism_R \{ \ast \}$, and it follows from
  \Cref{lem:Tor-ampl-prism-coh-quasi-syn} that the graded $R$-module
  $\tmop{Fil}_{\tmop{KO}}^r E_r$ has $p$-complete
  $\tmop{Tor}$-amplitude$\nosymbol \leq \alpha - r$.
  
  Now we control $\tmop{gr}_{\tmop{KO}}^s E_r$ for $0 \leq s < r$. As in the
  proof of \Cref{thm:dR-Hdg-cpl-rel-Nyg-cpl}, this is given by
  \[ \left( \bigoplus_{n \in \mathbb{Z}} \tmop{Fil}_H^{r - s} \tmop{dR}_{R /
     k} \otimes_{\mathbb{Z}}^{\mathbb{L}} \mathbb{F}_p \right)
     \otimes_{\bigoplus_{n \in \mathbb{Z}} \tmop{gr}_{\tilde{N}}^0 \left(
     \barPrism_k \{ n \} \otimes_{\mathbb{Z}}^{\mathbb{L}} \mathbb{F}_p
     \right)}^{\mathbb{L}} \bigoplus_{n \in \mathbb{Z}}
     \tmop{gr}_{\tilde{N}}^s \left( \barPrism_k \{ n \}
     \otimes_{\mathbb{Z}}^{\mathbb{L}} \mathbb{F}_p \right) . \]
  As a graded $\varphi_{k \otimes_{\mathbb{Z}}^{\mathbb{L}}
  \mathbb{F}_p}^{\ast}  (R \otimes_{\mathbb{Z}}^{\mathbb{L}}
  \mathbb{F}_p)$-module, it has $\tmop{Tor}$-amplitude$\nosymbol \leq (\beta -
  (r - s)) + (p - 1 - s) = p + \beta - r - 1$. Note that
  \[ \varphi_{k \otimes_{\mathbb{Z}}^{\mathbb{L}} \mathbb{F}_p}^{\ast} (R
     \otimes_{\mathbb{Z}}^{\mathbb{L}} \mathbb{F}_p) = R \otimes_{k,
     \psi_k}^{\mathbb{L}} (k \otimes_{\mathbb{Z}}^{\mathbb{L}} \mathbb{F}_p),
  \]
  and it follows from \Cref{lem:Nyg-conv-univ-Tor-Frob} that the $R$-module
  $\tmop{gr}_{\tmop{KO}}^s E_r$ has $p$-complete
  $\tmop{Tor}$-amplitude$\nosymbol \leq p + \alpha + \beta - 1 - r$.
\end{proof}

As a consequence, we deduce a version of {\cite[Prop~4.7]{Bhatt2023a}} without
boundedness of $\tmop{Tor}$-amplitude (or in particular, flatness) for maps $k
\rightarrow R$ whose modulo $p$ reduction is relatively perfect.

\begin{corollary}
  \label{cor:F-sm-rel-perf}Let $k \rightarrow R$ be a map of animated rings
  such that the induced map $k \otimes_{\mathbb{Z}}^{\mathbb{L}} \mathbb{F}_p
  \rightarrow R \otimes_{\mathbb{Z}}^{\mathbb{L}} \mathbb{F}_p$ of animated
  $\mathbb{F}_p$-algebras is relatively perfect, and that $R$ is
  $p$-completely $t$-bounded. If
  \begin{itemize}
    \item the animated ring $k$ is weakly $F$-smooth, the absolute prismatic
    cohomology $\Prism_k$ is Nygaard-complete, and the Nygaard filtration on
    $\hatbarPrism_k \{ \ast \}$ converges universally,
  \end{itemize}
  then the same is true for $R$. That is to say,
  \begin{itemize}
    \item the animated ring $R$ is weakly $F$-smooth, the absolute prismatic
    cohomology $\Prism_R$ is Nygaard-complete, and the Nygaard filtration on
    $\hatbarPrism_R \{ \ast \}$ converges universally.
  \end{itemize}
\end{corollary}

\begin{proof}
  It follows from
  \Cref{thm:weak-F-sm-quasism,thm:transitivity-F-sm-base-conv-univ,thm:transitivity-univ-conv-Nyg-filt,lem:rel-perf-Hdg-filt-cpl-conv-univ}.
\end{proof}

Finally, we compare the universal convergence of Hodge and Nygaard filtration
over a perfectoid base, an analogue of {\cite[Prop~4.12]{Bhatt2023a}} or
{\cite[Thm~2.16~\&~2.18]{Bouis2023}}.

\begin{theorem}
  \label{thm:Hdg-Nyg-conv-univ}Let $S$ be a perfectoid ring, and $R$ an
  animated $S$-algebra. If
  \begin{itemize}
    \item the Hodge filtration on the Hodge-completed derived de Rham
    cohomology $\widehat{\tmop{dR}}_{(R \otimes_{\mathbb{Z}}^{\mathbb{L}}
    \mathbb{F}_p) / (S \otimes_{\mathbb{Z}}^{\mathbb{L}} \mathbb{F}_p)}$
    converges universally,
  \end{itemize}
  then
  \begin{itemize}
    \item the Nygaard filtration on the Nygaard-completed absolute prismatic
    cohomology $\hatPrism_R \{ \ast \}$ converges universally.
  \end{itemize}
\end{theorem}

Actually, this admits a partial inverse, as we will see in
\Cref{prop:Hdg-Nyg-conv-univ}. As we fix a perfectoid base, we may suppress
all Breuil--Kisin twists $\{ \ast \}$. We start with some recollections on
absolute prismatic cohomology of rings over a perfectoid base.

\begin{remark}
  \label{rem:recol-rel-prism-coh}Let $(A, I)$ be a prism with $\overline{A}
  \assign A / I$, and $R$ an animated $\overline{A}$-algebra. Then we have the
  following statements.
  \begin{enumerate}
    \item When the prism $(A, I)$ is perfect, the Nygaard-filtered absolute
    prismatic cohomology $\tmop{Fil}_N  \Prism_R$ can be identified with the
    Nygaard-filtered relative prismatic cohomology $\tmop{Fil}_N
    \varphi_A^{\ast}  \Prism_{R / A}$, by {\cite[Thm~5.6.2]{Bhatt2022}}.
    
    \item The Hodge-filtered derived de Rham cohomology $\tmop{dR}_{R / S}$
    can be identified with
    \[ \tmop{Fil}_N \left( \varphi_A^{\ast}  \Prism_{R / A} \right)
       \otimes_{(I^{\geqslant 0} A)} \tmop{ins}^0  \overline{A} \]
    by {\cite[Prop~5.2.3 \& Cor~5.2.8]{Bhatt2022}}, where $(I^{\geqslant 0}
    A)$ is the ring $A$ with $I$-adic filtration, plus the observation that
    the inclusion $d \of I \rightarrow A$ induces a map $I
    \otimes_A^{\mathbb{L}} (I^{\geqslant 0} A) (1) \rightarrow I^{\geqslant 0}
    A$ of filtered $(I^{\geqslant 0} A)$-modules
    (cf.~{\cite[Rem~5.1.9]{Bhatt2022}}) whose cofiber is $\tmop{ins}^0 
    \overline{A}$. The same comparison for the Nygaard and Hodge-completed
    versions.
    
    \item The modulo $p$ reduction (or equivalently, modulo $\varphi
    (I)$-reduction) of de Rham comparison
    \[ \varphi_A^{\ast}  \Prism_{R / A} \otimes_A^{\mathbb{L}} \overline{A}
       \simeq \tmop{dR}_{R / \overline{A}} \]
    can be rewritten as
    \[ \varphi_A^{\ast}  \barPrism_{R / A} \otimes_{A / \varphi
       (I)}^{\mathbb{L}} (\overline{A} \otimes_{\mathbb{Z}}^{\mathbb{L}}
       \mathbb{F}_p) \simeq \tmop{dR}_{(R \otimes_{\mathbb{Z}}^{\mathbb{L}}
       \mathbb{F}_p) / (\overline{A} \otimes_{\mathbb{Z}}^{\mathbb{L}}
       \mathbb{F}_p)} \]
    which can be upgraded to an (increasingly) filtered equivalence, where
    both sides are equipped with conjugate filtration. In particular, this
    equivalence is an equivalence of derived $\varphi_{\overline{A}
    \otimes_{\mathbb{Z}}^{\mathbb{L}} \mathbb{F}_p}^{\ast} (R
    \otimes_{\mathbb{Z}}^{\mathbb{L}} \mathbb{F}_p)$-algebras, where the
    derived algebra structure on the left hand side is given by the
    equivalence
    \[ \varphi_A^{\ast} (R) \otimes_{A / \varphi (I)}^{\mathbb{L}}
       (\overline{A} \otimes_{\mathbb{Z}}^{\mathbb{L}} \mathbb{F}_p) \simeq
       \varphi_{\overline{A} \otimes_{\mathbb{Z}}^{\mathbb{L}}
       \mathbb{F}_p}^{\ast} (R \otimes_{\mathbb{Z}}^{\mathbb{L}} \mathbb{F}_p)
    \]
    of animated rings. When $(A, I)$ is a transversal prism and $R$ is a
    polynomial $\overline{A}$-algebra, both conjugate filtrations coincide
    with Whitehead towers, and the general case follows from left Kan
    extension\footnote{More precisely, we can pick a transversal prism $(A_0,
    I_0)$ equipped with a map $(A_0, I_0) \rightarrow (A, I)$ of prisms as in
    {\cite[Prop~.2.41]{Bhatt2022}}, and then we can left Kan extend from the
    case $(A, A / I \rightarrow R)$, where $A$ is a free
    $\delta$-$A_0$-algebra and $R$ is a polynomial $A / I$-algebra, to general
    animated $\delta$-$A_0$-algebras $A$ with an animated $A / I$-algebra $R$.
    Strictly speaking, this construction depends on the choice of $(A_0, I_0)
    \rightarrow (A, I)$, but in our case, we fix a base prism $(A, I)$, and we
    will not consider the functoriality with respect to the base prism $(A,
    I)$.}.
  \end{enumerate}
\end{remark}

The key to \Cref{thm:Hdg-Nyg-conv-univ} is the following observation.

\begin{lemma}
  \label{lem:tw-Nyg-fil-mod-d-p}Let $(A, I)$ be a perfect prism, and $R$ an
  animated $A / I$-algebra. We start with filtered object
  \[ \tmop{Fil}_N \left( \varphi_A^{\ast}  \Prism_{R / A} \right)
     \otimes_{(I^{\geqslant 0} A)} ((I^{\geqslant 0} A) / d^p) \]
  where $d^p$ is the map $I^p \otimes_A^{\mathbb{L}} (I^{\geqslant 0} A) (p)
  \rightarrow I^{\geqslant 0} A$, the $p$-th composite of the map $d$ in
  \Cref{rem:recol-rel-prism-coh}. Then the derived base change this filtered
  object along $A \rightarrow A / \varphi (I)$ can be identified with
  $\varphi_A^{\ast} \tmop{Fil}_{\tilde{N}} \left( \barPrism_R
  \otimes_{\mathbb{Z}}^{\mathbb{L}} \mathbb{F}_p \right)$. Concretely, we have
  a fiber sequence
  \[ I^p \otimes_A^{\mathbb{L}} \tmop{Fil}_N^{r - p} \left( \varphi_A^{\ast} 
     \barPrism_{R / A} \right) \xrightarrow{d^p} \tmop{Fil}_N^r \left(
     \varphi_A^{\ast}  \barPrism_{R / A} \right) \longrightarrow
     \varphi_A^{\ast} \tmop{Fil}_{\tilde{N}}^r \left( \barPrism_R
     \otimes_{\mathbb{Z}}^{\mathbb{L}} \mathbb{F}_p \right) \]
  of $R$-modules. The same for Nygaard-completed versions.
\end{lemma}

\begin{proof}
  This follows essentially from the computation in
  {\cite[Ex~3.2]{Bhatt2023a}}, cf.~\Cref{ex:perfd-twi-Nyg-fil}.
\end{proof}

\Cref{thm:Hdg-Nyg-conv-univ} then follows from
\Cref{lem:tw-Nyg-fil-mod-d-p,rem:recol-rel-prism-coh,prop:Hdg-Nyg-conv-univ}.
Actually, \Cref{prop:Hdg-Nyg-conv-univ} gives rise to a partial inverse of
\Cref{thm:Hdg-Nyg-conv-univ} as well.

\begin{proposition}
  \label{prop:Hdg-Nyg-conv-univ}Let $(A, I)$ be a prism with $\overline{A}
  \assign A / I$, and $R$ an animated $\overline{A}$-algebra. Let $\alpha \in
  \mathbb{N}$ be a natural number. Then the following statements are
  equivalent:
  \begin{enumerate}
    \item For every $r \in \mathbb{N}$, the $\varphi_A^{\ast} R$-module
    $\tmop{Fil}_H^r  \widehat{\tmop{dR}}_{R \otimes_{\mathbb{Z}}^{\mathbb{L}}
    \mathbb{F}_p / \overline{A} \otimes_{\mathbb{Z}}^{\mathbb{L}}
    \mathbb{F}_p}$ has $p$-complete $\tmop{Tor}$-amplitude$\nosymbol \leq
    \alpha - r$.
    
    \item For every $r \in \mathbb{N}$, the $r$-th filtered piece $F_r$ of
    \[ \tmop{Fil}_N \left( \varphi_A^{\ast}  \hatbarPrism_{R / A} \right)
       \otimes_{(I^{\geqslant 0} A)} ((I^{\geqslant 0} A) / d^p), \]
    as a $\varphi_A^{\ast} (R)$-module, has $p$-complete
    $\tmop{Tor}$-amplitude$\nosymbol \leq p - 1 + \alpha - r$.
  \end{enumerate}
\end{proposition}

\begin{proof}
  We start with supposing the first statement. By the second point of
  \Cref{rem:recol-rel-prism-coh}, we see that the $\varphi_A^{\ast} R$-module
  $F_r$ is a finite extension of $(\tmop{Fil}_H^i  \widehat{\tmop{dR}}_{R
  \otimes_{\mathbb{Z}}^{\mathbb{L}} \mathbb{F}_p / \overline{A}
  \otimes_{\mathbb{Z}}^{\mathbb{L}} \mathbb{F}_p})_{i = r - p + 1}^r$, thus it
  has $p$-complete $\tmop{Tor}$-amplitude$\nosymbol \leq \alpha - (r - p + 1)
  = p - 1 + \alpha - r$.
  
  Now we start with the second statement. Morally, we argue by realizing
  $\tmop{Fil}_N \left( \varphi_A^{\ast}  \hatbarPrism_{R / A} \right) / d$, up
  to a shift, as a retract of $\left( \tmop{Fil}_N \left( \varphi_A^{\ast} 
  \hatbarPrism_{R / A} \right) / d^p \right) / d$\footnote{This is not as
  obvious as it seems: the Moore spectrum $\mathbb{S}/ 2$ is not a retract of
  the spectrum $(\mathbb{S}/ 2) / 2$.}. Concretely, we have commutative
  diagram
  \begin{equation}
    \begin{array}{ccccc}
      I^{p + 1} \otimes_A^{\mathbb{L}} \tmop{Fil}_N^{r - p - 1} \left(
      \varphi_A^{\ast}  \hatbarPrism_{R / A} \right) & \longrightarrow & I
      \otimes_A^{\mathbb{L}} \tmop{Fil}_N^{r - 1} \left( \varphi_A^{\ast} 
      \hatbarPrism_{R / A} \right) & \longrightarrow & I
      \otimes_A^{\mathbb{L}} F_{r - 1}\\
      \longdownarrow &  & \longdownarrow &  & \longdownarrow\\
      I^p \otimes_A^{\mathbb{L}} \tmop{Fil}_N^{r - p} \left( \varphi_A^{\ast} 
      \hatbarPrism_{R / A} \right) & \longrightarrow & \tmop{Fil}_N^r \left(
      \varphi_A^{\ast}  \hatbarPrism_{R / A} \right) & \longrightarrow & F_r\\
      \longdownarrow &  & \longdownarrow &  & \\
      I^p \otimes_A^{\mathbb{L}} \tmop{Fil}_H^{r - p}  \widehat{\tmop{dR}}_{R
      \otimes_{\mathbb{Z}}^{\mathbb{L}} \mathbb{F}_p / \overline{A}
      \otimes_{\mathbb{Z}}^{\mathbb{L}} \mathbb{F}_p} & \longrightarrow &
      \tmop{Fil}_H^r  \widehat{\tmop{dR}}_{R \otimes_{\mathbb{Z}}^{\mathbb{L}}
      \mathbb{F}_p / \overline{A} \otimes_{\mathbb{Z}}^{\mathbb{L}}
      \mathbb{F}_p} &  & 
    \end{array} \label{eq:mod-dp-d}
  \end{equation}
  of $\varphi_A^{\ast} R$-modules, where rows and columns are fiber sequences.
  This gives rise to a fiber sequence
  \[ G_r \assign \tmop{fib} (I \otimes_A^{\mathbb{L}} F_{r - 1} \rightarrow
     F_r) \longrightarrow I^p \otimes_A^{\mathbb{L}} \tmop{Fil}_H^{r - p} 
     \widehat{\tmop{dR}}_{R \otimes_{\mathbb{Z}}^{\mathbb{L}} \mathbb{F}_p /
     \overline{A} \otimes_{\mathbb{Z}}^{\mathbb{L}} \mathbb{F}_p}
     \longrightarrow \tmop{Fil}_H^r  \widehat{\tmop{dR}}_{R
     \otimes_{\mathbb{Z}}^{\mathbb{L}} \mathbb{F}_p / \overline{A}
     \otimes_{\mathbb{Z}}^{\mathbb{L}} \mathbb{F}_p} \]
  of $\varphi_A^{\ast} R$-module. Now the key observation is that the second
  map is null-homotopic, which is a consequence of
  Construction~\ref{cons:null-htpy-lift-cofibs} and the fact that there exists
  a bottom-left to top-right map
  \[ d^{p - 1} \of I^p \otimes_A^{\mathbb{L}} \tmop{Fil}_N^{r - p} \left(
     \varphi_A^{\ast}  \hatbarPrism_{R / A} \right) \longrightarrow I
     \otimes_A^{\mathbb{L}} \tmop{Fil}_N^{r - 1} \left( \varphi_A^{\ast} 
     \hatbarPrism_{R / A} \right) \]
  in the top left square of \eqref{eq:mod-dp-d}, which makes the diagram
  commute. It follows that
  \[ I^p \otimes_A^{\mathbb{L}} \tmop{Fil}_H^{r - p}  \widehat{\tmop{dR}}_{R
     \otimes_{\mathbb{Z}}^{\mathbb{L}} \mathbb{F}_p / \overline{A}
     \otimes_{\mathbb{Z}}^{\mathbb{L}} \mathbb{F}_p} \]
  is a retract of $G_r$. The fiber sequence
  \[ F_r [- 1] \longrightarrow G_r \longrightarrow I \otimes_A^{\mathbb{L}}
     F_{r - 1} \]
  of $\varphi_A^{\ast} R$-modules implies that the $\varphi_A^{\ast} R$-module
  $G_r$ has $p$-complete $\tmop{Tor}$-amplitude
  \[ \nosymbol \leq \max \{ (p - 1 + \alpha - r) - 1, p - 1 + \alpha - (r - 1)
     \} = p + \alpha - r, \]
  and thus the $\varphi_A^{\ast} R$-module $\tmop{Fil}_H^{r - p} 
  \widehat{\tmop{dR}}_{R \otimes_{\mathbb{Z}}^{\mathbb{L}} \mathbb{F}_p /
  \overline{A} \otimes_{\mathbb{Z}}^{\mathbb{L}} \mathbb{F}_p}$ has
  $\tmop{Tor}$-amplitude$\nosymbol \leq \alpha - (r - p)$ for every $r \in
  \mathbb{Z}$. The result then follows.
\end{proof}

{\construction{\label{cons:null-htpy-lift-cofibs}Let $\mathcal{C}$ be a stable
$\infty$-category, and
\[ \begin{array}{ccc}
     A & \longrightarrow & B\\
     \longdownarrow &  & \longdownarrow\\
     C & \longrightarrow & D
   \end{array} \]
a commutative diagram in $\mathcal{C}$. Then a lift $C \rightarrow B$, namely
a map $C \rightarrow B$ along with homotopies as below,

\[\begin{tikzcd}[ampersand replacement=\&] A \& B \\ C \& D \arrow[from=1-1,
to=1-2] \arrow[from=1-1, to=2-1] \arrow[from=1-2, to=2-2] \arrow[from=2-1,
to=1-2] \arrow[from=2-1, to=2-2]
\end{tikzcd}\]

{\noindent}making two triangles $(A \rightarrow C \rightarrow B)$ and $(C
\rightarrow B \rightarrow D)$ commute, gives rise to a null-homotopy of the
canonical map
\[ \tmop{can} \of \tmop{cofib} (A \rightarrow C) \rightarrow \tmop{cofib} (B
   \rightarrow D) . \]
Indeed, such a lift splits the square into two triangles, which gives rise to
two ``octahedral'' fiber sequences
\[ \tmop{cofib} (A \rightarrow C) \longrightarrow \tmop{cofib} (A \rightarrow
   B) \longrightarrow \tmop{cofib} (C \rightarrow B) \]
and
\[ \tmop{cofib} (C \rightarrow B) \longrightarrow \tmop{cofib} (C \rightarrow
   D) \longrightarrow \tmop{cofib} (B \rightarrow D) . \]
The map $\tmop{can}$ is the composite of these two fiber sequences, which is
null-homotopic.}}

\appendix\section{Weak finite generation under countable generation
hypothesis}\label{app:weak-fg-countable-gen}

In this appendix, we show that, weak finite generation is a condition on
$\pi_0$ for countably generated connective modules
(\Cref{prop:countable-gen-weak-fg-pi-0}), as well as an easy characterization
of weak finite generation for countably generated modules
(\Cref{prop:countable-gen-crit-weak-fg}), which is comparable to relevant
results in {\cite[§6]{Fuchs1988}}. We also give counterexamples
(\Cref{ex:non-std-uniserial-mod,ex:hahn-val-field}) without countable
generation, which were discovered by ChatGPT-6 Astra.

\begin{proposition}
  \label{prop:countable-gen-weak-fg-pi-0}Let $R$ be a connective
  $\mathbb{E}_1$-ring, and $M$ a connective right $R$-module.
  \begin{enumerate}
    \item Suppose that the right $R$-module $M$ is weakly finitely generated
    (of rank$\nosymbol \leq r$). Then the right $\pi_0 (R)$-module $\pi_0 (M)$
    is weakly finitely generated.
    
    \item Conversely, suppose that the right $R$-module $M$ is countably
    generated (i.e. the right $\pi_0 (R)$-module $\pi_0 (M)$ is countably
    generated). If the right $\pi_0 (R)$-module $\pi_0 (M)$ is weakly finitely
    generated, then so is the right $R$-module $M$.
  \end{enumerate}
\end{proposition}

\begin{proof}
  Suppose that the right $R$-module $M$ is weakly finitely generated of
  rank$\nosymbol \leq r$. Then by definition, there exists a surjective map $P
  \rightarrow M$ of connective right $R$-modules with $\tmop{rk}_R (P) \leq
  r$. Taking $\pi_0$ and by \Cref{cor:rk-pi-0}, we see that the right $\pi_0
  (R)$-module $\pi_0 (M)$ is weakly finitely generated of rank$\nosymbol \leq
  r$.
  
  Conversely, suppose the right $\pi_0 (R)$-module $\pi_0 (M)$ is generated by
  a countable sequence $(x_i) \in \pi_0 (M)^{\mathbb{N}}$, and is weakly
  finitely generated of rank $r$. By definition, there exists a filtered
  system $(F_{\alpha})_{\alpha \in J}$ of free right $\pi_0 (R)$-modules of
  rank$\nosymbol \leq r$, along with a surjection $f \of \tmop{colim}_{\alpha}
  F_{\alpha} \twoheadrightarrow \pi_0 (M)$ of right $\pi_0 (R)$-modules.
  
  Then by dependent choice, we can pick a chain $(\alpha_0 \rightarrow
  \alpha_1 \rightarrow \cdots)$ in $J$, such that the image of the composite
  map $F_{\alpha_j} \rightarrow \tmop{colim}_{\alpha} F_{\alpha}
  \xrightarrow{f} \pi_0 (M)$ contains $\{ x_0, \ldots, x_j \}$ for every $j
  \in \mathbb{N}$. By choosing a base for each $F_{\alpha_j}$, we get a chain
  $(\tilde{F}_{\alpha_0} \rightarrow \tilde{F}_{\alpha_1} \rightarrow \cdots)$
  of free right $R$-modules of rank$\nosymbol \leq r$, along with an
  identification of the base-changed map $(\tilde{F}_{\alpha_j} \rightarrow
  \tilde{F}_{\alpha_{j + 1}}) \otimes_R^{\mathbb{L}} \pi_0 (R)$ with the map
  $F_{\alpha_j} \rightarrow F_{\alpha_{j + 1}}$ for every $j \in \mathbb{N}$,
  and compatible homotopies. This gives rise to a surjective map
  \[ \underset{j \in \mathbb{N}}{\tmop{colim}}  \tilde{F}_{\alpha_j}
     \nonconverted{longtwoheadrightarrow} M \]
  of connective right $R$-modules, whose source is a flat right $R$-module of
  rank$\nosymbol \leq r$.
\end{proof}

Countable generation hypothesis also allows us to deduce a criterion of weak
finite generation similar to \Cref{prop:crit-flat-rank-r}.

\begin{proposition}
  \label{prop:countable-gen-crit-weak-fg}Let $R$ be a connective
  $\mathbb{E}_1$-ring, and $M$ a connective right $R$-module.
  \begin{enumerate}
    \item Suppose that the right $R$-module $M$ is weakly finitely generated
    of rank$\nosymbol \leq r$. Then for every finite free right $R$-module
    $F$, every map $F \rightarrow M$ factors through some $R^s$ for $s \leq
    r$.
    
    \item Conversely, suppose that the right $R$-module $M$ is countably
    generated, and for every finite free right $R$-module $F$, every map $F
    \rightarrow M$ factors through some $R^s$ for $s \leq r$. Then the right
    $R$-module $M$ is weakly finitely generated of rank$\nosymbol \leq r$.
  \end{enumerate}
\end{proposition}

\begin{proof}
  The first statement follows from \Cref{prop:crit-flat-rank-r}. Now suppose
  that the right $R$-module $M$ is countably generated. We pick a sequence
  $(x_n)_{n \in \mathbb{N}_{> 0}}$ of elements of $\pi_0 (M)$ generating the
  $\pi_0 (R)$-module $\pi_0 (M)$. By dependent choice in $n \in \mathbb{N}$,
  we may pick sequences $(f_n \of R^r \rightarrow M)_{n \in \mathbb{N}}$,
  $(\theta_n \of R \rightarrow R^r)_{n \in \mathbb{N}_{> 0}}$ and $(g_n \of
  R^r \rightarrow R^r)_{n \in \mathbb{N}_{> 0}}$ of maps of right $R$-modules,
  such that $f_0 = 0$, and for $n \in \mathbb{N}_{> 0}$, the maps $f_n$,
  $\theta_n$, and $g_n$ of right $R$-modules are obtained by factoring $x_n
  \oplus f_{n - 1} \of R \oplus R^r \rightarrow M$ as
  \[ R \oplus R^r \xrightarrow{\theta_n + g_n} R^r \xrightarrow{f_n} M \]
  by assumption. We take $P \assign \tmop{colim} \left( R^r \xrightarrow{g_1}
  R^r \xrightarrow{g_2} R^r \xrightarrow{g_3} \cdots \right)$, and the
  previous data give rise to a map $f \of P \rightarrow M$ of right
  $R$-modules. As the map $x_n \of R \rightarrow M$ of right $R$-modules
  factors as $R \xrightarrow{\theta_n} R^r \xrightarrow{f_n} M$, we see that
  the map $f$ is surjective (on $\pi_0$), and the result follows.
\end{proof}

The countable generation in
\Cref{prop:countable-gen-weak-fg-pi-0,prop:countable-gen-crit-weak-fg} is
essential, and ChatGPT-6 Astra found counterexamples when the countable
generation is missing. We start with counterexample to the second statement of
\Cref{prop:countable-gen-crit-weak-fg} without countable generation
hypothesis.

Recall {\cite{Eklof1993}} that, letting $R$ be an integral domain with $K
\assign \tmop{Frac} (R)$, we say that an $R$-module is {\tmdfn{uniserial}} if
its submodules are totally ordered by inclusion, and such a module is
{\tmdfn{standard}} if it is isomorphic to a subquotient of the $R$-module $K$.

\begin{lemma}
  \label{lem:val-ring-uniserial-mod-weak-fg-rk-1}Let $V$ be a valuation ring
  with $K \assign \tmop{Frac} (V)$, and $M$ an $V$-module. Then it is standard
  uniserial if and only if it is weakly finitely generated of rank$\nosymbol
  \leq 1$.
\end{lemma}

\begin{proof}
  Every sub-$V$-module $N$ of $K$ is torsion-free thus flat, and $\tmop{rk}_V
  (N) \leq \tmop{rk}_V (K) = 1$ by \Cref{lem:rk-flat-Bezout}, every standard
  uniserial $V$-module is weakly finitely generated of rank$\nosymbol \leq 1$.
  
  Conversely, suppose that $M$ is weakly finitely generated of rank$\nosymbol
  \leq 1$. We pick a flat $V$-module $P$ of rank$\nosymbol \leq 1$ with a
  surjection $f \of P \rightarrow M$. By \Cref{lem:rk-flat-Bezout}, we have
  $\dim_K (P \otimes_V K) = \tmop{rk}_V (P) \leq 1$, and thus $P
  \hookrightarrow P \otimes_V K$ can be realized as a sub-$V$-module of $K$,
  and thus $P$ is standard uniserial, and so is its quotient $M$.
\end{proof}

\begin{example}
  \label{ex:non-std-uniserial-mod}By {\cite[Prop~8(i)]{Eklof1993}}, there
  exists a valuation ring $V$ with a non-standard uniserial $V$-module $M$.
  This, along with \Cref{lem:val-ring-uniserial-mod-weak-fg-rk-1}, shows the
  necessity of countable generation in \Cref{prop:countable-gen-crit-weak-fg}.
  Actually, the example $M$ there is also divisible, which even exclude the
  weak finite generation of $M$ by
  \Cref{prop:val-ring-div-uniserial-mod-weak-fg-std}.
\end{example}

\begin{proposition}
  \label{prop:val-ring-div-uniserial-mod-weak-fg-std}Let $V$ be a valuation
  ring with $K \assign \tmop{Frac} (V)$, and $M$ a uniserial $V$-module which
  is divisible. Suppose that the $V$-module $M$ is weakly finitely generated,
  then it is of rank$\nosymbol \leq 1$.
\end{proposition}

\begin{proof}
  We separate the proposition into two cases.
  \begin{enumerate}
    \item Suppose that the $V$-module $M$ contains an element $x \in M$ which
    is torsion-free. Then for every $y \in M$, since $M$ is uniserial, either
    $Vx \subseteq Vy$ or $Vy \subseteq Vx$. In both cases, we see that, if $y$
    is torsion, then $y = 0$. Consequently, the $V$-module $M$ is
    torsion-free, and $\dim_K (M \otimes_V K) = 1$, and thus $M$ can be
    embedded into the $V$-module $K$, thus it is standard uniserial, and we
    invoke \Cref{lem:val-ring-uniserial-mod-weak-fg-rk-1} to deduce.
    
    \item Now we suppose that $V$-module $M$ is torsion. Let $P$ be a flat
    $V$-module of rank $r > 1$ along with a surjection $f \of P \rightarrow M$
    of $V$-modules. We show that there exists a flat $V$-module along with a
    surjection $Q \rightarrow M$ such that $\tmop{rk}_V (Q) < r$. This implies
    what we want.
    
    Pick a non-zero element $x \in P$, and let $L \assign P \cap Kx$ which is
    contained in the $K$-vector space $P \otimes_V K$. Then $L$ is
    torsion-free thus flat, and it is of rank $1$ by
    \Cref{lem:rk-flat-Bezout}. Moreover, the kernel of the map $P \rightarrow
    P \otimes_V K \rightarrow (P \otimes_V K) / (Kx)$ is precisely $L$, which
    gives rise to an injection $P / L \hookrightarrow (P \otimes_V K) / (Kx)$
    into a $K$-vector space. It follows that $P / L$ is torsion-free thus
    flat, and it is of rank$\nosymbol < r$. We again separate into two cases.
    \begin{enumerate}
      \item Suppose that the restriction $\nobracket f |_L \of L \rightarrow
      M$ is surjective. Then we may take $Q \assign L$.
      
      \item There exists $y \in M$ which is not contained in $f (L)$. Then by
      uniseriality of $M$, we have $f (L) \subseteq Vy$. As the $V$-module $M$
      is torsion, there exists a non-zero element $\alpha \in V$ such that
      $\alpha y = 0$, and thus $\nobracket \alpha f |_L = 0$. It follows that
      the map $\alpha f \of P \rightarrow M$ of $V$-modules induces a map $g
      \of P / L \rightarrow M$ of $V$-modules. Since the $V$-module $M$ is
      divisible, the map $\alpha f \of P \rightarrow M$ is surjective, thus so
      is the induced map $g$. Thus we may take $Q \assign P / L$.
    \end{enumerate}
  \end{enumerate}
\end{proof}

Now we describe the counterexample to the second statement of
\Cref{prop:countable-gen-weak-fg-pi-0}. Recall that, letting $R$ be an
integral domain with $K \assign \tmop{Frac} (R)$, an $R$-module $M$ is
{\tmdfn{Matlis cotorsion}} if $\tmop{Ext}_R^1 (K, M) = 0$, it is
{\tmdfn{$h$-reduced}} if $\tmop{Hom}_R (K, M) = 0$ (this is the case if $M$ is
reduced, as explained in {\cite{Bazzoni2020}}), and a torsion-free $R$-module
$M$ is {\tmdfn{$R$-complete}} if the canonical map
\[ M \longrightarrow \lim_{s \in R \setminus 0}  (M / sM) \]
is an isomorphism on $\pi_0$. The following statement is mentioned in
{\cite[§1]{Fuchs2007}}, and we include a proof for convenience of readers.

\begin{lemma}
  \label{lem:crit-R-cpl}Let $R$ be an integral domain with $K \assign
  \tmop{Frac} (R)$. Then a torsion-free $R$-module $M$ is $R$-complete if and
  only if the following two conditions are satisfied.
  \begin{enumerate}
    \item The $R$-module $M$ is $h$-reduced.
    
    \item The $R$-module $M$ is Matlis cotorsion.
  \end{enumerate}
\end{lemma}

\begin{proof}
  We rewrite the $R$-module $K$ as the filtered colimit $\tmop{colim}_{s \in R
  \setminus 0} Rs^{- 1}$. Thus
  \begin{eqnarray*}
    \tmop{RHom}_R (K, M) & = & \lim_{s \in R \setminus 0} \tmop{RHom}_R (Rs^{-
    1}, M)\\
    & = & \lim_{s \in R \setminus 0} sM
  \end{eqnarray*}
  where the limit is derived, and consequently, we get a fiber sequence
  \[ \tmop{RHom}_R (K, M) \longrightarrow M \longrightarrow \lim_{s \in R
     \setminus 0}  (M / sM) \]
  in $D (R)$, which gives rise to a long exact sequence
  \[ 0 \longrightarrow \tmop{Hom}_R (K, M) \longrightarrow M \longrightarrow
     \pi_0 (\lim_{s \in R \setminus 0}  (M / sM)) \longrightarrow
     \tmop{Ext}_R^1 (K, M) \longrightarrow 0. \]
  The result then follows.
\end{proof}

\begin{lemma}
  \label{lem:Pruefer-flat-fin-rk-fin-fil}Let $R$ be a Bézout domain, and $P$
  a flat $R$-module of rank$\nosymbol \leq r$. Then there exists a finite
  filtration on $P$ whose associated graded pieces are flat $R$-modules of
  rank$\nosymbol \leq 1$.
\end{lemma}

\begin{proof}
  Let $K \assign \tmop{Frac} (R)$, and $P_K \assign P \otimes_R K$. As $P$ is
  a filtered colimit of finite free $R$-modules of rank$\nosymbol \leq r$, the
  $P_K$ is a filtered colimit of $K$-vector spaces of dimension$\nosymbol \leq
  r$, thus $\dim_K (P_K) \leq r$. It follows that $P_K$ admits a finite
  filtration
  \[ 0 = F_0 \subseteq F_1 \subseteq \cdots \subseteq F_r = P_K \]
  of $K$-vector spaces with $\dim_K (F_j / F_{j - 1}) \leq 1$ for $j = 1,
  \ldots, r$. We take $H_j \assign P \cap F_j$. Then the induced map $H_j /
  H_{j - 1} \rightarrow F_j / F_{j - 1}$ is injective, and thus the $R$-module
  $H_j / H_{j - 1}$ is torsion-free, therefore flat as $R$ is a Prüfer
  domain.
  
  Now we show that $\tmop{rk}_R (H_j / H_{j - 1}) \leq 1$. If $F_j = F_{j -
  1}$, then $H_j = H_{j - 1}$ and the result is true. Otherwise, $\dim_K (F_j
  / F_{j - 1}) = 1$, and we get an injection $H_j / H_{j - 1} \hookrightarrow
  K$ of $R$-modules. Now for every finite collection $\{ x_1 / y_1, x_2 / y_2,
  \ldots, x_m / y_m \} \subseteq H_j / H_{j - 1}$, as $R$ is Bézout, the
  sub-$R$-module generated by them is principal, generated by
  \[ (y_1 \cdots y_m)^{- 1} \gcd (x_1 y_2 \cdots y_m, y_1 x_2 \cdots y_m,
     \ldots, y_1 y_2 \cdots x_m) . \]
  It follows that $\tmop{rk}_R (H_j / H_{j - 1}) \leq 1$.
\end{proof}

\begin{remark}
  \Cref{lem:Pruefer-flat-fin-rk-fin-fil} gives rise to another proof of
  \Cref{lem:rk-flat-Bezout}.
\end{remark}

\begin{corollary}
  \label{cor:val-ring-flat-fin-rk-Matlis-cotor}Let $V$ be a valuation ring
  which is complete with respect to the valuation, and $K \assign \tmop{Frac}
  (V)$. Then for every flat $V$-module $P$ of finite rank, we have
  $\tmop{Ext}_V^1 (K, P) = 0$.
\end{corollary}

\begin{proof}
  The statement is true for $P = 0$. Suppose that $P \neq 0$ and $\tmop{rk}_V
  (P) = 1$. Then the $V$-module $P$ embeds into $P \otimes_V K$, which is a
  $K$-vector space of dimension $1$ as in the proof of
  \Cref{lem:Pruefer-flat-fin-rk-fin-fil}, thus we may view $P$ as a
  sub-$V$-module of $K$. We separate two cases:
  \begin{enumerate}
    \item If $P = K$, then as the $V$-module $K$ is injective, we have
    $\tmop{Ext}_V^1 (K, P) = 0$.
    
    \item Suppose that $P \neq K$. We may rescale the module $P$ so that there
    exists $d \in V \setminus 0$ such that $Vd \subseteq P \subseteq V$. This
    implies that the ``$V$-topology'' on $P$ is the same as the subspace
    topology of $P$ in $V$, and that $P$ is clopen in $V$, thus it is complete
    with respect to the topology as well, and therefore $V$-complete. It
    follows that $\tmop{Ext}_V^1 (K, P) = 0$ by \Cref{lem:crit-R-cpl}.
  \end{enumerate}
  The general case then follows from \Cref{lem:Pruefer-flat-fin-rk-fin-fil}.
\end{proof}

\begin{corollary}
  \label{cor:cpl-val-ring-weak-fin-gen-surj-quot-field-sect}Let $V$ be a
  valuation ring which is complete with respect to the valuation, and $K
  \assign \tmop{Frac} (V)$. Let $M$ be a weakly finitely generated connective
  $V$-module along with a map $\varphi \of M \rightarrow K$ in $D (V)$ which
  is surjective on $\pi_0$. Then the map $M \rightarrow K$ admits a section in
  $D (V)$.
\end{corollary}

\begin{proof}
  Let $P$ be a flat $V$-module of finite rank along with a surjective map $P
  \rightarrow M$ of connective $V$-modules. Then the composite map $f \of P
  \rightarrow M \rightarrow K$ is a surjective map of $V$-modules concentrated
  in degree $0$. Since the $V$-module $\tmop{Ker} (f) \subseteq P$ is a
  torsion-free thus flat, and $\dim_K (\tmop{Ker} (f) \otimes_V K) \leq \dim_K
  (P \otimes_V K) < \infty$, by \Cref{lem:rk-flat-Bezout}, the $V$-module
  $\tmop{Ker} (f)$ is also flat of finite rank, thus by
  \Cref{cor:val-ring-flat-fin-rk-Matlis-cotor}, we have $\tmop{Ext}_V^1 (K,
  \tmop{Ker} (f)) = 0$, thus the surjection $f$ admits a section, and so does
  the surjective map $\varphi$.
\end{proof}

\begin{example}[ChatGPT-6 Astra]
  \label{ex:hahn-val-field}Let $\Gamma \assign \bigoplus_{\alpha < \omega_1}
  \mathbb{Z}e_{\alpha}$ be the totally ordered abelian group, where $\sum
  a_{\alpha} e_{\alpha} < \sum b_{\alpha} e_{\alpha}$ if there exists $\beta <
  \omega_1$ such that, for every $\beta > \alpha$, we have $a_{\alpha} =
  b_{\alpha}$, and $a_{\beta} < b_{\beta}$. Let $k$ be a field. We consider
  the Hahn's valued field $K \assign k \pparen{t^{\Gamma}}$, which is complete
  with respect to its valuation $v$, with $V \assign \{ x \in K \barsuchthat v
  (x) \geq 0 \}$. Then by {\cite[Thm~A \& the last paragraph of
  §2]{Osofsky1967}}, and that $K$ is not $\aleph_0$-generated (as $\omega_1$
  is $\aleph_1$-filtered) but $\aleph_1$-generated, the projective dimension
  $\tmop{pd}_V (K)$ is $2$.
  
  Now we pick a connective $V$-module $M$ with $\pi_0 (M) = K$ with projective
  dimension $\tmop{pd}_V (M) \leq 1$ (e.g. let $F_1 \rightarrow F_0
  \rightarrow K$ be a free presentation of $K$, and take $M \assign
  \tmop{cofib} (F_1 \rightarrow F_0)$). If $M$ is weakly finitely generated,
  then by \Cref{cor:cpl-val-ring-weak-fin-gen-surj-quot-field-sect}, the
  Postnikov truncation map $M \rightarrow K$ admits a section in $D (V)$,
  which is impossible as $\tmop{pd}_V (M) < \tmop{pd}_V (K)$. Thus the
  connective $V$-module $M$ is not weakly finitely generated, while the
  $V$-module $\pi_0 (M) = K$ is weakly finitely generated of rank$\nosymbol
  \leq 1$.
\end{example}

\

\end{document}